\documentclass[12pt]{iopart}
\usepackage[dcucite]{harvard}
\usepackage{graphicx}
\usepackage{float}
\usepackage{iopams}

\usepackage{subfigure}

\usepackage[normalem]{ulem}
\usepackage[usenames,dvipsnames]{color}

\begin{document}
\title[Iterative image reconstruction in 3D OAT]{Investigation of 
iterative image reconstruction in three-dimensional 
optoacoustic tomography}

\author{Kun Wang$^1$, Richard Su$^2$, Alexander A Oraevsky$^2$ \\
and Mark A Anastasio$^1$}
\address{
$^1$ Department of Biomedical Engineering, Washington University in St. Louis,
St. Louis, MO 63130}
\address{
$^2$ TomoWave Laboratories, 675 Bering drive, Suite 575, Houston,
Texas, Houston, TX 77057}

\ead{\mailto{anastasio@wustl.edu}}

\begin{abstract}
Iterative image reconstruction algorithms for optoacoustic tomography (OAT),
also known as photoacoustic tomography,
have the ability to improve image quality over analytic algorithms
due to their ability to incorporate accurate models of the imaging physics, instrument
response, and measurement noise.
However, to date,  there have been few reported attempts to employ advanced
 iterative image reconstruction algorithms
for improving image quality in three-dimensional (3D) OAT.
In this work, we implement and investigate two iterative
 image reconstruction methods for use with a 3D OAT small animal imager: 
namely, a penalized least-squares (PLS) method employing a quadratic smoothness penalty 
and a PLS method employing a total variation norm penalty.
The reconstruction algorithms employ accurate  models of the ultrasonic transducer
impulse responses. 
Experimental data sets are employed to compare the performances of the iterative reconstruction 
algorithms to that of a 3D filtered backprojection (FBP) algorithm.
By use of quantitative measures of image quality,
we demonstrate that the iterative reconstruction algorithms
can mitigate image artifacts and preserve spatial resolution more effectively than
FBP algorithms. 
These features suggest that the use of advanced image reconstruction
algorithms can improve the effectiveness of 3D OAT 
while reducing the amount of data required for biomedical applications.
\end{abstract}
\submitto {Physics in Medicine and Biology}
\maketitle


\section{INTRODUCTION}
\label{sec:intro}
Optoacoustic tomography (OAT), also known as photoacoustic tomography, 
is a rapidly emerging imaging modality that has
great potential for a wide range of biomedical imaging applications 
\cite{OraBook,WangTutorial,KrugerMP1999,CoxOSA2006}.
OAT is a hybrid imaging method in which biological tissues are illuminated with short laser pulses, 
which results in the generation of internal acoustic wavefields via the thermoacoustic effect.
The initial amplitudes of the induced acoustic wavefields are proportional to the
spatially variant absorbed optical energy density in the tissue.
The propagated acoustic wavefields are detected by use of a collection of wide-band
ultrasonic transducers that are located outside the object.
From knowledge of these acoustic data, an image reconstruction
algorithm is employed to estimate the absorbed optical energy density within the tissue.

A variety of analytic image reconstruction algorithms for three-dimensional (3D) OAT have been developed
\cite{Kunyansky:07,Finch:02,Xu:2005bp,Xu:planar}.
These algorithms are of filtered backprojection forms
and assume that the underlying model that relates
the object function to measured data is a 
spherical Radon transform.
Analytic image reconstruction algorithms
generally possess several limitations that impair their performance.
For example, analytic algorithms are often based on discretization
of a continuous reconstruction formula and
require the measured data to be densely sampled on
an aperture that encloses the object. This is problematic for 3D
OAT, in which acquiring densely sampled acoustic measurements on a
two-dimensional (2D) surface can require expensive transducer arrays
and/or long data-acquisition times if a mechanical scanning is employed.
Moreover, the simplified forward models, such as the spherical
Radon transform, upon which analytic image reconstruction algorithms are based, 
do not comprehensively describe the imaging physics
or response of the detection system \cite{TMI:transmodel}.
Finally, analytic methods ignore measurement noise and will
generally yield images that have suboptimal trade-offs between
image variances and spatial resolution.  The use of iterative
image reconstruction algorithms \cite{Fessler:94,AnastasioTATHT,Wernickbookchap21,XPanIP2009}  
can circumvent all of these shortcomings.

When coupled with suitable OAT imager designs, iterative image reconstruction algorithms can
improve image quality and permit reductions in data-acquistion times as compared with 
those yielded by use of analytic reconstruction algorithms.
Because of this,
the development and investigation of iterative image reconstruction
algorithms for OAT \cite{GPaltauf:2002} is an important research topic of current interest.  
Recent studies have sought to develop improved discrete imaging models
\cite{HJiang:2007,JCarson:2008,Ntziachristos:2011,TMI:transmodel}
as well as  advanced reconstruction algorithms
\cite{OATTV09:Provost,GuoCS:2010,KunSPIETV:2011}. 
The majority of these studies utilize approximate 2D imaging models
and measurement geometries in which focused transducers are employed to suppress out-of-plane
acoustic signals and/or a thin object embedded in an acoustically homogeneous
background is employed.
Because image reconstruction of extended objects in OAT is inherently a 3D problem, 
 2D image reconstruction approaches may not yield
accurate values of the absorbed optical energy density 
 even when the measurement data are densely sampled.
This is due to the fact that  simplified 2D imaging models cannot
accurately describe transducer focusing and out-of-plane
acoustic scattering effects;  this results in inconsistencies between
the imaging model and the measured data that can result in artifacts
and loss of accuracy in the reconstructed images.

Several 3D OAT imaging systems have been constructed
and investigated \cite{KrugerMP2010,JCarson:2008,petermice:jbo}.
These systems employ unfocused ultrasonic transducers 
and analytic 3D image reconstruction algorithms.
Only limited studies of the use of iterative 3D 
algorithms for reconstructing extended objects have been conducted; and these 
studies employed simple phantom objects \cite{GPaltauf:2002,KunSPIETV:2011,KunSPIETV:2012,JCarson:2008}. 
Moreover, iterative image reconstruction in 3D OAT can be extremely computationally
burdensome, which can require the use of high performance computing platforms.
Graphics processing units (GPUs) can now be employed to accelerate 3D iterative image reconstruction
algorithms to the point where they are feasible.
However, there remains an important need for the development of accurate discrete imaging
models and image reconstruction algorithms for 3D OAT  and an investigation
of their ability to mitigate different types of measurement errors found in 
real-world implementations. 

In this work, we implement and investigate two 3D iterative
 image reconstruction methods for use with a small animal OAT imager.
Both reconstruction algorithms compensate for the ultrasonic transducer
responses but employ different regularization strategies.  
We compare the different regularization strategies by use of quantitative measures of 
image quality. 
Unlike previous studies, we apply the 3D image reconstruction algorithms not only
to experimental phantom data but also to the data from a mouse whole-body imaging experiment.
The aim of this study is to demonstrate the feasibility and efficacy 
of iterative image reconstruction in 3D OAT and to identify current limitations 
in its performance. 

The remainder of the article is organized as follows: 
In \Sref{sec:imagingmodel}, we derive the numerical imaging model that is employed 
by the iterative image reconstruction algorithms
and briefly review the three image reconstruction algorithms. 
\Sref{sec:methods} describes the experimental studies including
the data acquisition, implementation details, 
and approaches for image quality assessment. 
The numerical results are presented in \Sref{sec:results}, 
and a discussion of our findings is presented in \Sref{sec:summary}. 

\section{Background: imaging models and reconstruction algorithms for 3D OAT}
\label{sec:imagingmodel}
Iterative image reconstruction algorithms commonly employ a discrete imaging model that relates the
measured data to an estimate of the sought-after object function.
We previously proposed a general procedure for constructing discrete OAT imaging models
that incorporate the spatial and acousto-electric 
impulse responses of an ultrasonic transducer \cite{TMI:transmodel}.
We review the salient features of this procedure in \Sref{sec:ima}.   For use in the studies
presented in this work, in \Sref{sec:imb} we reformulate 
the discrete imaging model in the temporal-frequency space
for the case of flat rectangular ultrasonic transducers.

\subsection{Discrete imaging model in the time-domain}
\label{sec:ima}
A canonical OAT imaging model in its continuous form is expressed as
\cite{LVbook:optics,OraBook,TMI:transmodel}:
\begin{equation}\label{eqn:ccmodel}
 p(\mathbf r',t) = {\frac{\beta}{4\pi C_p} } 
   \int_V\!\! d^3\mathbf{r}\, A(\mathbf{r}) 
    {\frac{d}{dt}} 
    \frac{\delta\left(t-{
         \frac{\vert \mathbf r'-\mathbf{r}\vert}{ c_0}
       }\right)}
  {\vert \mathbf r'-\mathbf r\vert} *_t I(t), 
\end{equation}
where $p(\mathbf r',t)$ denotes the acoustic pressure
measured at location $\mathbf r'$ and time $t$,
$A(\mathbf{r})$ denotes the sought-after absorbed optical energy density,
$I(t)$ describes the normalized temporal profile of the illumination pulse,
$\delta (t)$ is the Dirac delta function,
$V$ denotes the object's support volume, 
$*_t$ denotes 1D temporal convolution, and
$\beta$, $c_0$, and $C_p$ denote the thermal coefficient
of volume expansion, (constant) speed-of-sound,
and the specific heat capacity of the medium at constant pressure, respectively. 
Because many OAT applications employ a laser pulse of nano-seconds in duration, we assume $I(t)\approx\delta(t)$ in this study. 
In accordance, we drop the last temporal convolution in \eref{eqn:ccmodel} hereafter.
This model assumes an idealized data-acquisition process and neglects
finite sampling effects.

In practice, the acoustic pressure is converted to a voltage  signal 
by use of ultrasonic transducers that is subsequently sampled and recorded.
Consider that the ultrasonic transducers collect data at $Q$  locations
that are specified by the index  $q=0,\cdots,Q-1$ and 
$K$ temporal samples, specified by the index $k=0,\cdots, K-1$, are acquired at each  location with a sampling
 interval $\Delta T$. 
A continuous-to-discrete (C-D) imaging model \cite{BarrettBook,MarkBookChapter} for OAT can be generally expressed as 
\cite{TMI:transmodel}: 
\begin{equation}\label{eqn:cdmodel}
  [\mathbf u]_{qK+k} = u_q(t)\Big|_{t=k\Delta T}=
        h^e(t)*_t \frac{1}{\Omega_q}
    \int_{\Omega_q}\!\! d^2\mathbf r'\,
        p(\mathbf r',t)\Big|_{t=k\Delta T},
\end{equation}
where 
$u_q(t)$ is the pre-sampled electric voltage signal corresponding to
location index $q$,
the surface integral is over the detecting area of the $q$-th
transducer denoted by $\Omega_q$, 
and 
$h^e(t)$ denotes the acousto-electric impulse response (EIR) of transducers. 
The  $QK\times 1$ vector  $\mathbf u$ denotes a lexiographically ordered version
of the sampled data.
The notation $[\mathbf u]_{qK+k}$ is employed to denote the $(qK+k)$-$th$
element of  $\mathbf u$.
The pressure data function $p(\mathbf r',t)$ is determined
by $A(\mathbf r)$  via \eref{eqn:ccmodel}.
Accordingly,  the C-D mapping given by \eref{eqn:cdmodel} maps the function $A(\mathbf r)$
to the measurement vector $\mathbf u$.

To obtain a discrete-to-discrete (D-D) imaging model for use with
iterative image reconstruction algorithms,
a finite-dimensional approximate representation of the object function $A(\mathbf r)$
can be introduced. 
We have previously employed \cite{TMI:transmodel} the representation 
\begin{equation}\label{eqn:sph_vox}
  A^a(\mathbf r)=\sum_{n=0}^{N-1} 
   \big[{\boldsymbol \theta}\big]_n\phi_n({\mathbf{r}}),
\end{equation}
where the superscript $a$ indicates that $A^a(\mathbf r)$ is an approximation
of $A(\mathbf r)$ and
 $\{\phi_n(\mathbf r)\}_{n=0}^{N-1}$ are expansion functions defined as
\begin{equation}\label{eqn:basisfunc}
  \phi_n(\mathbf r) = \left\{\begin{array}{ll}
                1, & {\rm if}\quad|\mathbf r - \mathbf r_n| \leq \epsilon \\
                0, & {\rm otherwise} 
                \end{array}\right. .
\end{equation}
Here, $\mathbf r_n=(x_{n},y_{n},z_n)^T$ specifies the coordinate of the
$n$-th grid point of a uniform Cartesian lattice
and $\epsilon$  is the half spacing between lattice points.
 The $n$-th component
of the coefficient vector $\boldsymbol \theta$ is defined as
\begin{equation}
  \big[\boldsymbol \theta\big]_n = \frac{V_{\rm cube}}{V_{\rm sph}}
              \int_{V}\!\! d^3\mathbf{r} \;
               \phi_n (\mathbf{r})A(\mathbf{r}), 
\end{equation}
where $V_{\rm cube}$ and $V_{\rm sph}$ are the volumes of
a cubic voxel of dimension $2\epsilon$ and
of a spherical voxel of radius $\epsilon$ respectively.

Let $u_q^a(t)$ denote the pre-sampled voltage signal that would
be produced by the absorbed optical energy density $ A^a(\mathbf r)$.
Note that  $u_q^a(t)$  is an approximation of  $u_q(t)$, 
which would be produced by $ A(\mathbf r)$.
By use of \eref{eqn:ccmodel}-\eref{eqn:sph_vox}, it can be verified that
\begin{equation}\label{eqn:ua2}
\eqalign{
  \fl  u_q^a(t) &= 
    \underbrace{-\frac{\beta c^3 \pi}{C_p}t\Big[  H(t+\frac{\epsilon}{c_0})
            -H(t-\frac{\epsilon}{c_0}) \Big]}_{ \equiv p_0(t)} *_t
            h^e(t) *_t
            \frac{1}{\Omega_q}\sum_{n=0}^{N-1} 
     \big[{\boldsymbol \theta}\big]_n
       \underbrace{\int_{\Omega_q}\!\!d^2\mathbf r'\, 
          \frac{\delta(t-\frac{|\mathbf r' -\mathbf r_n|}{c_0})}
                 {2\pi|\mathbf r' -\mathbf r_n|}}_{\equiv h^s_q(\mathbf r_n, t)}\\
  \fl &= p_0(t) *_t h^e(t) *_t
        \frac{1}{\Omega_q}   \sum_{n=0}^{N-1} 
     \big[{\boldsymbol \theta}\big]_n h^s_q(\mathbf r_n, t)},
\end{equation}
where 
$H(t)$ is Heaviside step function, $p_0(t)$ is the `N'-shape pressure profile  produced
by a uniform sphere of radius $\epsilon$,
and $h^s_q(\mathbf r_n,t)$  denotes the spatial impulse response (SIR) of
the $q$-th transducer. 
By temporally sampling \eref{eqn:ua2} and employing
the approximation $[\mathbf u]_{qK+k}\approx u_q^a(t)\big|_{t=k\Delta T}$, 
one can establish \cite{TMI:transmodel} a D-D imaging model as
\begin{equation}
\mathbf u = \mathbf H_t \boldsymbol\theta,
\label{eqn:syst}
\end{equation}
where the system matrix $\mathbf H_t$ maps the coefficient vector $\boldsymbol \theta$, which
determines $A^a(\mathbf r)$, to the measured temporal samples of the voltage signals.

\subsection{Temporal frequency-domain version of the discrete imaging model}
\label{sec:imb}
Because a transducer's EIR $h^e(t)$ must typically be measured, it generally cannot 
be described by a simple analytic expression.
Accordingly, the two temporal convolutions in \eref{eqn:ua2}
must be approximated by use of discrete time convolution operations.
However, both  $p_0(t)$ and $h_q^s(\mathbf r_n, t)$ are very narrow
functions of time, and therefore temporal sampling can result in strong aliasing 
artifacts unless very large sampling rates are employed. 
As described below, to circumvent this we reformulated the D-D imaging model in the 
temporal frequency domain.

\begin{figure}[h]
\centering
\includegraphics[width=8.0cm]{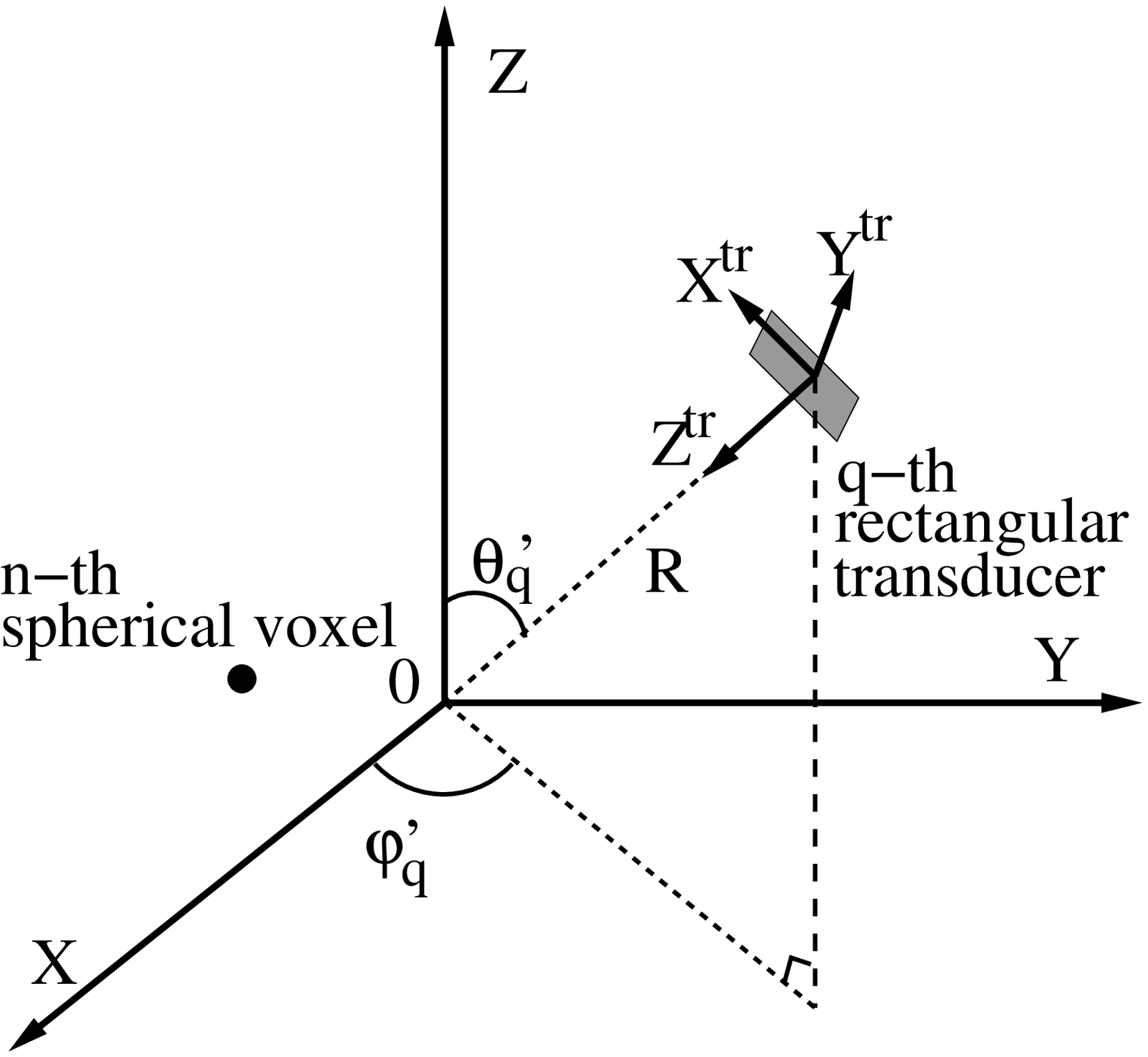}
\caption{Schematic of the local coordinate system for the $q$-th 
transducer where the  
$z^{\rm tr}$-axis points to the origin of the global coordinate system,
the $x^{\rm tr}$ and $y^{\rm tr}$-axes are along the two edges of the rectangular 
transducer respectivley. 
\label{fig:loccoordsys}
}
\end{figure}

Consider  \eref{eqn:ua2} expressed in the temporal frequency domain:
\begin{equation}\label{eqn:uaf}
  \tilde u^a_q(f)
                   = \tilde p_0(f)
                    \tilde h^e(f) 
               \frac{1}{\Omega_q}
                    \sum_{n=0}^{N-1} 
     \big[{\boldsymbol \theta}\big]_n \tilde h^s_q(\mathbf r_n, f),
\end{equation}
where $f$ is the temporal frequency variable conjugate to $t$
 and a `$\;\tilde{\;} \;$' above a function denotes the Fourier transform of that function defined as:
\begin{equation}
  \tilde x(f) = \int_{-\infty}^\infty\!\!dt\;
    x(t) \exp(-\hat\jmath2\pi f t). 
\end{equation}
For $f\neq 0$, the temporal Fourier transform of $p_0(t)$ is given by
\begin{equation}
  \tilde p_0(f) = -\hat\jmath \frac{\beta c_0^3}{C_pf}
    \Bigg[
    \frac{\epsilon}{c_0}\cos\big(\frac{2\pi f\epsilon}{c_0}\big)
   -\frac{1}{2\pi f}\sin\big(\frac{2\pi f\epsilon}{c_0}\big)
    \Bigg].
\end{equation}
When the transducer has a flat and rectangular detecting surface 
of area $a\times b$, under the far-field assumption, the temporal
Fourier transform of 
the SIR is given by \cite{PStepanishen:1971}: 
\begin{equation}\label{eqn:sir}
\fl\quad
  \tilde h_q^s(\mathbf r_n, f)
 = \frac{ab \exp(-\hat{\jmath}\,2\pi f\frac{|\mathbf r'_q-\mathbf r_n|}{c_0})}
      {2\pi |\mathbf r'_q-\mathbf r_n|}
  {\rm sinc}\Big(\pi f
    \frac{ax^{\rm tr}_{nq}}{c_0|\mathbf r'_q-\mathbf r_n|}\Big)
  {\rm sinc}\Big(\pi f
    \frac{by^{\rm tr}_{nq}}{c_0|\mathbf r'_q-\mathbf r_n|}\Big),
\end{equation}
where $x^{\rm tr}_{nq}$ and $y^{\rm tr}_{nq}$ specify the 
transverse coordinates in a local coordinate system that is centered at the 
$q$-th transducer, as depicted in \Fref{fig:loccoordsys},
corresponding to the location  of a point source described by a 3D Dirac delta function.
The SIR does not depend on the third coordinate ($z^{\rm tr}$)
specifying the point-source location due to the far-field assumption.
Given the voxel location $\mathbf r_n=(x_n,y_n,z_n)$
and transducer location $\mathbf r'_q=(r'_q, \theta'_q,\phi'_q)$,
expressed in spherical coordinates as  shown in \Fref{fig:loccoordsys},
the corresponding values of the local coordinates can be computed as:
\begin{eqnarray}
    &x^{\rm tr}_{nq}=-x_n\cos\theta'_q\cos\phi'_q
        -y_n\cos\theta'_q\sin\phi'_q+z_n\sin\theta'_q,\\
    &y^{\rm tr}_{nq}=-x_n\sin\phi'_q+y_n\cos\phi'_q.
\end{eqnarray}

\Eref{eqn:uaf} can be 
discretized by considering $L$ temporal frequency samples specified by a sampling 
interval $\Delta f$ that are referenced by the index $l=0,1,\cdots,L-1$. 
Utilizing the approximation 
$[\tilde{\mathbf u}]_{qL+l}\approx \tilde u_q^a(f)\big|_{f=l\Delta T}$
yields the D-D imaging model:
\begin{equation}\label{eqn:imagmodeldd}
  \tilde {\mathbf u} = \mathbf H \boldsymbol \theta, 
\end{equation}
where $\mathbf H$ is the system matrix of dimension $QL\times N$,
whose elements are defined by
\begin{equation}\label{eqn:SysMatrix}
  \big[\mathbf H\big]_{qL+l,n} = 
              \tilde p_0(f)
              \tilde h^e(f)  
               \frac{\tilde h^s_q(\mathbf r_n, f)}{ab} 
          \Big\vert_{f=l\Delta f}. 
\end{equation}
The imaging model in \eref{eqn:imagmodeldd} will form the
basis for the iterative image reconstruction studies described
in the remainder of the article.

\subsection{Reconstruction algorithms}
We investigated a 3D filtered backprojection algorithm (FBP) and 
two iterative reconstruction algorithms that employed 
different forms of regularization. 

{\bf Filtered backprojection:}
A variety of FBP type OAT image reconstruction algorithms have been developed
based on the continuous imaging model described by \eref{eqn:ccmodel}
\cite{Kunyansky:07,Finch:02,Xu:2005bp,Xu:planar}. 
If sampling effects are not considered and a closed measurement surface 
is employed,  these algorithms provide a 
mathematically exact mapping from the acoustic pressure function 
$p(\mathbf r',t)$ to the obsorbed energy density function $A(\mathbf r)$. 
Since we only have direct access to electric signals in practice, 
in order to apply FBP algorithms, we need to first estimate the sampled values of the 
acoustic pressure data from the measured electric signals. 
In this study, we considered a spherical scanning geometry. 
When the object is near the center of the measurement sphere, 
the surface integral over $\Omega_q$ in \eref{eqn:cdmodel},
i.e., SIR effect, is negligible. 
The remaining EIR effect is described by a temporal convolution. 
We employed linear regularized Fourier deconvolution \cite{KrugerMP1999}
to estimate the pressure data, expressed in temporal frequency domain as: 
\begin{equation}
  \tilde p(\mathbf r', f) = \frac{\tilde u(\mathbf r', f)}{\tilde h^e(f)}
    \tilde W(f),
\end{equation}
where $\tilde W(f)$ is a window function for noise suppression. 
In this study, we adopted the Hann window function defined as: 
\begin{equation}
  \tilde W(f) = \frac{1}{2}\Big[1-\cos(\pi\frac{f_c-f}{f_c})\Big],
\end{equation}
where $f_c$ is the cutoff frequency. 
We implemented the following FBP reconstruction formula for a spherical measurement 
geometry \cite{Finch:02}:
\begin{equation}\label{eqn:fbp}
  A(\mathbf r) = -\frac{C_p}{2\pi \beta c_0^2R_s}
   \int_S\!\! d^2\mathbf r'\, 
    \frac{2p(\mathbf r', t)+t\frac{\partial}{\partial t} p(\mathbf r', t)
           }{|\mathbf r-\mathbf r'|}\Bigg|
      _{t=\frac{|\mathbf r-\mathbf r'|}{c_0}},
\end{equation}
where $R_s$ and $S$ denote the radius and surface area of the measurement sphere
respectively. 
Note that the value of the cutoff frequency $f_c$ controls the degree of 
noise suppression during the deconvolution, 
thus indirectly regularizing the FBP algorithm. 

{\bf Penalized least-squares with quadratic penalty:}
PLS reconstruction methods seek to minimize a cost-function of the form as: 
\begin{equation}\label{eqn:PLSObjFunc}
  \hat {\boldsymbol \theta} = \arg\min_{\boldsymbol \theta}
    \Vert \tilde{\mathbf u} - \mathbf H\boldsymbol\theta\Vert^2
    +\alpha R(\boldsymbol\theta),
\end{equation}
where $R(\boldsymbol\theta)$ is a regularizing penalty term whose impact
is controlled by the regularization parameter $\alpha$.
We employed the conventional quadratic smoothness 
Laplacian penalty given by \cite{Fessler:94}: 
\begin{equation}\label{eqn:penalty}
\fl  R(\boldsymbol \theta) = 
  \sum_{n=0}^{N-1}
   \Big(2[\boldsymbol \theta]_n-[\boldsymbol \theta]_{k_{x_1}}
                               -[\boldsymbol \theta]_{k_{x_2}}\Big)^2
  +\Big(2[\boldsymbol \theta]_n-[\boldsymbol \theta]_{k_{y_1}}
                               -[\boldsymbol \theta]_{k_{y_2}}\Big)^2
  +\Big(2[\boldsymbol \theta]_n-[\boldsymbol \theta]_{k_{z_1}}
                               -[\boldsymbol \theta]_{k_{z_2}}\Big)^2,
\end{equation}
where $k_{x_1}$ and $k_{x_2}$ were the indices of the two neighbor 
voxels before and after the $n$-th voxel along x-axis. 
Similarly, $k_{y_1}$, $k_{y_2}$ and  $k_{z_1}$, $k_{z_2}$ 
were the indices of the neighbor voxels along y- and z-axis respectively. 
The reconstruction algorithm for solving \eref{eqn:PLSObjFunc}
was based on the  Fletcher Reeves version of conjugate gradient (CG) method
\cite{Wernickbookchap21}, and will be referred to as the PLS-Q algorithm. 

{\bf Penalized least-squares with total variation norm penalty:}
We also investigated the PLS algorithm regularized by a TV-norm penalty. 
This method seeks to minimize a cost-function of the form as: 
\begin{equation}\label{eqn:TVcost}
  \hat {\boldsymbol \theta} = \arg\min_{\boldsymbol \theta \geq 0}
    \Vert \tilde{\mathbf u}_a- \mathbf H\boldsymbol\theta\Vert^2
    +\beta \vert \boldsymbol\theta\vert_{\rm TV},
\end{equation}
where $\beta$ is the regularization parameter, 
and a non-negativity constraint is employed. 
The TV-norm is defined as 
\begin{equation}
  |\boldsymbol\theta|_{\rm TV}
     =\sum_{n=0}^{N-1}\sqrt{
      \Big( [\boldsymbol\theta]_n
              -[\boldsymbol\theta]_{k_{x_1}}\Big)^2 \\
            +\Big( [\boldsymbol\theta]_n
              -  [\boldsymbol\theta]_{k_{y_1}}\Big)^2\\
            +\Big(  [\boldsymbol\theta]_n
              -  [\boldsymbol\theta]_{k_{z_1}}\Big)^2
               }.
\end{equation}
We implemented the fast iterative shrinkage/thresholding algorithm (FISTA) 
to solve \eref{eqn:TVcost} \cite{FISTATIP:2009}, 
which will be referred to as PLS-TV algorithm. 

\section{Descriptions of numerical studies}
\label{sec:methods}

\subsection{Experimental data acquisition} 
\begin{figure}[h]
\centering
\subfigure[]{\includegraphics[width=8.0cm]{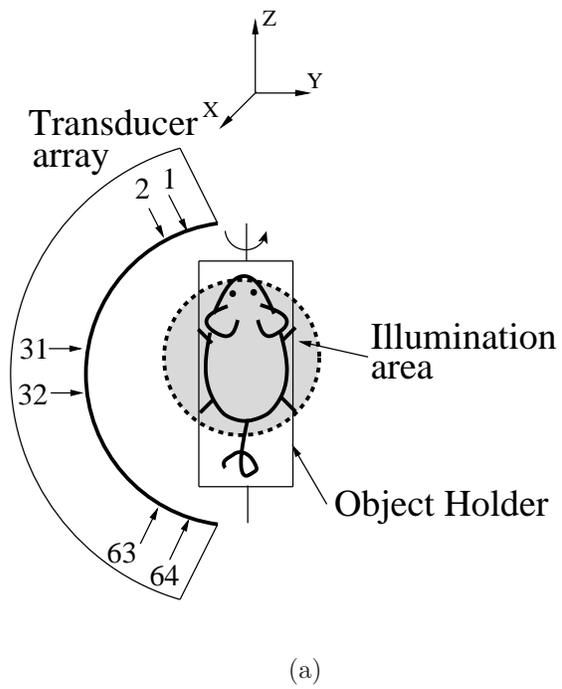}}
\hskip 1cm
\subfigure[]{\includegraphics[width=5.0cm]{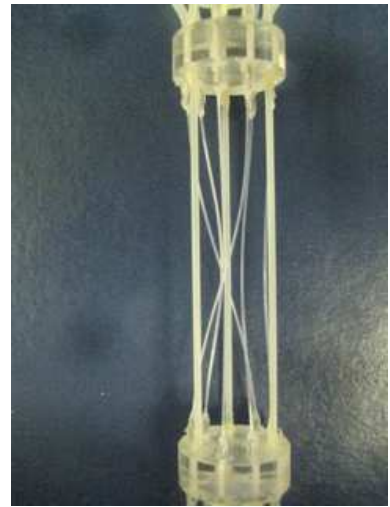}}
\caption{
(a) Schematic of the 3D OAT scanning geometry; 
(b) Photograph of the six-tube phantom.
\label{fig:geo}
}
\end{figure}

{\bf Scanning geometry:}
The small animal OAT imager employed in our studies has been described in previous publications
\cite{ermilov:2009,petermice:spie,petermice:jbo}.
As illustrated in \Fref{fig:geo}-(a), the arc-shaped probe  
consisted of $64$ transducers that spanned $152$ degrees over a circle of radius $65$-mm. 
Each transducer possessed a square detecting area of size $2\times 2$-mm$^2$. 
The laser illuminated the object from rectangular illumination bars 
in orthogonal mode. 
During scanning, the object was mounted on the object holder and rotated 
over the full $360$ degrees while the probe and light illumination stayed
stationary. 
Scans were set to sample at 20MHz over 1536 samples with an amplification 
of 60dB and 64 averages per acquisition. 

{\bf Six-tube phantom:}
A physical phantom was created that contained three pairs of polytetrafluoroethylene 
thin walled tubing of $0.81$-mm in diameter
that were filled with different concentrations of nickel sulfate solution having 
absorption coefficient values of $5.681$-cm$^{-1}$, $6.18$-cm$^{-1}$, and  $6.555$cm$^{-1}$. 
The tubes were held within a frame of two acrylic discs of 1'' diameter  
that were separated at a height of 3.25'' and kept attached by three garolite 
rods symmetrically spaced $120^\circ$ apart. 
The tubing was such that each pair would contain a tube that would follow 
along the outside of the phantom and the second would be diagonally 
inside.
A photograph of the phantom is shown in Fig. \ref{fig:geo}-(b).
The entire phantom was encased inside a thin latex membrane that was filled 
with skim milk to create an optically scattering medium.
A titanium sapphire laser with a peak at $765$-nm and a pulse  
width of 16ns (Quanta Systems) were employed to irradiate the phantom. 
The temperature of the water bath was kept at approximately $ 29.5^\circ$C with a water pump and 
heater. 
A complete tomographic data set was acquired by rotating the object about $360^\circ$ in 
$0.5^\circ$ steps. Accordingly, data were recorded by each transducer on the probe at 
$720$ tomographic view angles about the vertical axis. 

{\bf Mouse whole-body imaging:}
A 6 to 7 week old athymic Nude-Foxn1$^{\rm nu}$ live mouse (Harlan, Indianapolis, Indiana) 
was imaged with a similar setup to the phantom scan with a customized holder that 
provided air to the mouse when it was submerged in water. 
The holder was essentially comprised of three parts: 1) a hollow acrylic cylinder for breathing, 
2) an acrylic disc with hole for mouse tail and an apparatus to attach the legs, 
and 3) pre-tensioned fiber glass rods to connect the two acrylic pieces. 
The mouse was given pure oxygen with a flow rate of 2L/min with an additional 
2\% isoflurane concentration for initial anesthesia. 
During scanning the isoflurane was lowered to be around 1.5\%. 
The temperature of water was held constant at $34.7^\circ C$ with the use of a 
PID temperature controller connected to heat pads (Watlow Inc., Columbia, MO) 
underneath the water tank. A bifurcated optical fiber was attached to a 
ND:YAG laser (Brilliant, Quantel, Bozeman, MT) operating at $1064$-nm
wavelength with a energy pulse of about $100$-mJ during scans and a pulse duration of 
$15$-ns. The optical fiber outputs were circular beams of approximately $8$-cm at 
the target with an estimated $25$-mJ directly out of each fiber. 
Illumination was done in orthogonal mode along the sides of the water tank with 
in width of $16$". 
A complete tomographic data set was acquired by rotating the object 
about $360^\circ$ in $2^\circ$ steps.
Accordingly, data were recorded by each transducer on the probe at $180$ tomographic 
view angles about the vertical axis. 
More details regarding the data acquisition procedure can be found in 
\cite{petermice:spie,petermice:jbo}.

\subsection{Implementation of reconstruction algorithms}
{\bf Six-tube phantom:}
The region to-be-reconstructed was of size $19.8\times19.8\times50.0$-mm$^3$ and 
centered at $(-1.0, 0, -3.0)$-mm. 
The FBP algorithm was employed to determine estimates of $A(\mathbf r)$ that 
were sampled on a 3D Cartesian grid with spacing $0.1$-mm by use of a discretized 
form of  \eref{eqn:fbp}. 
The iterative reconstruction algorithms employed spherical voxels of $0.1$-mm in diameter
inscribed in the cuboids of the Cartesian grid. 
Accordingly, the reconstructed image matrices for all three algorithms 
were of size $198\times198\times500$. 
The speed-of-sound was set at $c_0=1.47$-mm/$\mu$s. 
We selected the Gr\"uneisen coefficient as $\Gamma=\beta c^2/C_p = 2,000$
of arbitrary units for all implementations. 
Since the top and bottom transducers received mainly noise for elongated structures 
that were aligned along z-axis, we utilized only the data that were acquired by the 
central 54 transducers. 

{\bf Mouse whole-body imaging:}
The region to-be-reconstructed was of size $29.4\times29.4\times61.6$-mm$^3$ and 
centered at $(0.49, 2.17, -2.73)$-mm.
The FBP algorithm was employed to determine estimates of $A(\mathbf r)$ that
were sampled on a 3D Cartesian grid with spacing $0.14$-mm by use of
\eref{eqn:fbp}.
The iterative reconstruction algorithms adopted spherical voxels of $0.14$-mm in diameter
inscribed in the cuboids of the Cartesian grid.
Accordingly, the reconstructed image matrices for all three algorithms
were of size $210\times210\times440$.
The speed-of-sound was chosen as $c_0=1.54$-mm/$\mu$s.
We selected the Gr\"uneisen coefficient as $\Gamma=\beta c^2/C_p = 2,000$
of arbitrary units for all implementations.
We utilized only the data that were acquired by the central 54 transducers.

{\bf Parallel programming by CUDA GPU computing:}
Three-dimensional iterative image reconstruction is computationally burdensome
in general. 
It demands even more computation when utilizing the system matrix defined 
by \eref{eqn:SysMatrix}, as opposed to the conventional 
spherical Radon transform model, mainly because the former accumulates contributions 
from more voxels to compute a single data sample. 
In addition, calculation of the SIR defined by \eref{eqn:sir} introduces
extra computation. 
It can take weeks to accomplish a single iteration by sequential programming 
using a single CPU, which is infeasible for practical applications. 
Because the calculation of SIR for each pair of transducer and voxel 
is mutually independent, we parallelized the SIR calculation by use of 
GPU computing with CUDA \cite{MRIGPU:2008,ChouMP:GPU}  such that multiple SIR samples were computed 
simultaneously, dramatically reducing the computational time. 
The six-tube phantom data were processed by use of 3 NVIDIA Tesla C2050 GPU cards, 
taking $4.52$-hours for one iteration from the data set containing $144$ tomographic views,  
while the mouse-imaging data were processed by use of 6 NVIDIA Tesla C1060 GPU cards, 
taking $5.73$-hours for one iteration from the data set containing $180$ tomographic views.
Though for testing we let the reconstruction algorithms iterate for over 100 iterations, 
both PLS-Q and PLS-TV usually converged within 20 iterations. 

\subsection{Image quality assessment}

{\bf Visual inspection:}
We examined both the 3D images and 2D sectional images. 
To avoid misinterpretations due to display colormap, we compared grayscale images. 
Also, for each comparison, we varied the grayscale window to ensure 
the observations are minimally dependent on the display methods. 
For each algorithm we reconstructed a series of images
corresponding to different values of regularization parameter over a wide range. 
To understand how image intensities are affected by the choice of regularization parameter, 
each 2D sectional image was displayed in the grayscale window that spanned from the minimum 
to the maximum of the determined image intensities.  

It is more challenging to fairly compare 3D images by visual inspection. 
Hence we intended not to draw conclusions on which algorithm was superior, 
but instead to reveal the similarities among algorithms when data were densely sampled. 
Although for each reconstruction algorithm we reconstructed a series of images corresponding to 
the values of regularization parameter over a wide range, 
the main structures within the images appeared very similar in general. 
Thus we selected a representitive 3D image for each reconstruction algorithms. 
These representative images correspond to the regularization 
parameters whose values were near the center of the full ranges 
and have a very close range of image intensities.
We displayed these images in the same grayscale window. 
For a prechosen grayscale window $[\theta_{\rm low}, \theta_{\rm up}]$,
the reconstructed images were truncated as:
\begin{equation}
  [\boldsymbol\theta^{\rm disp}]_n=\left\{\begin{array}{ll}
    \theta_{\rm low},& {\rm if}\,\, [\boldsymbol\theta]_n<\theta_{\rm low}\\
    \theta_{\rm up} ,& {\rm if}\,\, [\boldsymbol\theta]_n>\theta_{\rm up} \\
    {[\boldsymbol\theta]_n}, &{\rm otherwise} .
    \end{array}\right.
\end{equation}
The truncated data were linearly projected to the 
range $[0,255]$ as 8-bit unsigned integers:
\begin{equation}
  [\boldsymbol\theta^{\rm int8}]_n={\rm round}\Big(
         -\frac{255}{\theta_{\rm up}-\theta_{\rm low}}
           ( [\boldsymbol\theta^{\rm disp}]_n-\theta_{\rm low})
         \Big).
\end{equation}
The 3D image data $\boldsymbol\theta^{\rm int8}$ were visualized by computing
maximum intensity projection (MIP) images by use of the Osirix software \cite{OsiriX}.

{\bf Quantitative metrics:} 
Because the six-tube phantom contained nickel sulfate solution as the only optical absorber,
the tubes were interpreted as signals in the reconstructed images, which were 
contaminated by random noise, e.g., the electronic noise. 
Since the tubes were immersed in nearly pure scattering media, the reconstructed 
images were expected to have zero-mean background. 
In contrary, the mouse whole-body imaging possessed a nonzero-mean background 
because the absorbing capillaries within blood-rich structures were beyond the $0.5$-mm 
resolution limit \cite{petermice:jbo} of the imaging system, resulting
a diffuse background. 
Consequently, we interpreted the arteries and veins as signals, which 
were immersed in nonzero-mean background plus random noise. 

{\it Image resolution:}
Because both the tubes and blood vessels were fine threadlike objects, 
we quantified the spatial resolution by their thickness. 
To estimate the thickness of a threadlike object lying along z-axis at certain height, 
we first selected the 2D sectional image at that height. 
Subsequently, we truncated the 2D image into dimension of 
$(2N_r+1)$-by-$(2N_r+1)$ pixels; and adjusted the location of the truncated image
such that only a continuous group of pixels corresponding to the structure of interest, 
or hot spot, was present at the center. 
We then fitted  the 2D sectional image to a 2D Gaussian function given by: 
\begin{equation}
  G[n_1,n_2] = G[0,0]\exp(
         -\frac{n_1^2+n_2^2}{2\sigma_r^2}),
\end{equation}
where $n_1$ and $n_2$ denoted the indices of pixels in the 2D digital image 
with $n_1, n_2 = -N_r, -N_r+1, \cdots, N_r$ in units of pixel size, 
$G[0,0]$ was the peak value of the Gaussian function located in the center,
and $\sigma_r$ was the standard deviation of the Gaussian function to be 
estimated. 
We let $N_r=15$ and $N_r= 10$ for the six-tube phantom and the mouse 
imaging respectively. 
Finally, the estimated $\sigma_r$ was converted to full width at half maximum (FWHM)
as the spatial resolution measure, i.e.,  
\begin{equation}\label{eqn:sgres}
 {\rm FWHM}=2\sqrt{2\ln 2}\sigma_r.
\end{equation}

{\it Contrast-to-noise ratio (CNR):} 
For a prechosen structure, a series of adjacent 2D sectional images were selected 
along the structure (i.e, along z-axis) as described above. 
We collected the central voxel of each 2D image, forming the signal
region-of-interest (s-ROI). 
The signal intensity was calculated as:
\begin{equation}\label{eqn:sgmean}
  \bar{\boldsymbol\theta}^s=\frac{1}{N^s}\sum_{n=0}^{N^s-1}[\boldsymbol\theta^s]_n,
\end{equation}
where $N^s$ denoted the total number of voxels within the s-ROI. 
For the six-tube phantom, the s-ROI for each tube contained $N^s=200$ voxels
that extended from $z=-10.4$-mm to $z=9.6$-mm, 
while for the mouse-imaging study, the s-ROI for the vessel under study 
contained $N^s=20$ voxels that extended from $z=7.0$-mm to $z=9.8$-mm. 
For each s-ROI, we defined a background region-of-interest (b-ROI) 
that has the same dimension along z-axis as the s-ROI. 
For the six-tube phantom, we randomly selected $50$ voxels at 
every height that were within a circle of radius $5$-mm centered 
at the hot spot of the signal. 
Similarly, for the mouse-imaging study, we randomly selected $15$ voxels 
at every height that were within a circle of radius $2.1$-mm centered 
at the hot spot of the signal. 
Correspondingly, the b-ROI contained $N^b=10,000$ and $N^b=300$ voxels 
for the six-tube phantom and the mouse-imaging study, respectively. 
The background intensity was calculated by: 
\begin{equation}
  \bar{\boldsymbol\theta}^b=\frac{1}{N^b}\sum_{n=0}^{N^b-1}[\boldsymbol\theta^b]_n. 
\end{equation}
Also, the background standard deviation was calculated by: 
\begin{equation}\label{eqn:bgstd}
  \sigma_b = \sqrt{\frac{1}{N^b-1}\sum_{n=0}^{N^b-1} 
     \Big([\boldsymbol\theta^b]_n-\bar{\boldsymbol\theta}^b\Big)^2}.
\end{equation} 
Because the reconstructed image is not a realization of an ergodic random process, 
the value of $\sigma_b$ estimated from a single reconstructed image is not equivalent 
to the standard deviation of the ensemble of images reconstructed by use of a 
certain reconstruction algorithm. 
Nevertheless, the $\sigma_b$ can be employed as a summary measure of 
 the noise level in the reconstructed 
images. Consequently, the contrast-to-noise ratio (CNR) was calculated by: 
\begin{equation}
  {\rm CNR} = \frac{|\bar {\boldsymbol\theta}^s-\bar{\boldsymbol\theta}^b|}
   {\sigma_b}.
\end{equation}

{\it Plot of resolution against standard deviation:}
All three reconstruction algorithms possess regularization parameters 
that control the trade-offs among multiple aspects of image quality. 
A plot of resolution against standard deviation evaluates how much spatial resolution
is degraded by a regularization method during its noise suppression. 
To obtain this plot for each reconstruction algorithm, we swept the value of 
the regularization parameter. 
For each value, we reconstructed 3D images and quantified the spatial resolution
and noise level by use of \eref{eqn:sgres} and \eref{eqn:bgstd} respectively. 
The FWHM values calculated along the structure of interest were averaged 
as a summary measure of resolution, denoted by $\overline{\rm FWHM}$. 
The average was conducted over $20$-mm and $2.8$-mm for the six-tube phantom 
and the mouse imaging respectively. 
The $\overline{\rm FWHM}$ was plotted against the standard deviation ($\sigma_b$).

{\it Plot of signal intensity against standard deviation:} 
In addition to the trade-off between resolution and standard deviation, 
regularization parameters also control the trade-off between bias and standard deviation. 
In general, stronger regularization may introduce higher bias 
while more effectiviely reducing the variance of the reconstructed image. 
Because the true values of absorbed energy density were unavailable, 
we plotted the signal intensity against the image standard variation
that were calculated by use of \eref{eqn:sgmean} and \eref{eqn:bgstd}. 
From this plot, we compared the noise level of the reconstructed images with comparable 
image intensities and hence with comparable biases.  

\section{Experimental results}
\label{sec:results}
The data for the six-tube phantom and mouse whole-body imaging were collected at 
$720$ and $180$ view angles respectively, referred to as full data sets. 
In order to emulate the scans with reduced numbers of views, we undersampled the full data sets 
to subsets with different numbers of view angles that were evenly distributed 
over $2\pi$. 
These subsets will be referred to as `$N$-view data' sets, where $N$ is the number of view angles. 

\subsection{Six-tube phantom}
\begin{figure}[ht]
\centering
\subfigure[]{\includegraphics[width=5cm, trim=3cm 2.0cm 4cm 2cm, clip]
{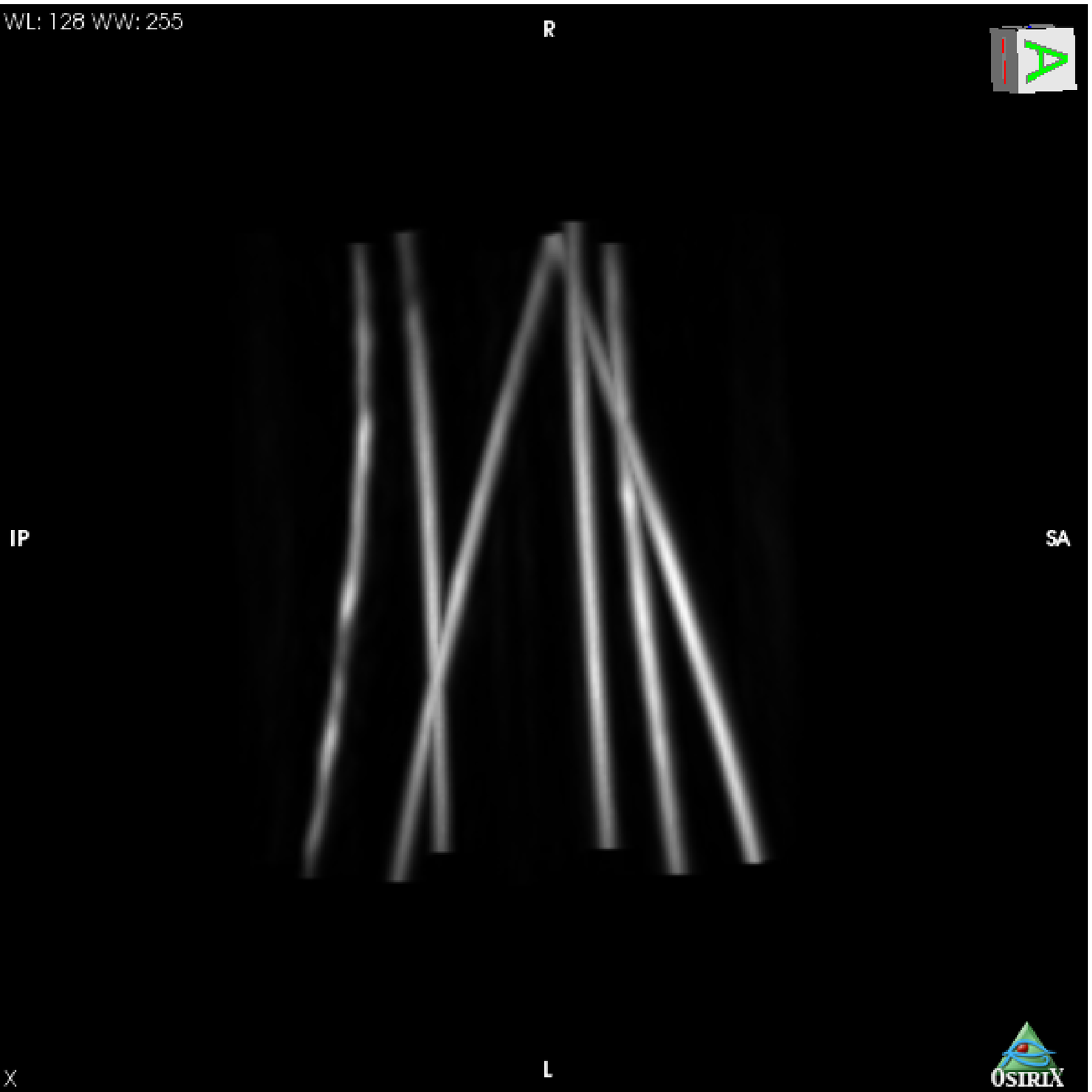}}
\hskip 1cm
\subfigure[]{\includegraphics[width=5cm, trim=3cm 2.0cm 4cm 2cm, clip]
{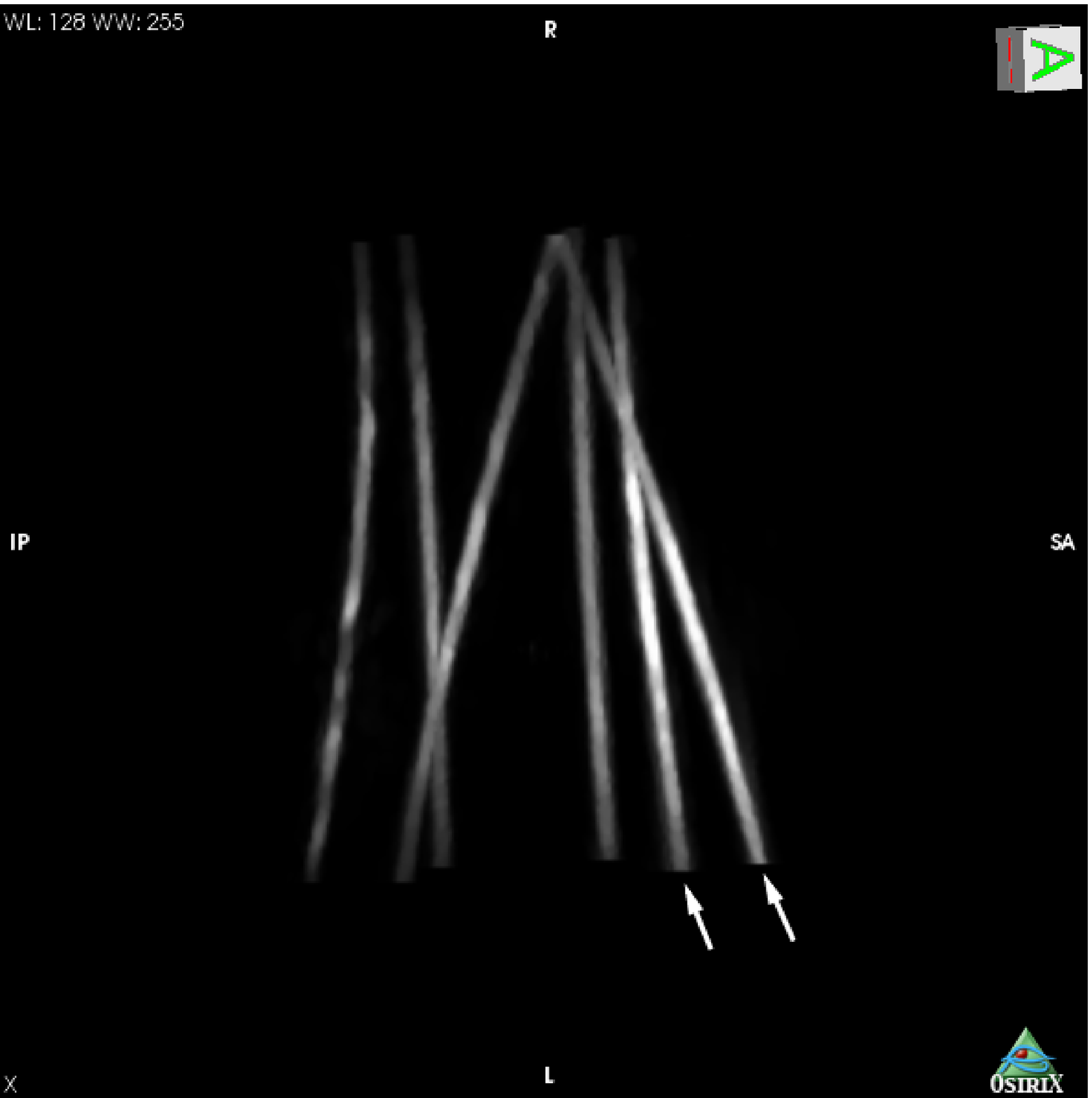}}
\caption{
MIP renderings of the six-tube phantom images reconstructed 
(a) from the 720-view data by use of the FBP algorithm with $f_c=6$-MHz; 
and
(b) from the 144-view data by use of the PLS-TV algorithm with $\lambda=0.1$. 
The grayscale window is [0,7.0]. 
Two arrows indicate the two tubes that were filled with the solution of the highest 
absorption coefficient $6.555$-cm$^{-1}$. (QuickTime)
\label{fig:SixTube3D}
}
\end{figure}

{\bf Visual inspection of reconstructed images from densely-sampled data sets:}
From densely-sampled data sets, the MIP images corresponding to the FBP and the PLS-TV 
algorithms appear to be very similar as shown in \Fref{fig:SixTube3D}. 
Note that the two images were reconstructed from different data sets: 
The image reconstructed by use of the FBP algorithm is from the full data set,
i.e., the $720$-view data set, 
while the one reconstructed by use of the PLS-TV algorithm is from the $144$-view data set. 
We did not apply iterative reconstruction algorithms to the $720$-view data set mainly 
because of the computational burden. 
Moreover, the images reconstructed from the $144$-view data set by use of the PLS-TV algorithm 
already appear to be at least comparable with those reconstructed by use of the FBP algorithm 
from full data set. 
Certain features are shared by both images. For example, 
both images contain two tubes (indicated by white arrows) that are brighter 
than the others, which is consistent with the 
fact that these two tubes are filled with the solution of higher absorption coefficient 
($\mu_a=6.555$cm$^{-1}$). 
The similarities between the two images are not surprising for two reasons:
Firstly, when the pressure function is densely sampled and the object is near the center of 
the measurement sphere, where the SIR can be neglected, we would expect all three algorithms 
to perform similarly because the imaging models they are based upon are equivalent in the 
continuous limit;  
Also the process of forming the MIP images strongly attenuates the background artifacts.

\begin{figure}[ht]
\centering
\subfigure[]{\includegraphics[width=5.1cm]{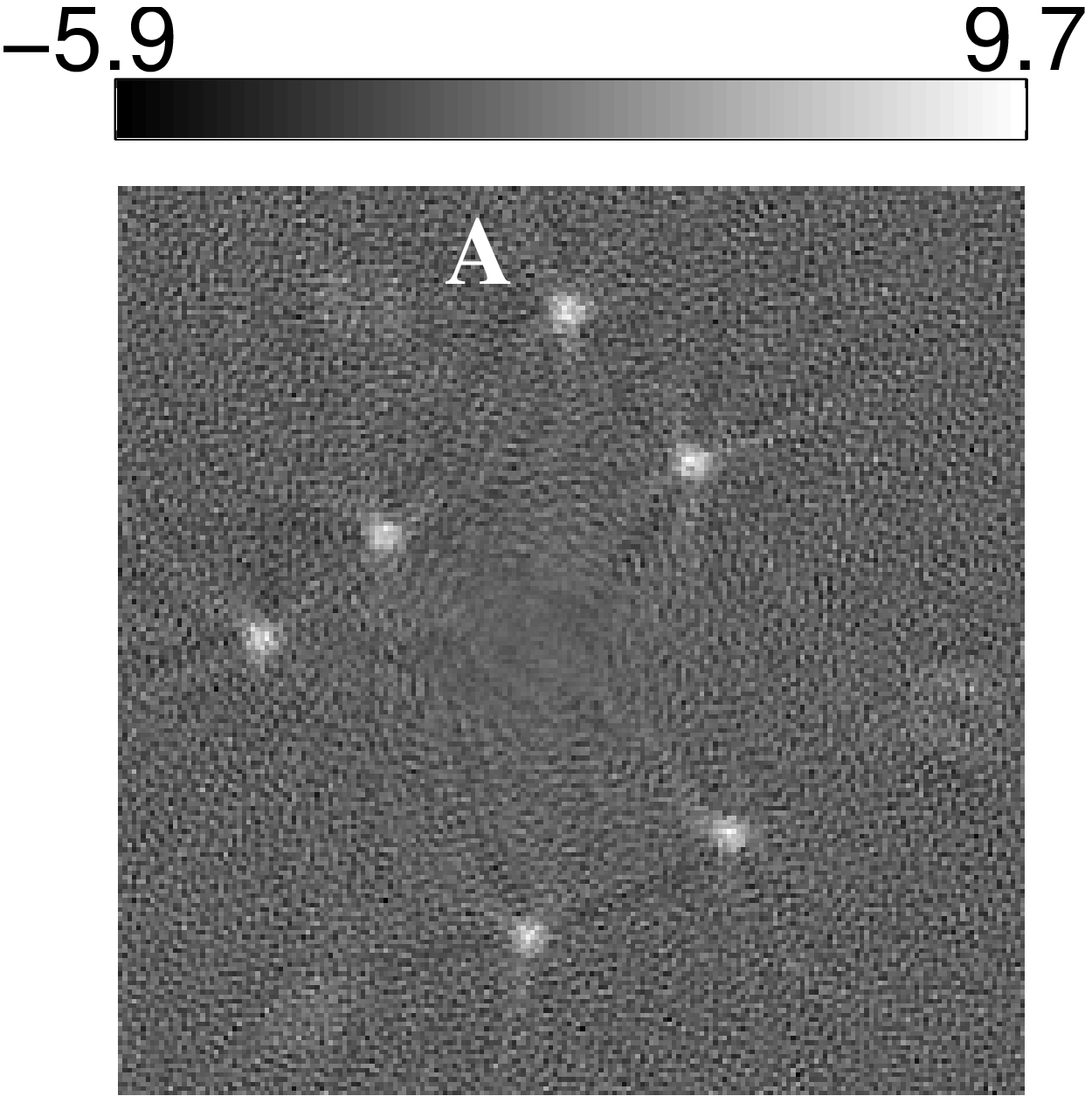}}
\subfigure[]
{\includegraphics[width=5.1cm]{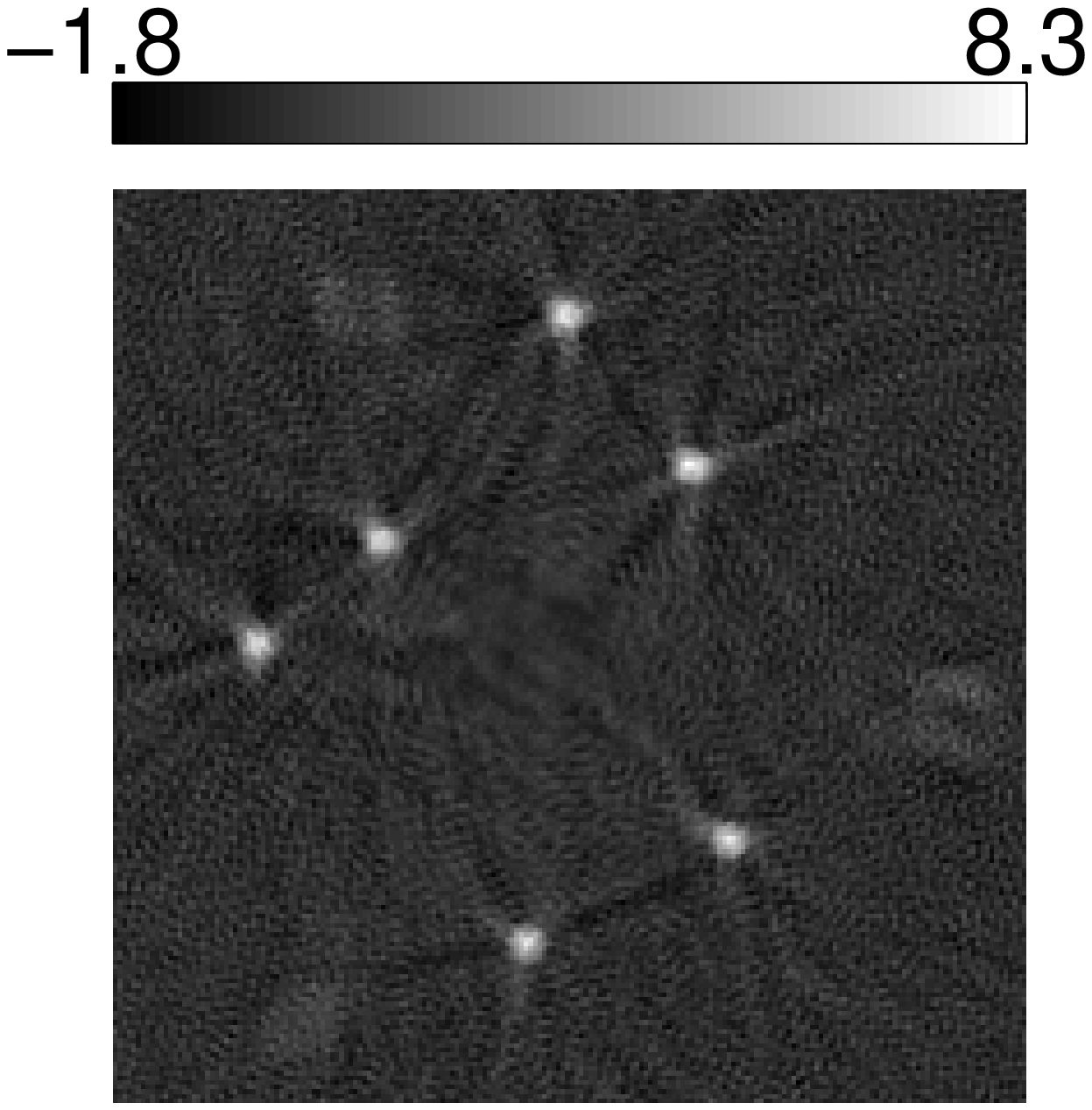}}
\subfigure[]
{\includegraphics[width=5.1cm]{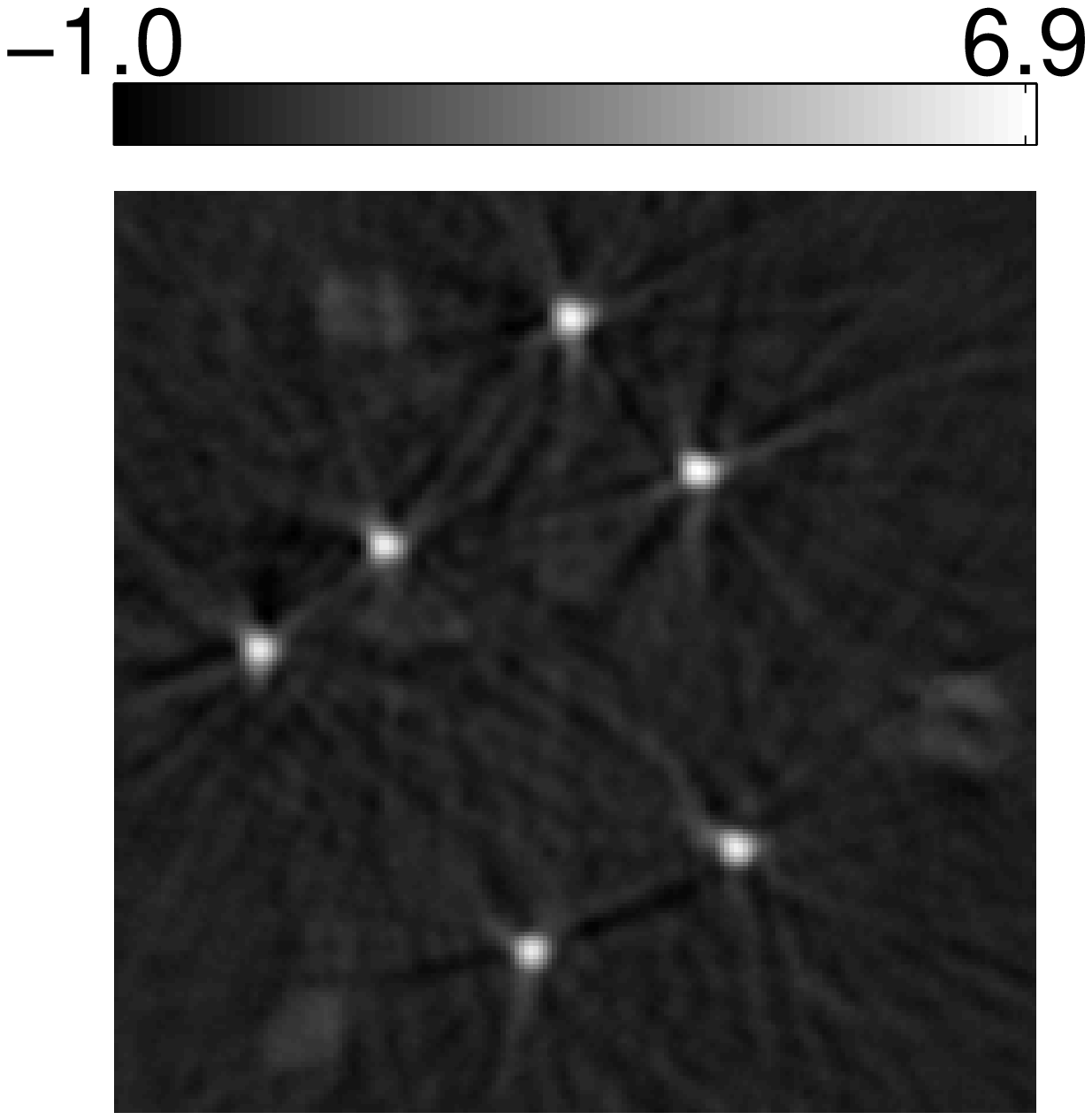}}\\
\subfigure[]
{\includegraphics[width=5.1cm]{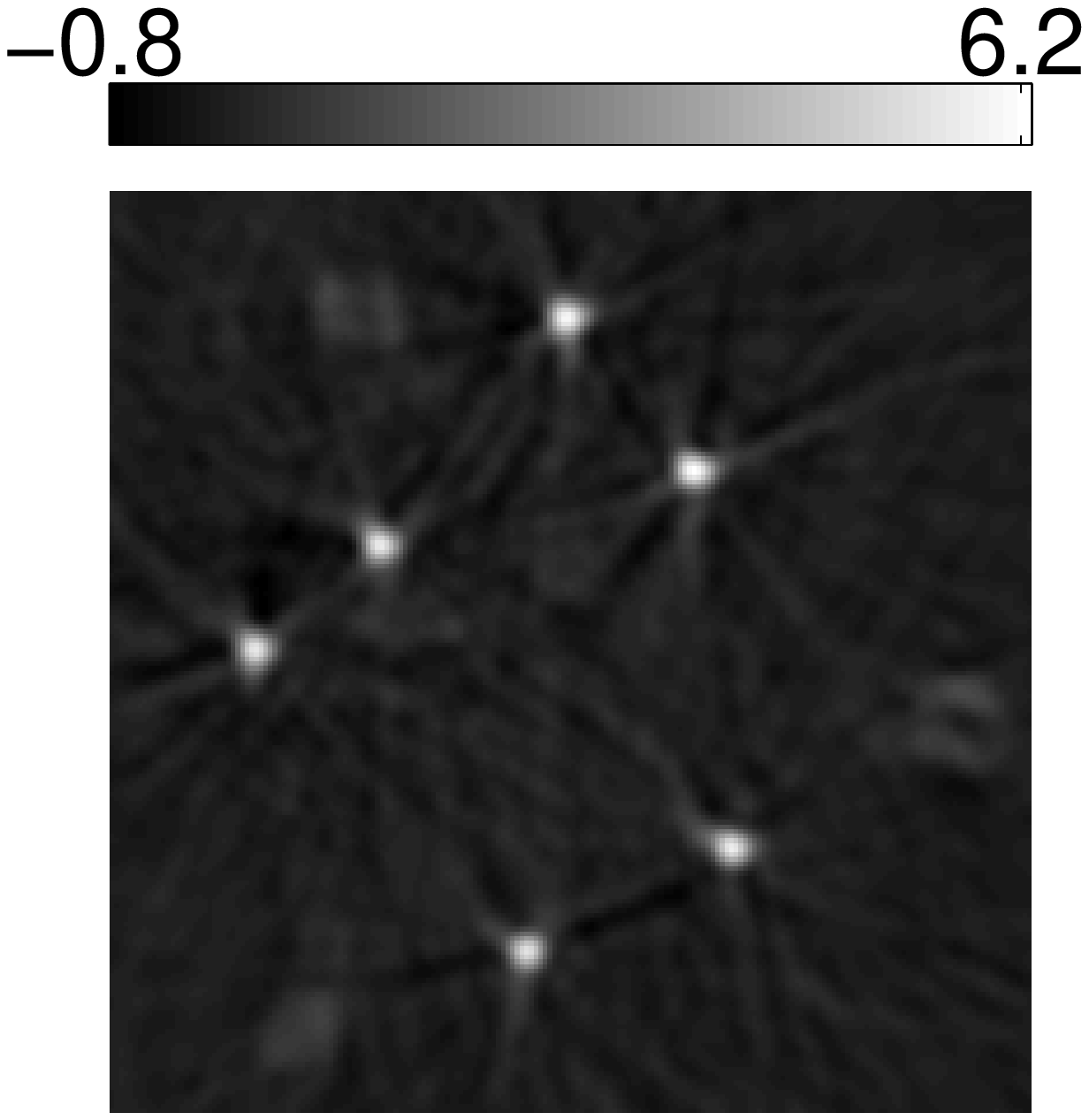}}
\subfigure[]
{\includegraphics[width=5.1cm]{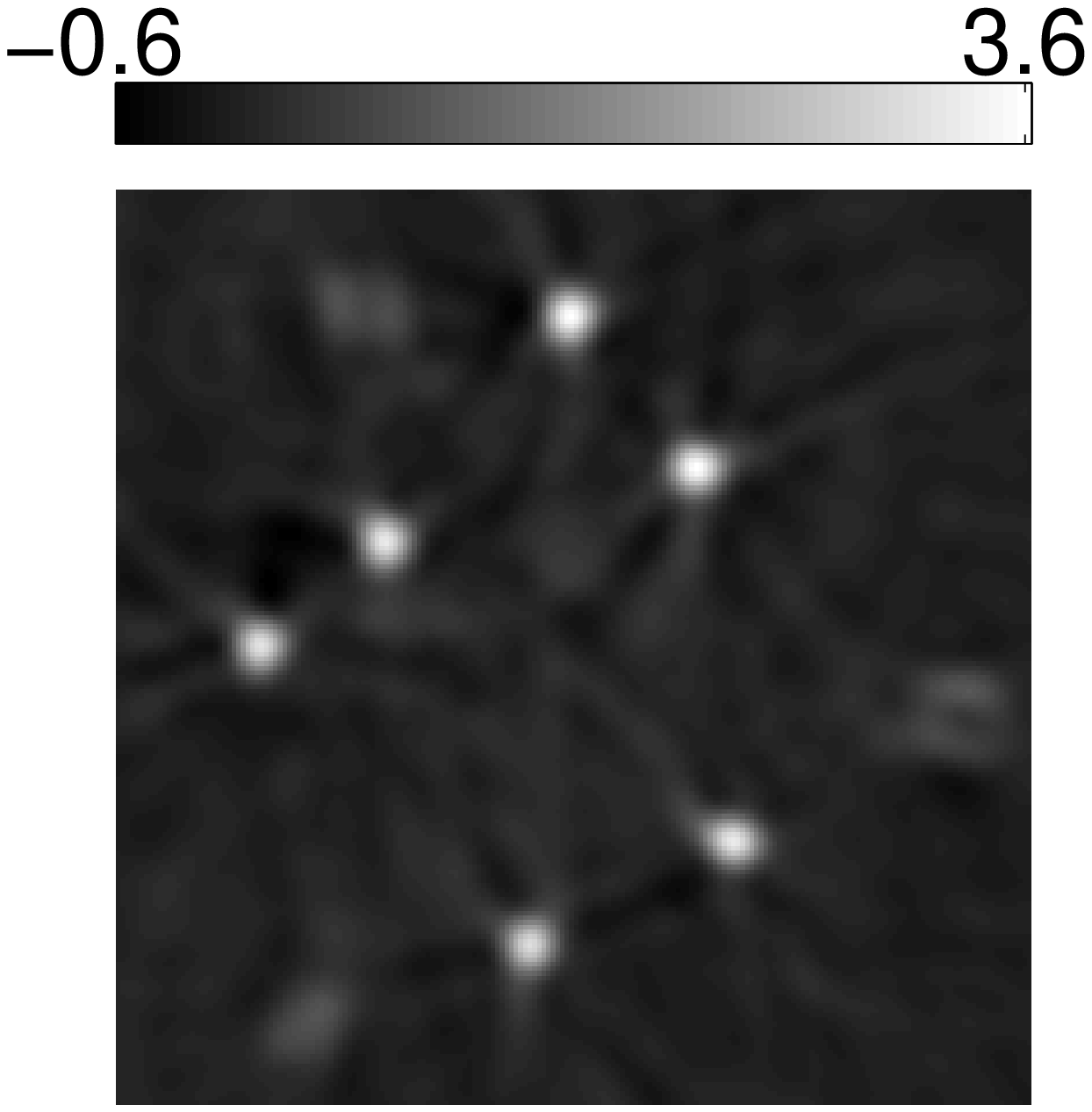}}
\subfigure[]
{\includegraphics[width=5.1cm]{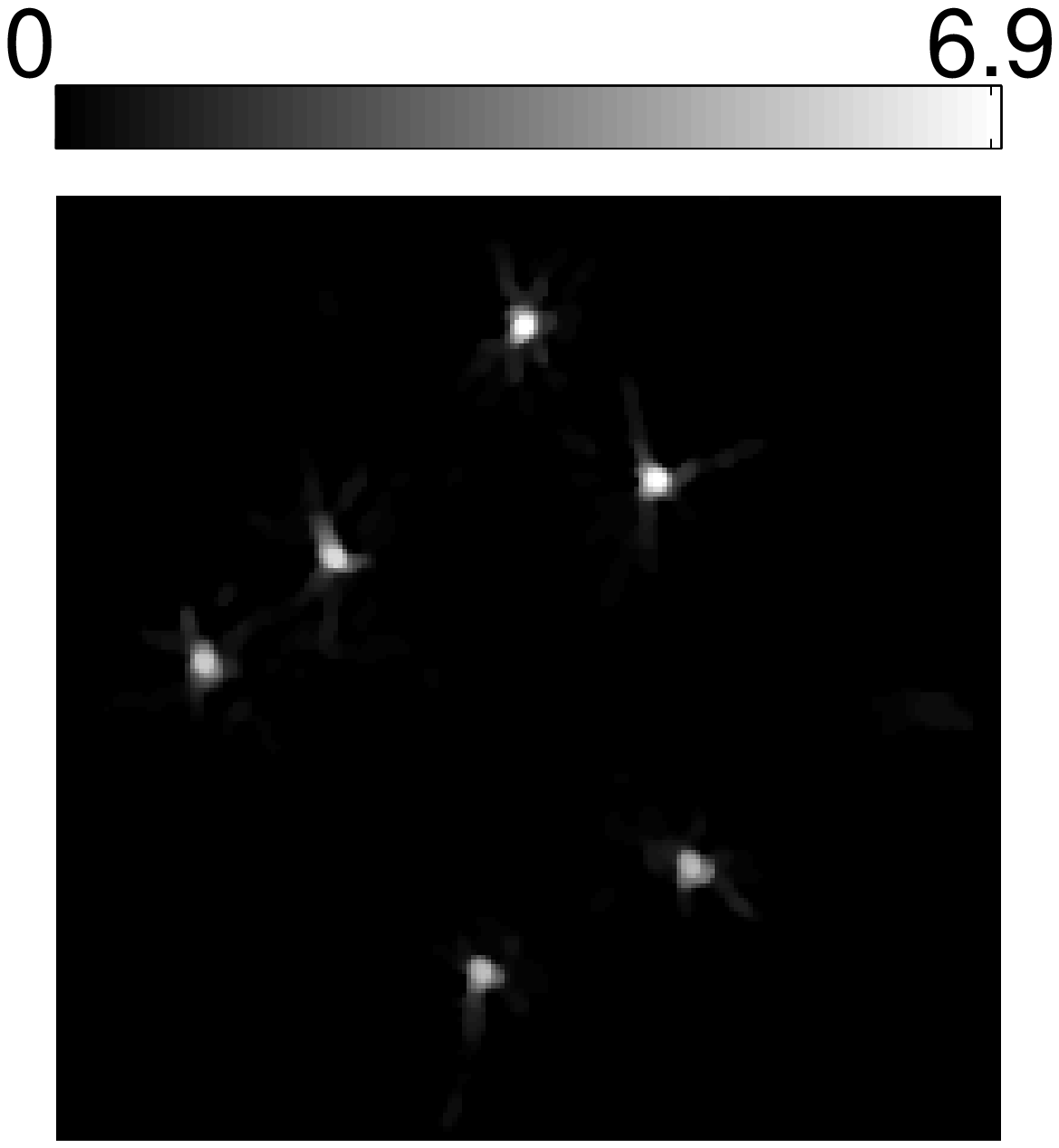}}
\caption{
Slices corresponding to the plane $z=-2.0$-mm through the 3D images of the six-tube phantom 
reconstructed  from 
(a) the $720$-view data by use of the FBP algorithm with $f_c=10$-MHz;
(b) the $720$-view data by use of the FBP algorithm with $f_c=8$-MHz;
(c) the $720$-view data by use of the FBP algorithm with $f_c=6$-MHz;
(d) the $720$-view data by use of the FBP algorithm with $f_c=4$-MHz;
(e) the $720$-view data by use of the FBP algorithm with $f_c=2$-MHz;
and 
(f) the $144$-view data by use of the PLS-TV algorithm with $\beta=0.1$.
All images are of size $19.8\times 19.8$-mm$^2$. 
The ranges of the grayscale windows were determined by the minimum and the maximum values in each 
image. 
\label{fig:SixTubeDenseSample}
}
\end{figure}

However, 2D sections of the 3D images reveal cerrtain favorable 
characteristics of the PLS-TV algorithm, as shown in \Fref{fig:SixTubeDenseSample}. 
Though we varied the cutoff frequency $f_c$ over a wide range for the FBP algorithm, 
none of these images has background as clean as the image reconstructed by the PLS-TV 
algorithm. 
We notice two types of artifacts in the images reconstructed by use of the FBP algorithm:
the random noise and the radial streaks centered at the tubes. 
The former is caused by measurement noise 
while the latter is likely due to certain unmodeled system inconsistencies that are 
referred to as systematic artifacts and will be addressed in \Sref{sec:summary}. 
The regularizing low-pass filter mitigates the random noise but also degrades 
the spatial resolution (\Fref{fig:SixTubeDenseSample}-b-e). 
The TV-norm regularization mitigates the background artifacts with minimal sacrifice in 
spatial resolution. 
The image reconstructed by use of the PLS-TV algorithm shown in \Fref{fig:SixTubeDenseSample}-(f) has 
at least comparable resolution as that of the FBP image 
with $f_c=6$-MHz (\Fref{fig:SixTubeDenseSample}-c).

\begin{figure}[ht]
\begin{center}
\subfigure[]{\includegraphics[width=5.1cm]{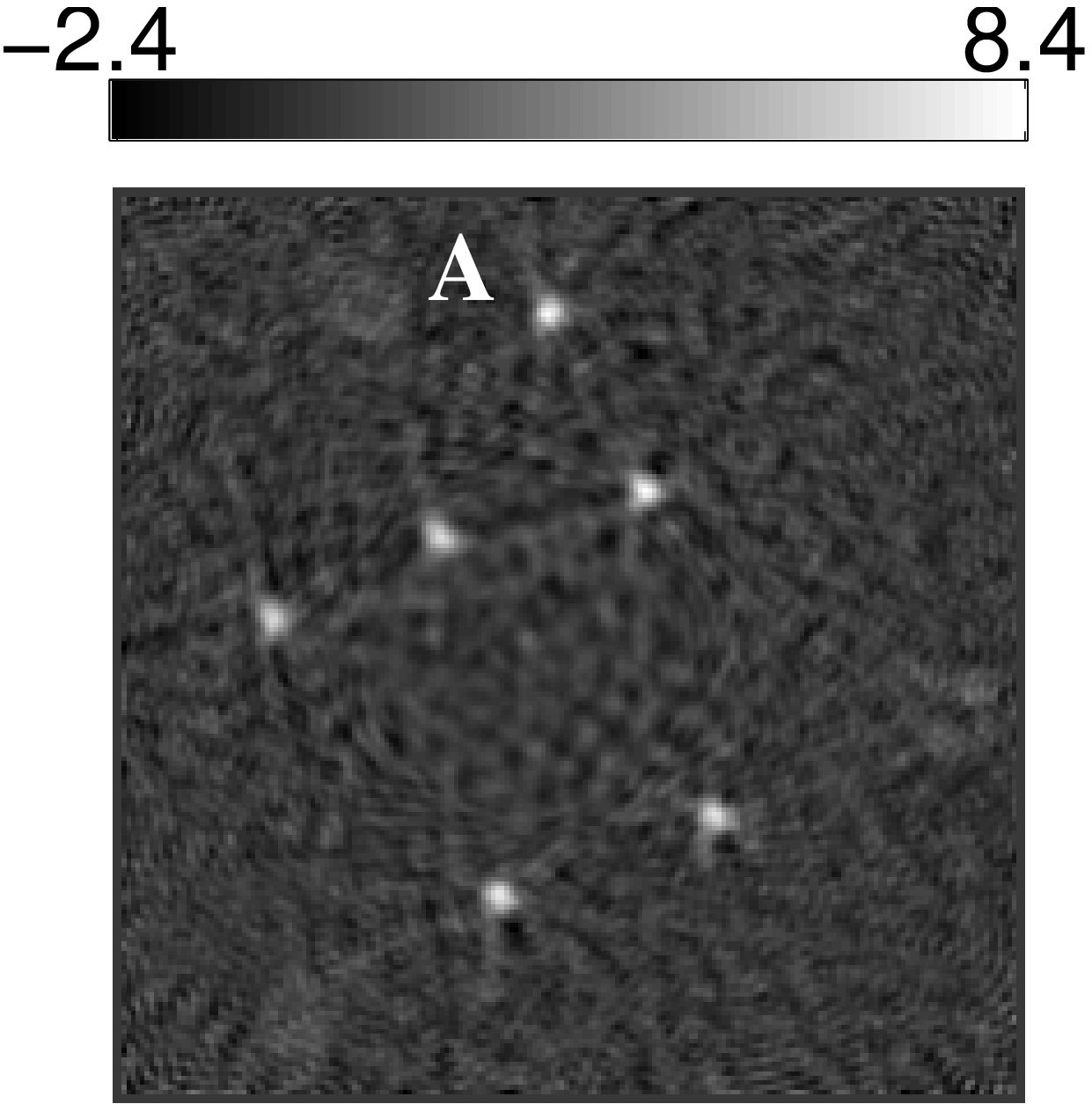}}
\subfigure[]
{\includegraphics[width=5.1cm]{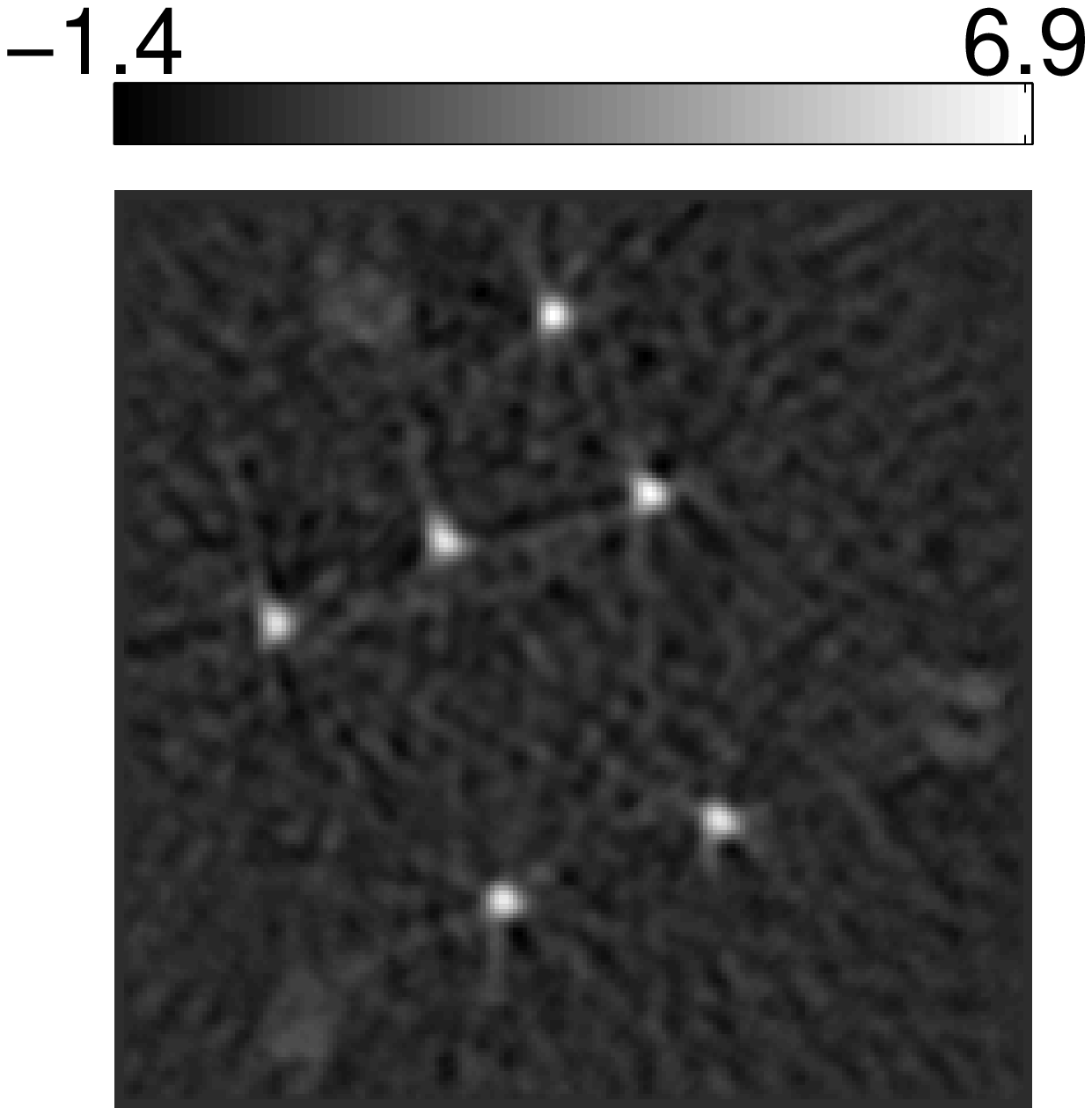}}
\subfigure[]
{\includegraphics[width=5.1cm]{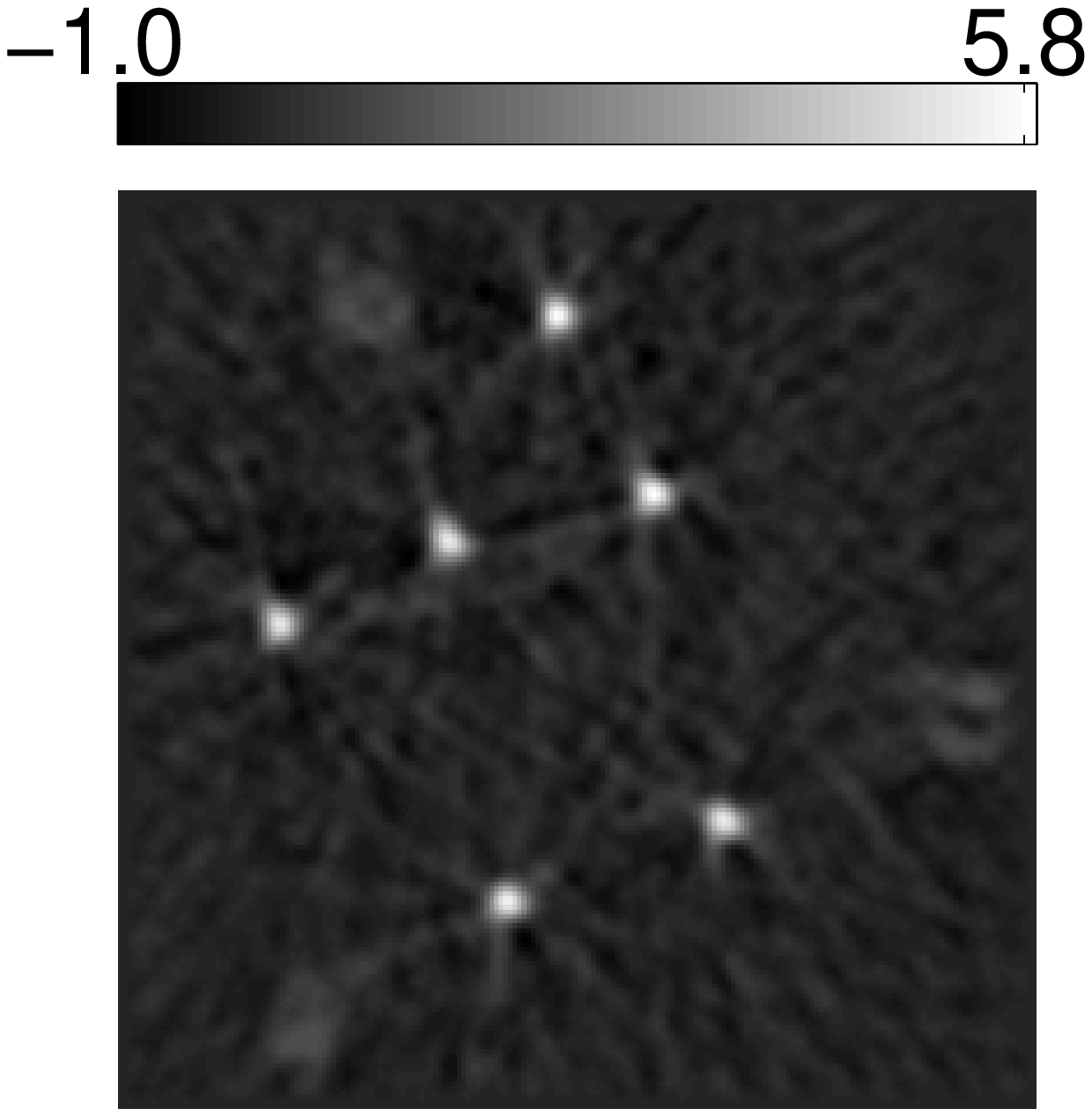}}\\
\subfigure[]
{\includegraphics[width=5.1cm]{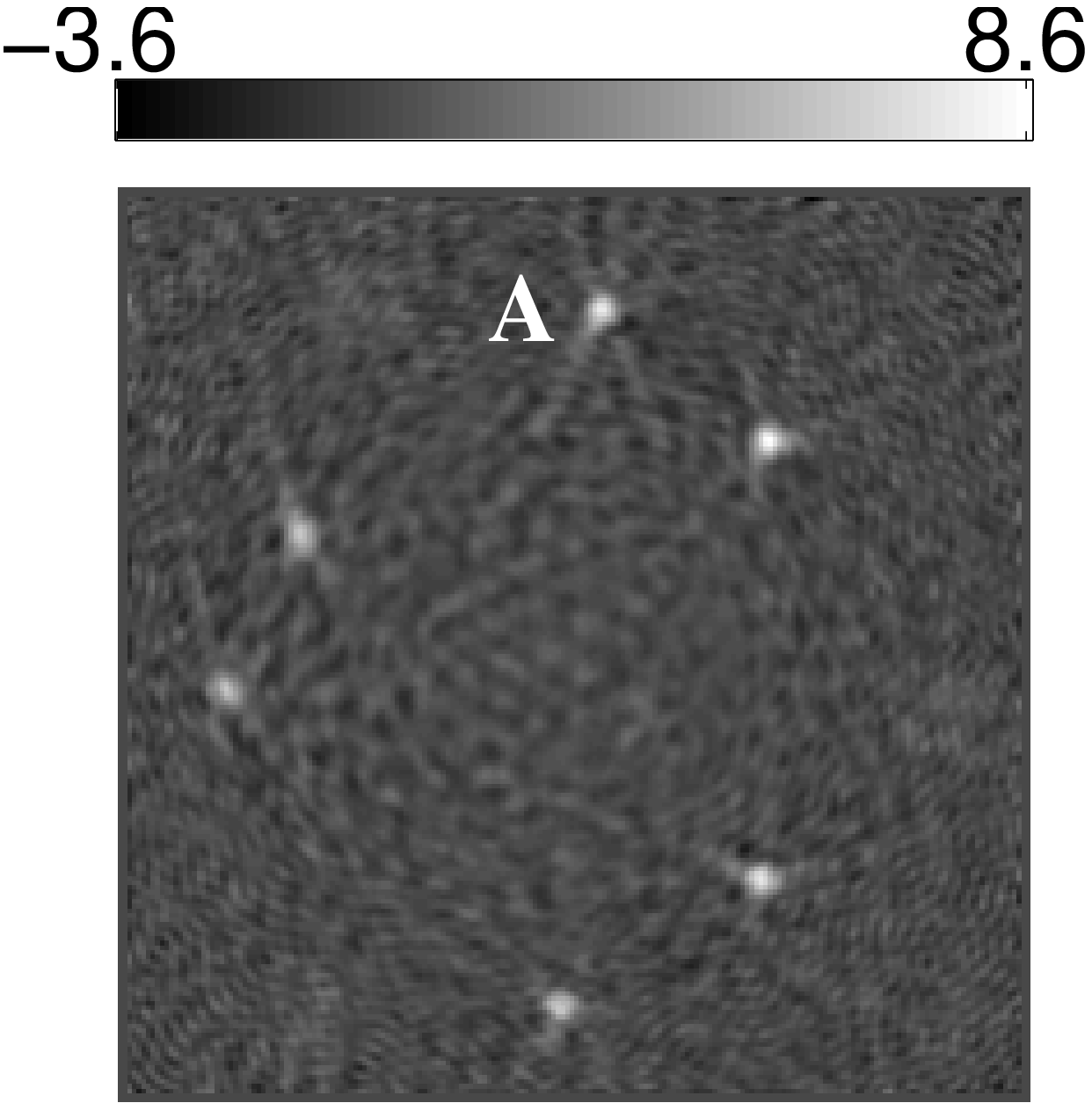}}
\subfigure[]
{\includegraphics[width=5.1cm]{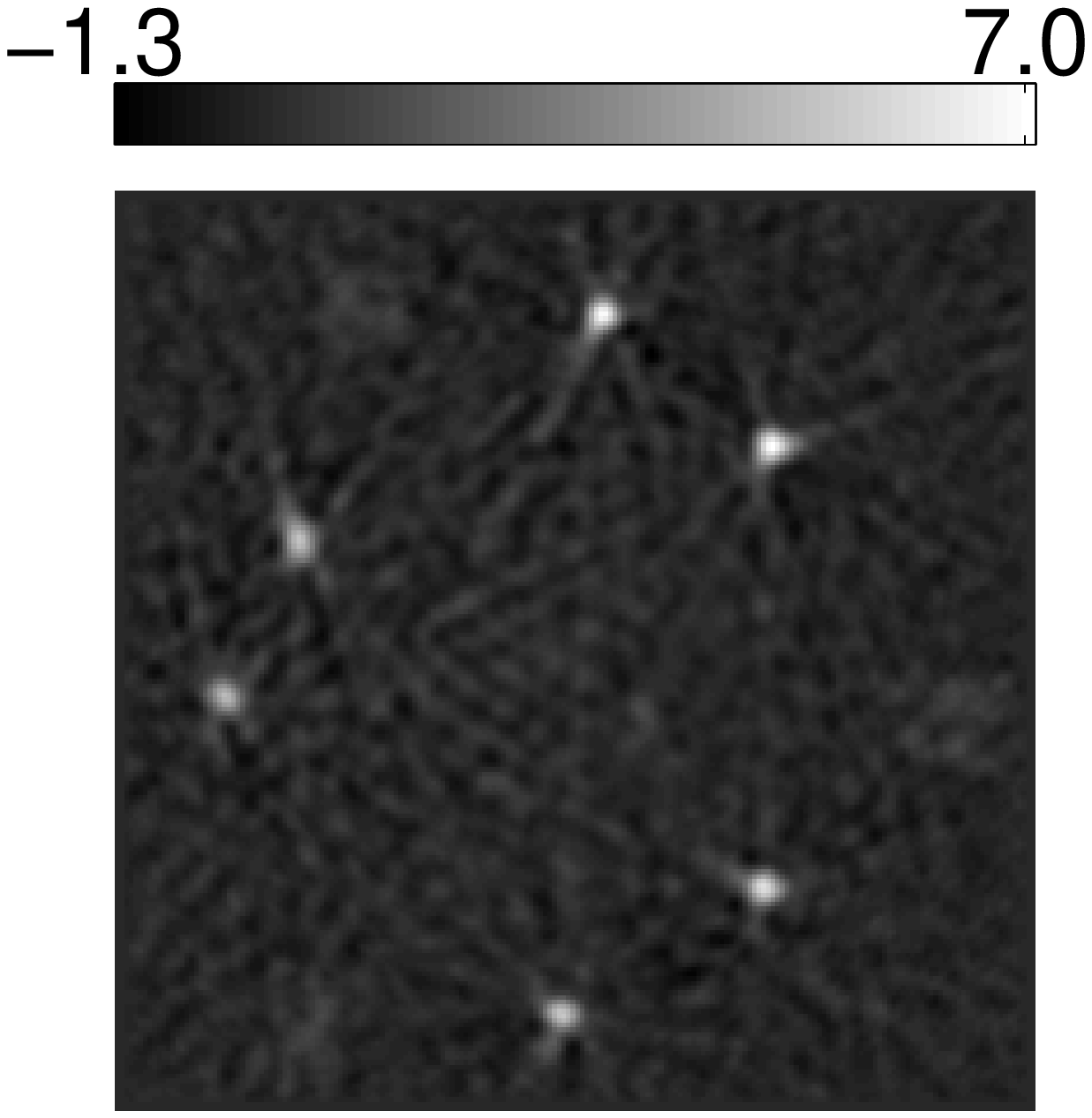}}
\subfigure[]
{\includegraphics[width=5.1cm]{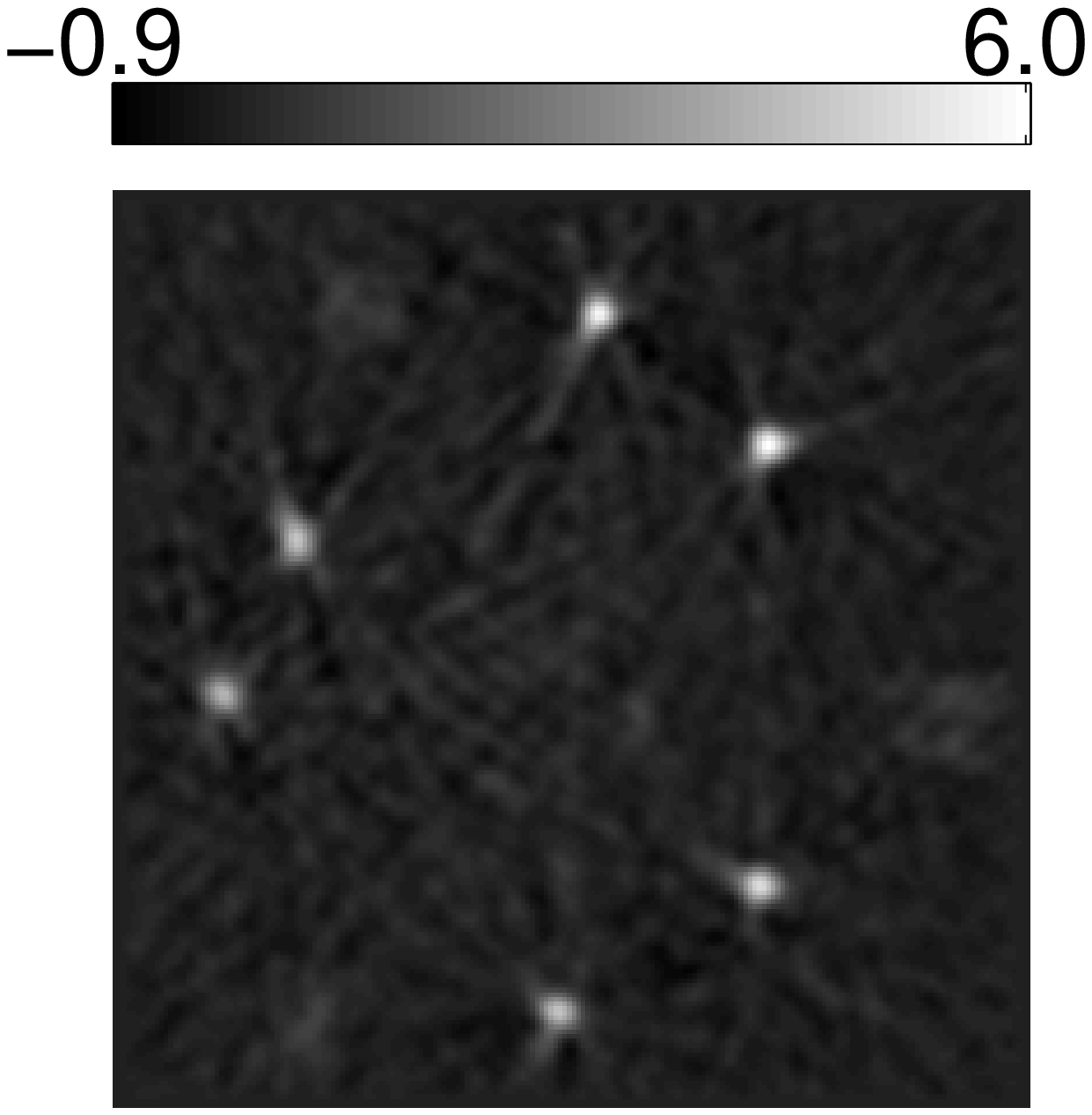}}
\end{center}
\caption{
Slices corresponding to the plane $z=-6.0$-mm (top row: a-c) and the plane $z=4.5$-mm 
(bottom row: d-f)
through the 3D images of the six-tube phantom reconstructed 
from the $144$-view data by use of the PLS-Q algorithm with varying regularization
parameter $\alpha$:  
(a), (d) $\alpha=0$;
(b), (e) $\alpha=1.0\times 10^3$;
and
(c), (f) $\alpha=5.0\times 10^3$;
All images are of size $19.8\times 19.8$-mm$^2$. 
The ranges of the grayscale windows were determined by the minimum and the maximum values 
in each image. 
\label{fig:SixTubePLSReg}
}
\end{figure}

\begin{figure}[ht]
\begin{center}
\subfigure[]{\includegraphics[width=5.1cm]{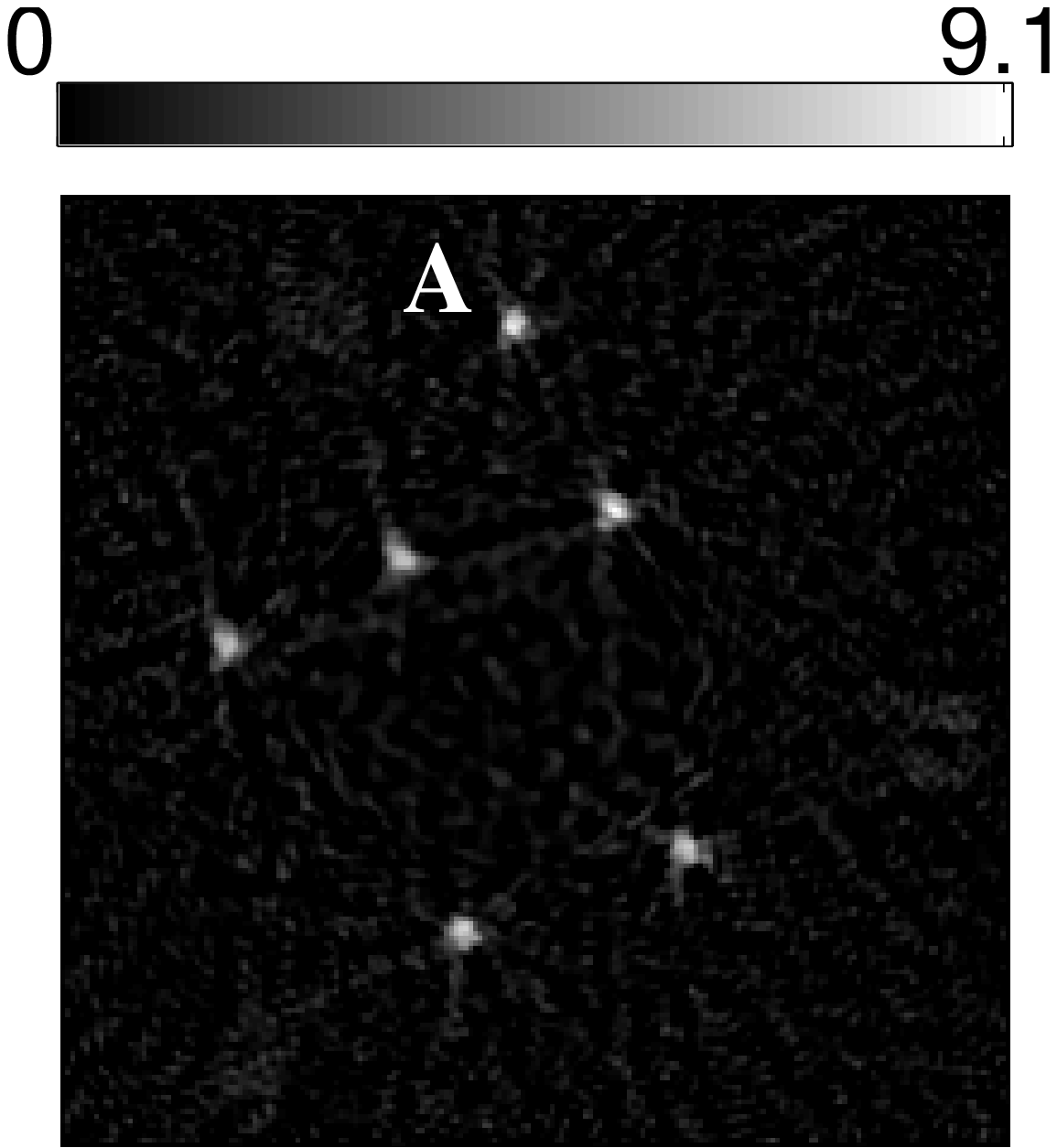}}
\subfigure[]
{\includegraphics[width=5.1cm]{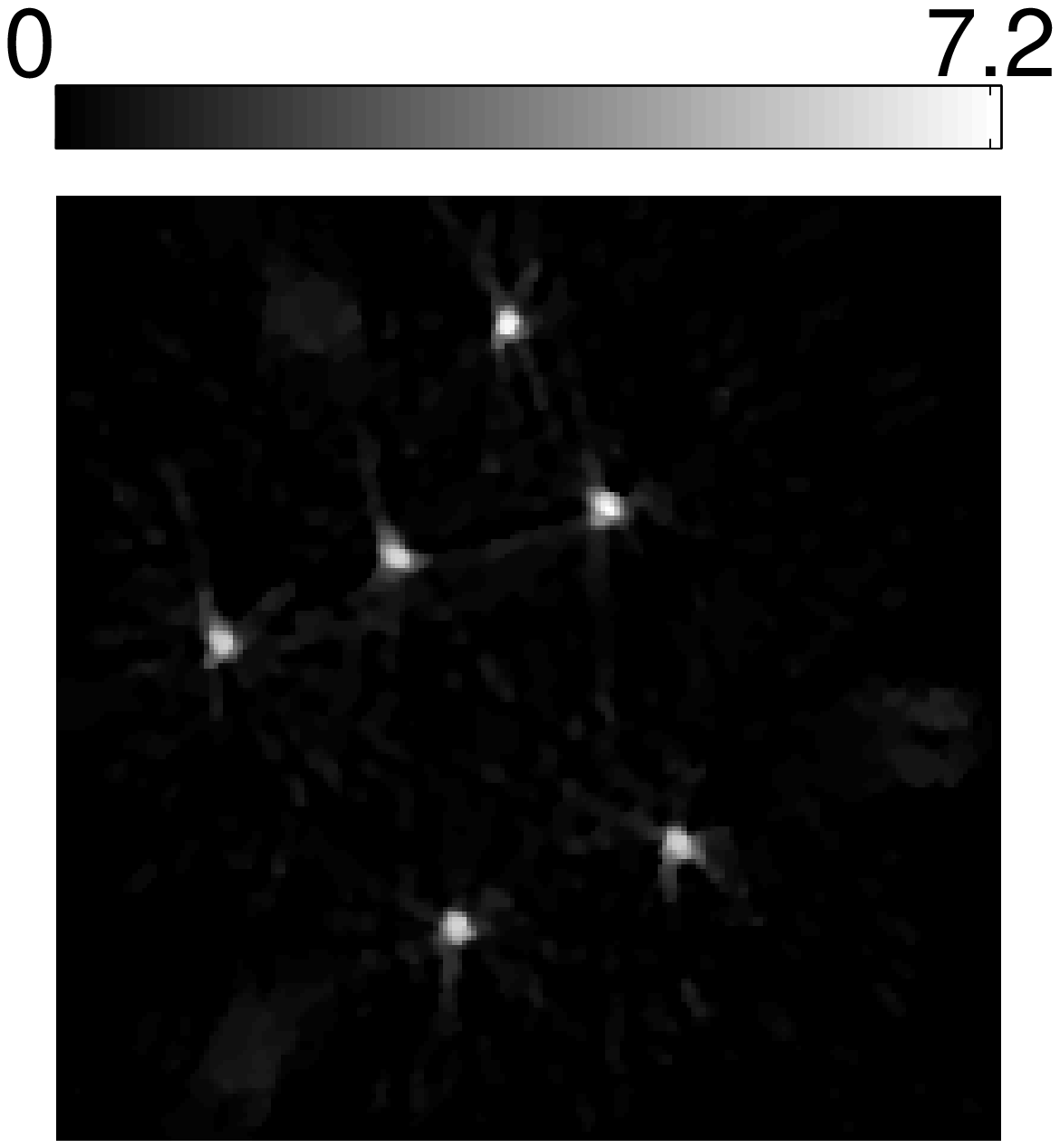}}
\subfigure[]
{\includegraphics[width=5.1cm]{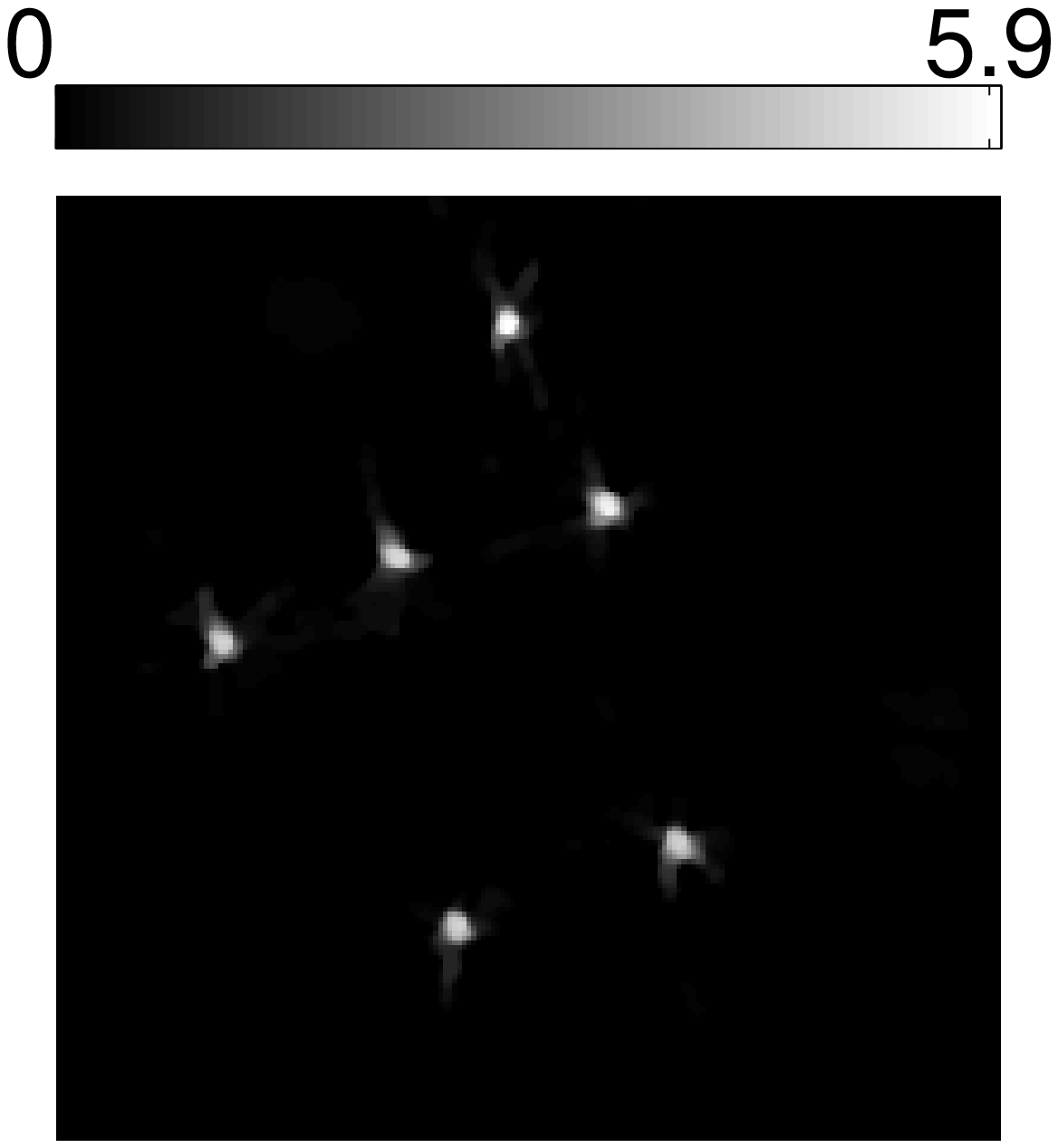}}\\
\subfigure[]
{\includegraphics[width=5.1cm]{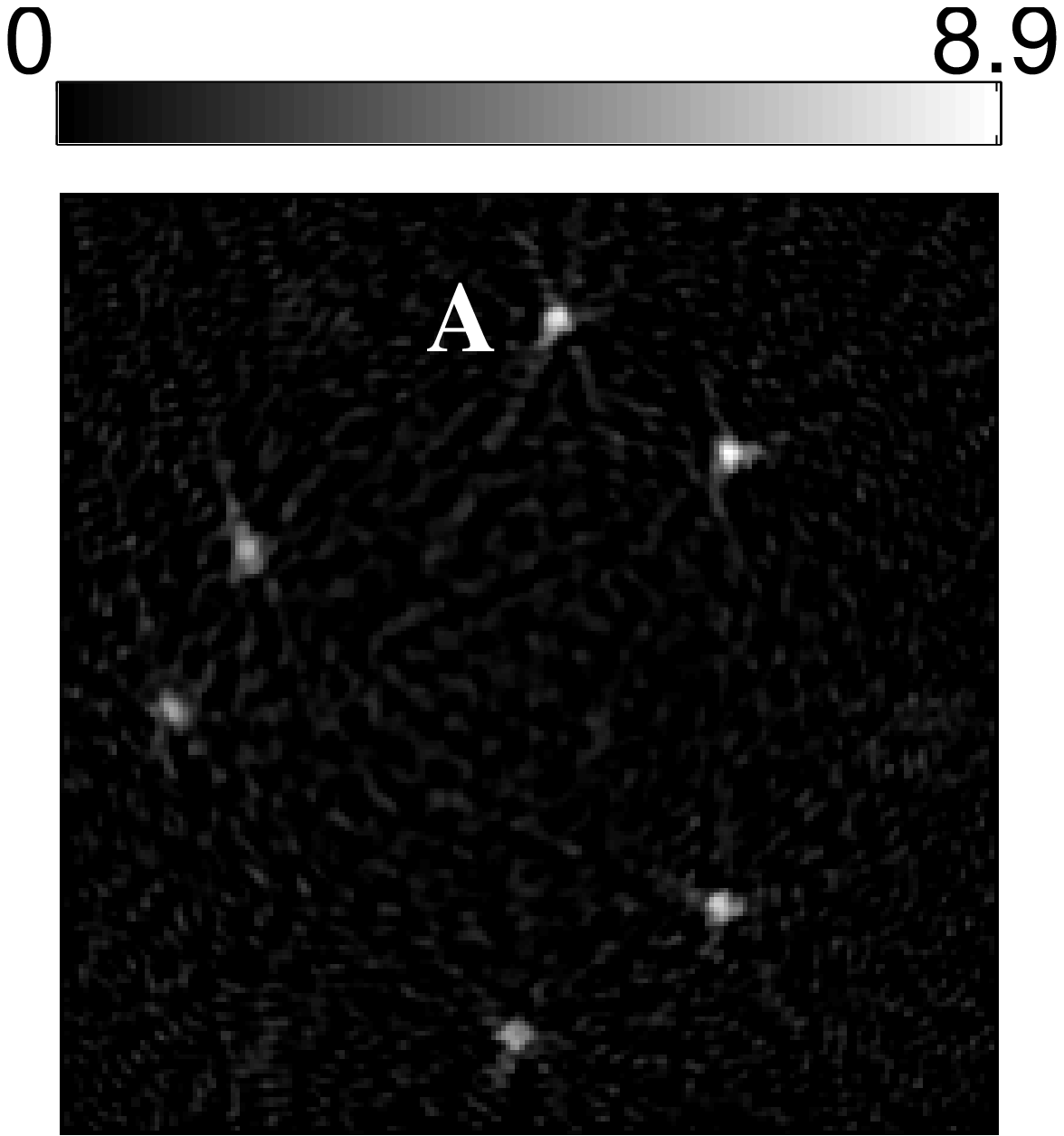}}
\subfigure[]
{\includegraphics[width=5.1cm]{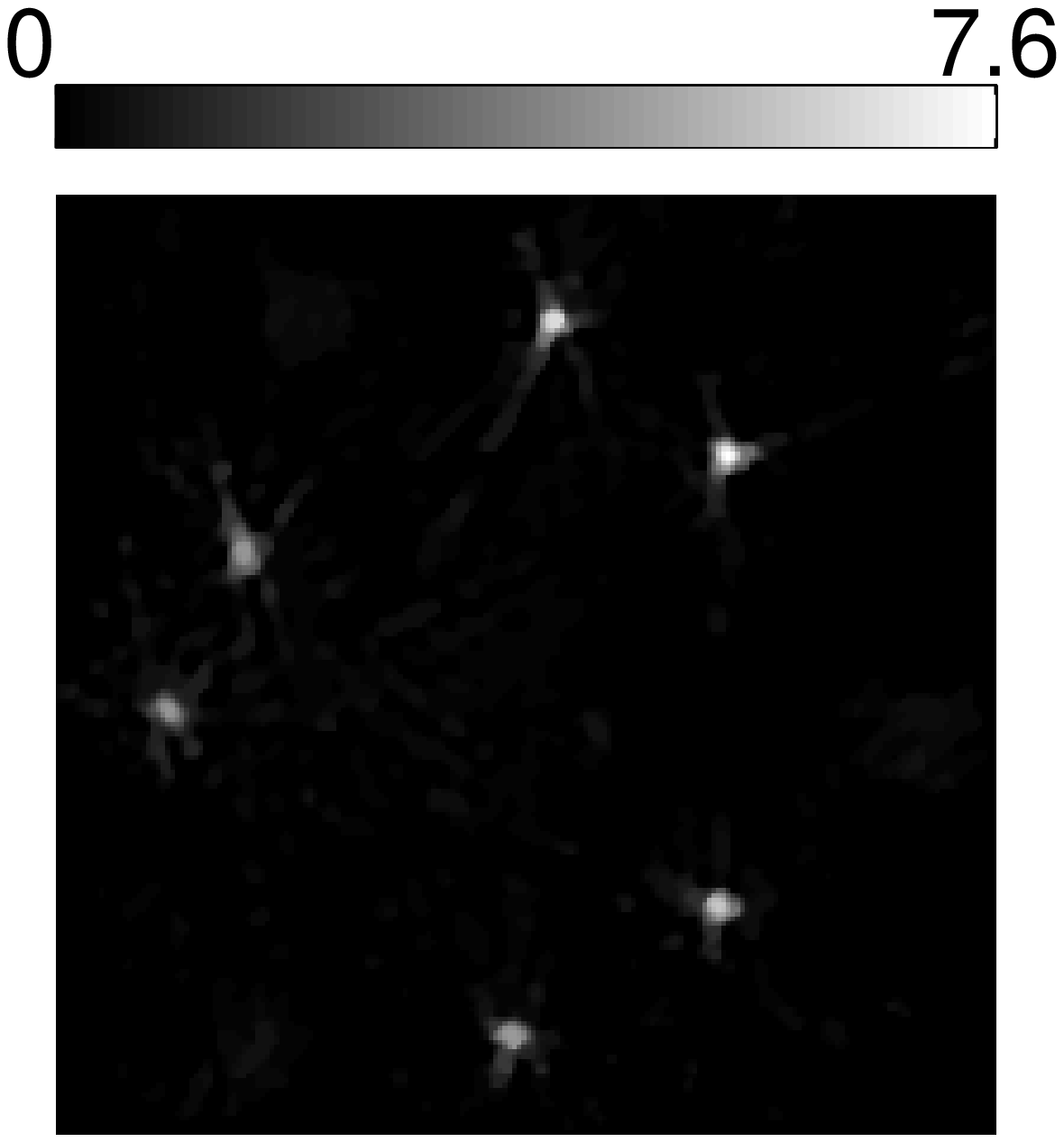}}
\subfigure[]
{\includegraphics[width=5.1cm]{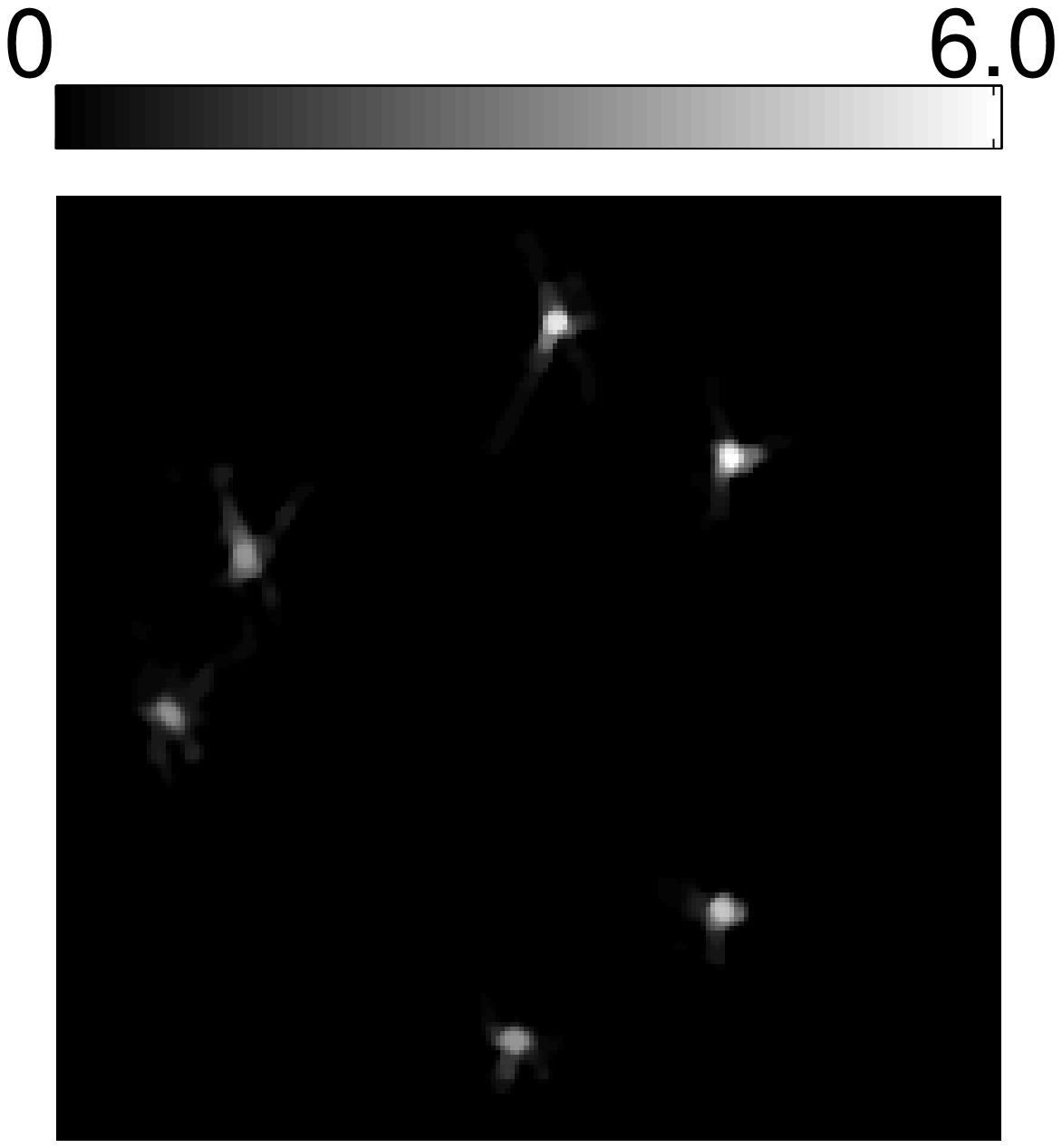}}
\end{center}
\caption{
Slices corresponding to the plane  $z=-6.0$-mm (top row: a-c) and 
the plane $z=4.5$-mm (bottom row: d-f) through the 3D images of the six-tube phantom reconstructed 
from the $144$-view data by use of the PLS-TV algorithm with varying regularization parameter 
$\beta$: 
(a), (d) $\beta=0.001$;
(b), (e) $\beta=0.05$;
and
(c), (f) $\beta=0.1$;
All images are of size $19.8\times 19.8$-mm$^2$. 
The ranges of the grayscale windows were determined by the minimum and the maximum values in each 
image. 
\label{fig:SixTubeTVReg}
}
\end{figure}

{\bf Qualitative comparison of regularization methods:}
The three reconstruction algorithms are regularized by use of the low-pass filter, 
the $\ell_2$-norm smoothness penalty and the TV-norm penalty, respectively. 
The impacts of the low-pass filter are revealed in \Fref{fig:SixTubeDenseSample}.
We observe that a slight regularization (i.e., a large value of $f_c$) 
results in sharp but noisy images 
while a heavy regularization (i.e., a small value of $f_c$) 
produces clean but blurry images. 
Also, the intensities of the tubes are lower for a smaller value 
of $f_c$. 
Similar effects are observed for the PLS-Q algorithm with the $\ell_2$-norm smoothness penalty 
as shown in \Fref{fig:SixTubePLSReg}. 
These observations agree with the conventional understandings of 
the impacts of regularization summarized as two trade-offs: 
resolution versus variance and bias versus variance \cite{Fessler:94}. 
The TV-norm regularization, however, mitigates the image variance with minimal 
sacrifice in image resolution as shown in \Fref{fig:SixTubeTVReg}. 
Though the signal intensity is reduced at $\beta=0.1$ 
(\Fref{fig:SixTubeTVReg}-c and -f), the spatial resolution appears to be  
very close to that of the images corresponding to smaller values of $\beta$ 
(\Fref{fig:SixTubeTVReg}-a and -d). 
In addition, both the low-pass filter and the $\ell_2$-norm penalty 
have little effects on the systematic artifacts 
while the TV algorithm effectively mitigates both the systematic artifacts 
and the random measurement noise.

\begin{figure}[ht]
\begin{center}
\subfigure[]{\includegraphics[width=6.0cm]{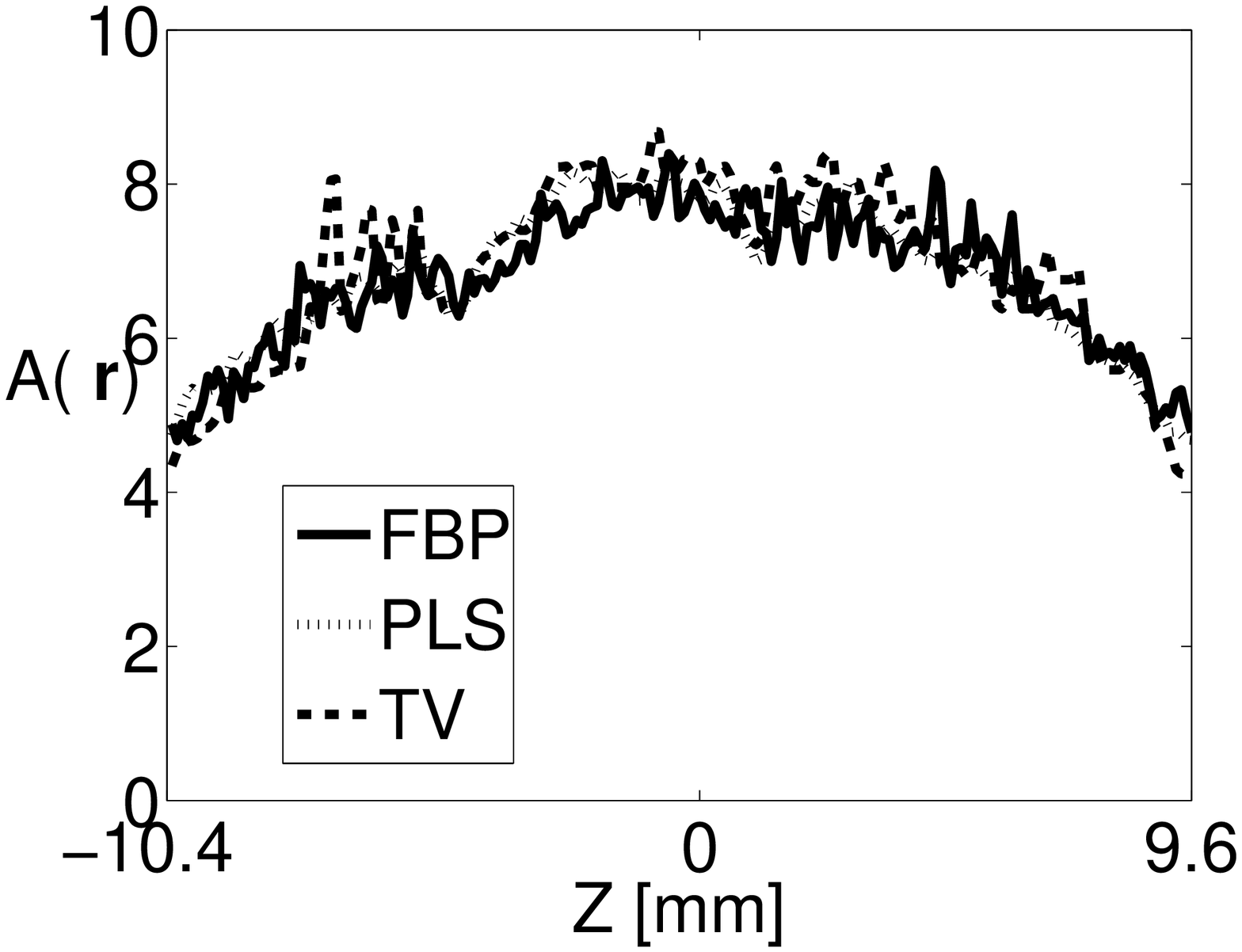}}
\hskip 1.0cm
\subfigure[]
{\includegraphics[width=6.0cm]{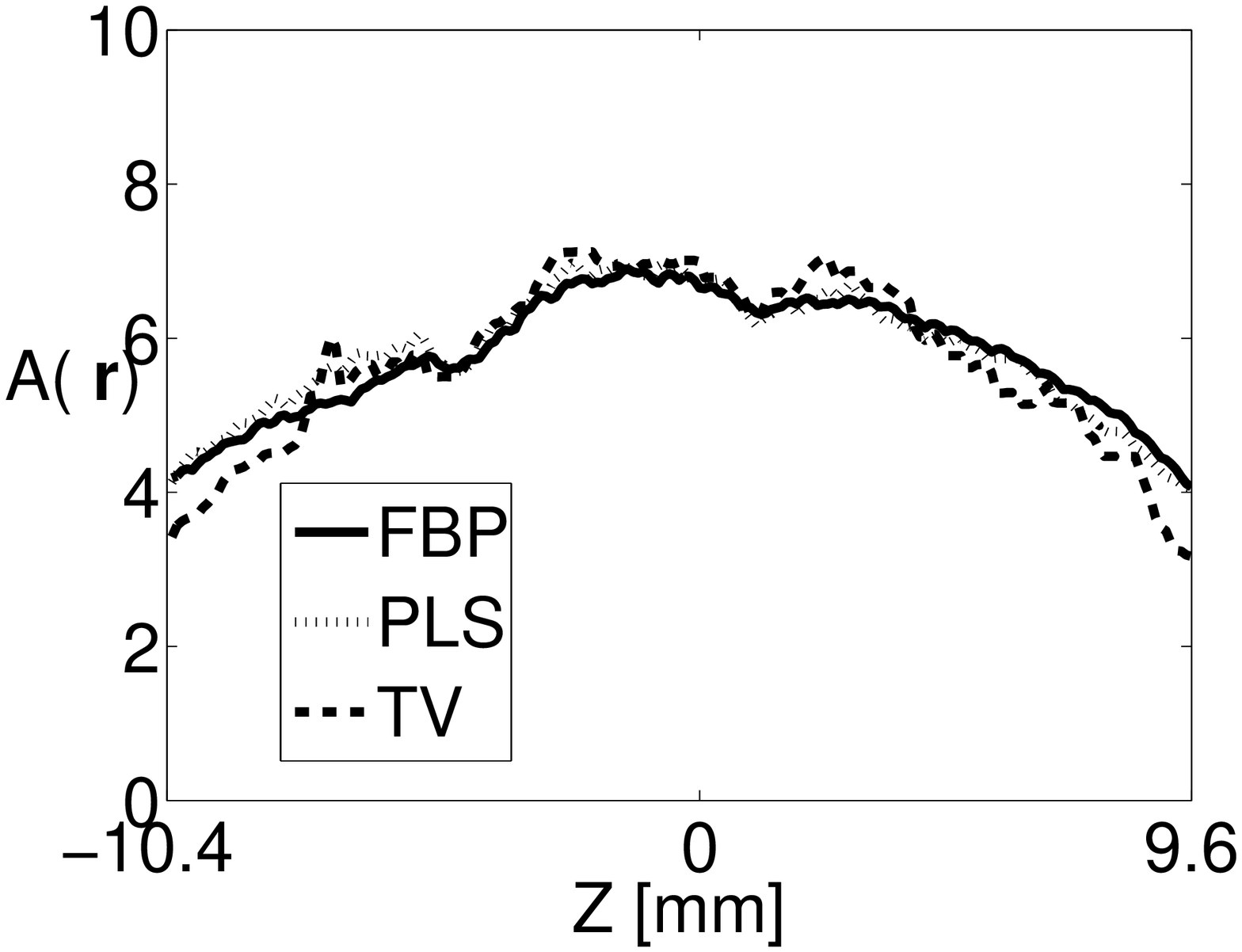}}
\end{center}
\caption{
Image profiles along the z-axis through the center of Tube-A extracted from 
images reconstructed by use of
(a) the FBP algorithm with $f_c=10$-MHz from the $720$-view data (solid line)
the PLS-Q algorithm with $\alpha=1.0\times 10^3$ from the $144$-view data (dotted line), 
and 
the PLS-TV algorithm with $\beta=0.05$ from the $144$-view data (dashed line).
Subfigure (b) shows the corresponding profiles for the case where each 
algorithm employed stronger regularization specified by the parameters
$f_c=5$-MHz, $\alpha=5.0\times 10^3$, and $\beta=0.1$, respectively. 
\label{fig:SixTubeAprofiles}
}
\end{figure}

{\bf Tradeoff between signal intensity and noise level of reconstructed images:}
The image intensities in tube-A are plotted as a function of z, as shown 
in \Fref{fig:SixTubeAprofiles}
where the location of tube-A is indicated in the 2D image slices as shown in 
\Fref{fig:SixTubeDenseSample}-\Fref{fig:SixTubeTVReg}. 
The profiles corresponding to the FBP algorithm were extracted from images reconstructed from 
the 720-view data set while the profiles corresponding to iterative algorithms 
were extracted from images reconstructed from the $144$-view data set. 
For each reconstruction algorithm, two profiles are plotted that correspond to 
moderate and strong regularization as displayed in \Fref{fig:SixTubeAprofiles}-(a) and (b) 
respectively. 
As expected, the quantitative values are smaller when the algorithms are heavily regularized. 
More importantly, images reconstructed by use of iterative image reconstruction algorithms 
quantitatively match with those reconstructed by use of FBP algorithm from densely sampled data. 
In addition, the signal intensities vary gradually along z-axis because the 
laser was illuminated from the side resulting a higher energy distribution near the center 
of z-axis. 
These plots demonstrate the effectiveness of PLS-TV algorithm when the object is not piecewise constant.

\begin{figure}[ht]
\centering
\subfigure[]{\includegraphics[width=7.0cm]{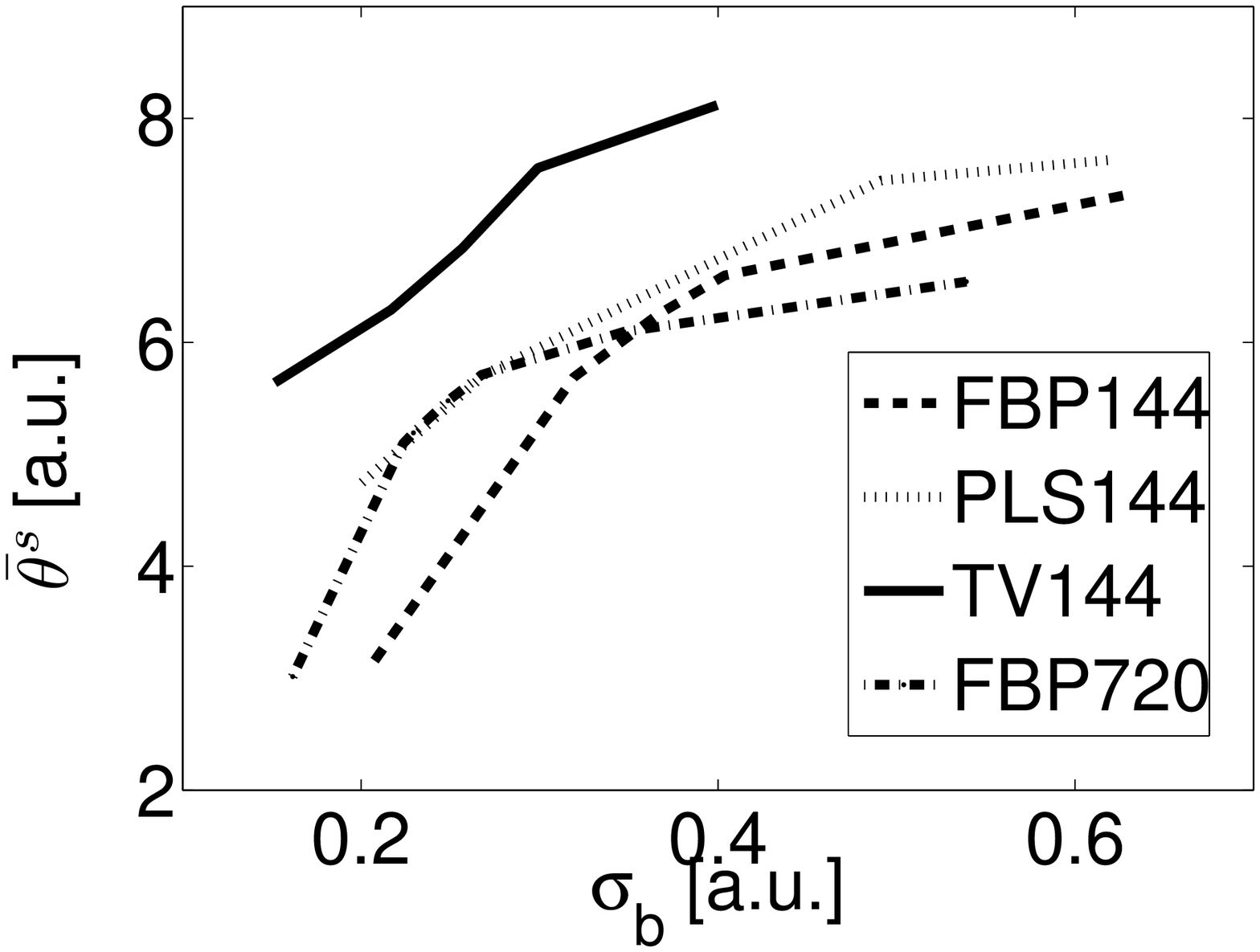}}
\hskip 0.5cm
\subfigure[]{\includegraphics[width=7.0cm]{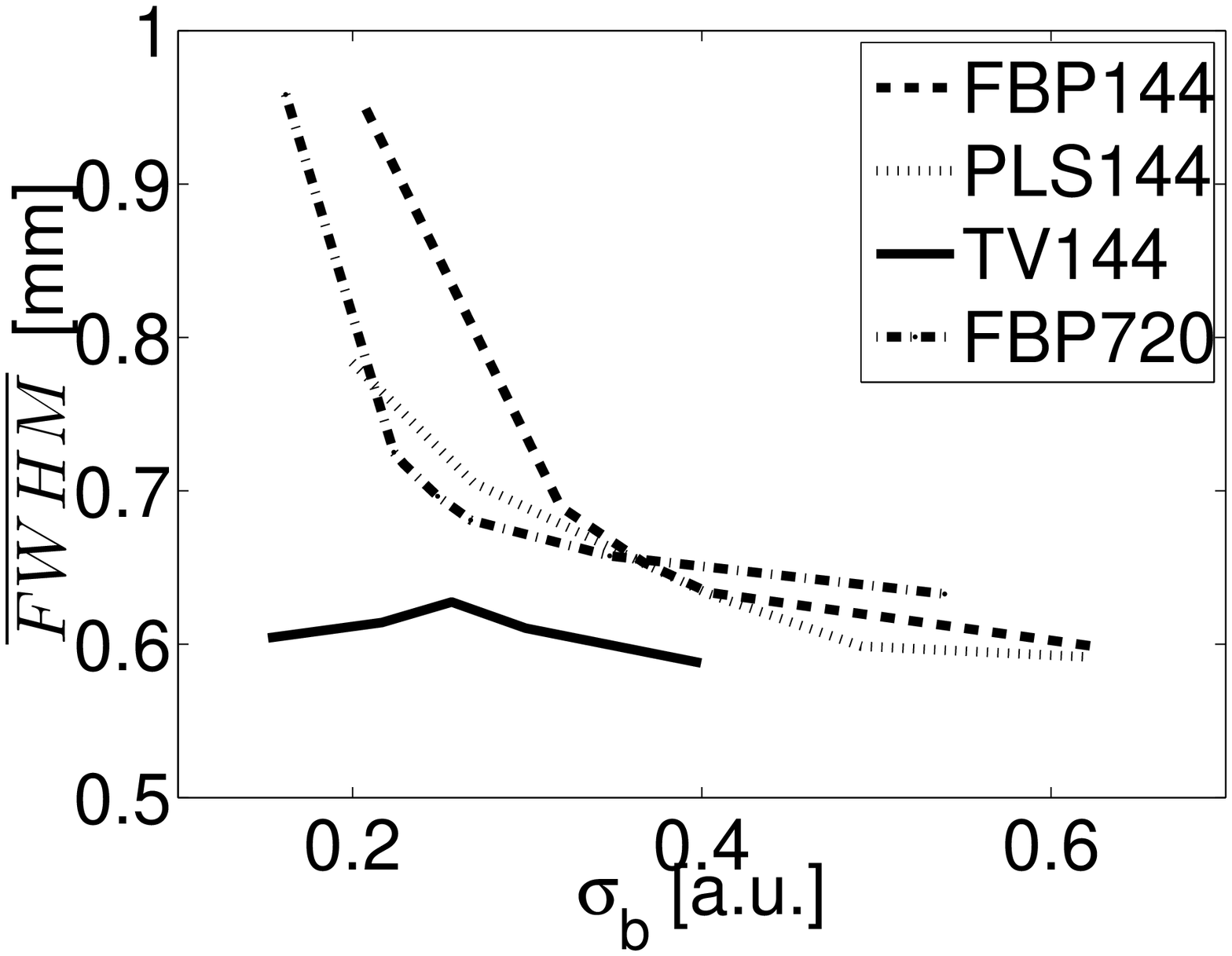}}
\caption{
(a) Signal intensity vs. standard deviation curves and 
(b) image resolution vs. standard deviation curves 
for the images reconstructed by use of the FBP algorithm from the $144$-view data (FBP144), 
the PLS-Q algorithm from the $144$-view data (PLS144), 
the PLS-TV algorithm from the $144$-view data (TV144), 
and 
the FBP algorithm from the $720$-view data (FBP720).
\label{fig:SixTubeTradeoff}
}
\end{figure}

From the same data sets, the signal intensities are plotted against the image standard 
deviations in \Fref{fig:SixTubeTradeoff}-(a). 
This plot suggests that the images reconstructed by use of the PLS-TV algorithm 
have smaller standard deviation while sharing the same bias as those of 
images reconstructed by use of the FBP and the PLS-Q algorithms 
because the same signal intensity indicates the same bias. 
Note that these curves were obtained from densely sampled data. 
Visual inspections suggest the systematic artifacts contribute more to the background 
standard deviation measure than does the measurement random noise. 
Hence, to be more precise, this plot demonstrates the PLS-TV algorithm outperforms 
the FBP and the PLS-Q algorithms in the sense of balancing the tradeoff between 
bias and standard deviation when the signal is present in a uniform background. 

{\bf Tradeoff between resolution and noise level of reconstructed images:}
The plots of resolution ($\overline{\rm FWHM}$) against background standard 
deviation measure ($\sigma_b$) are shown in \Fref{fig:SixTubeTradeoff}-(b). 
Clearly, the spatial resolution of the images reconstructed by use of the PLS-TV algorithm 
is higher than that of the images reconstructed by the FBP and the PLS-Q algorithms
while the images having the same background standard deviation.  
In addition, the curve corresponding to the PLS-TV algorithm is flatter than those corresponding to 
the FBP and PLS-Q algorithm, suggesting that TV regularization mitigates image noise 
with minimal sacrifice in spatial resolution. 
This observation is consistent with our earlier visual inspections of the reconstructed 
images. 
It is also interesting to note that the curve corresponding to the PLS-Q algorithm
intersects the one corresponding to the FBP algorithm, indicating that 
conventional iterative reconstruction algorithms may not always outperform the FBP algorithm. 

\begin{figure}[ht]
\begin{center}
\subfigure[]
{\includegraphics[width=5.1cm]{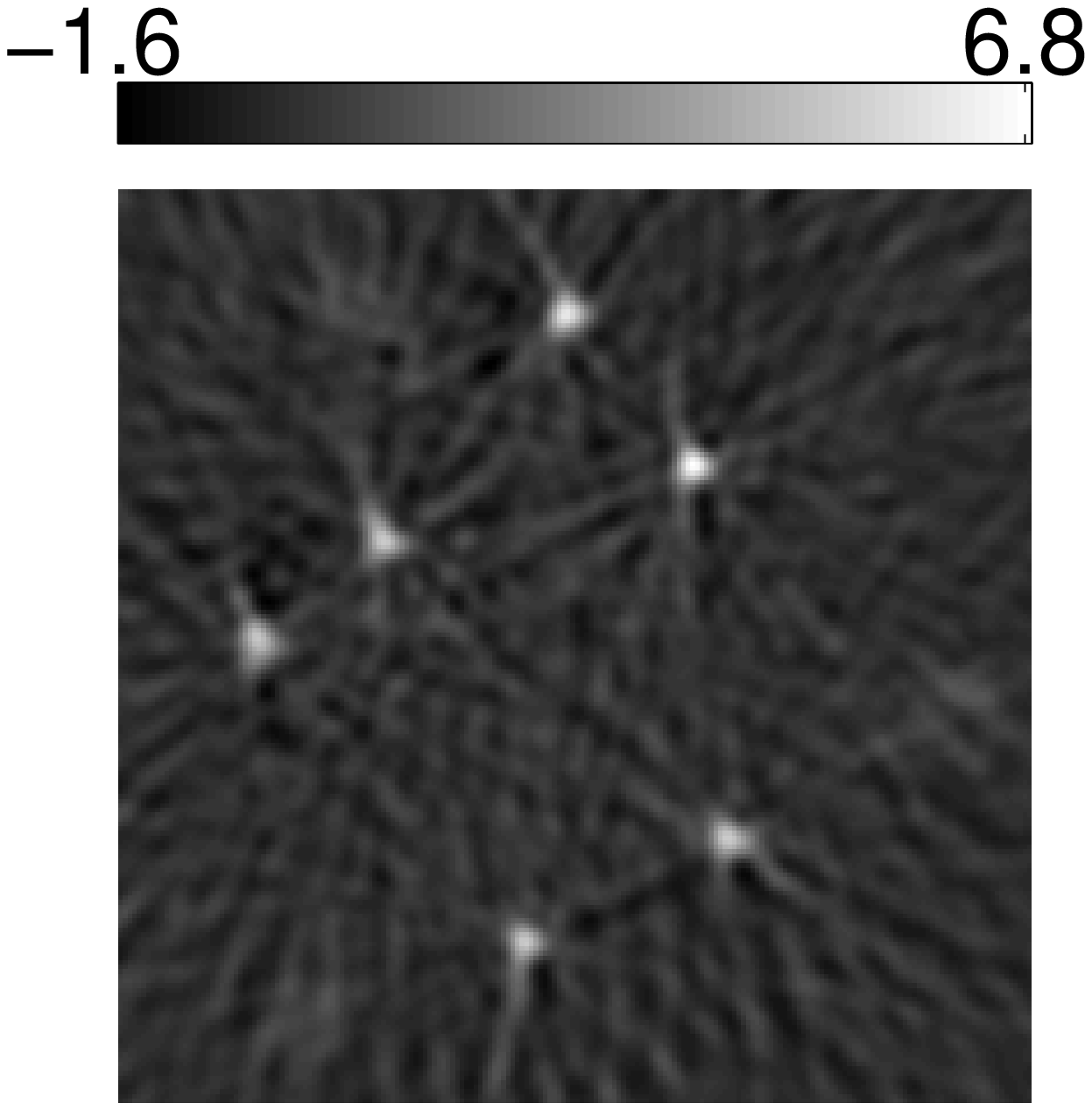}}
\subfigure[]
{\includegraphics[width=5.1cm]{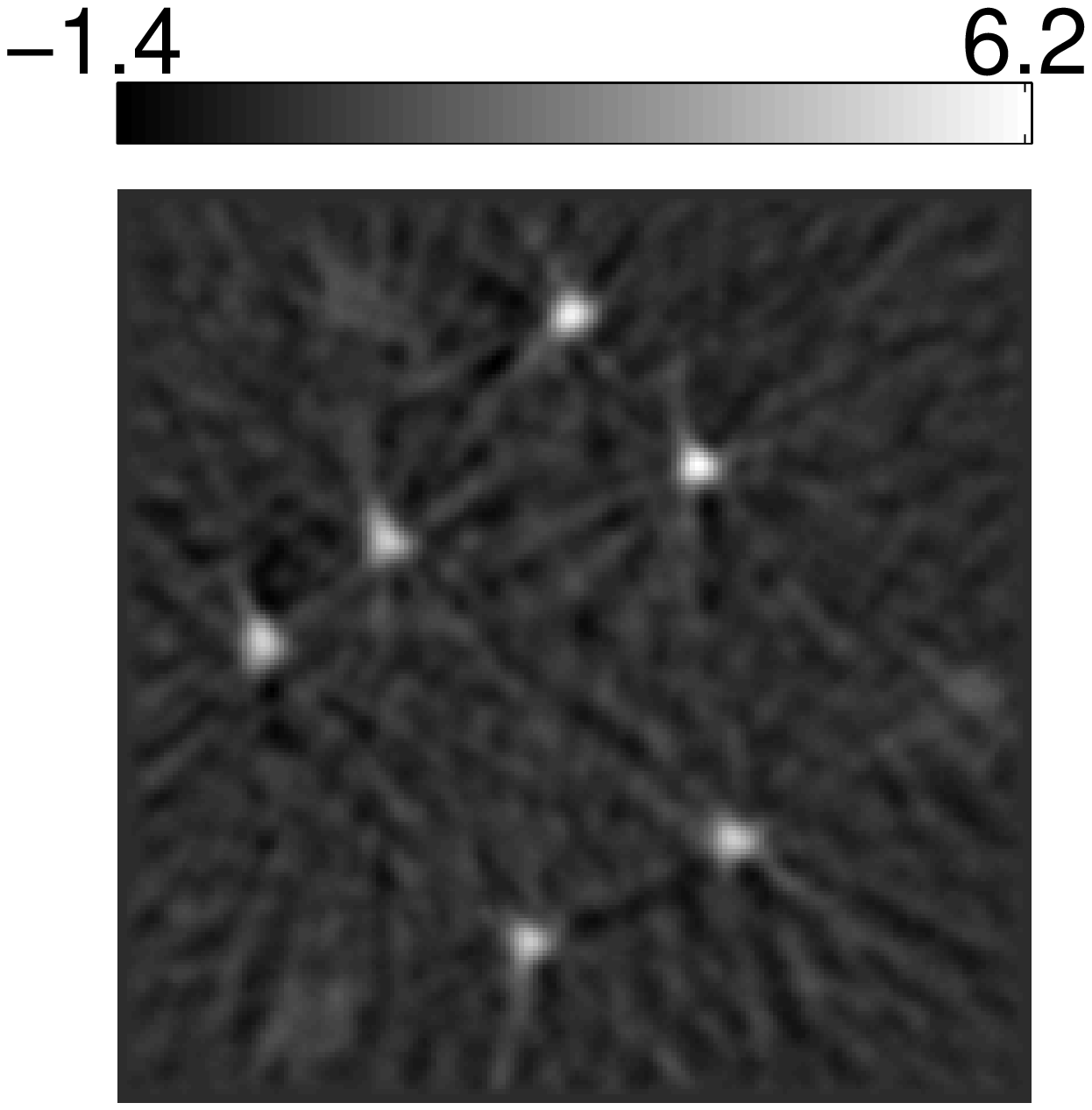}}
\subfigure[]
{\includegraphics[width=5.1cm]{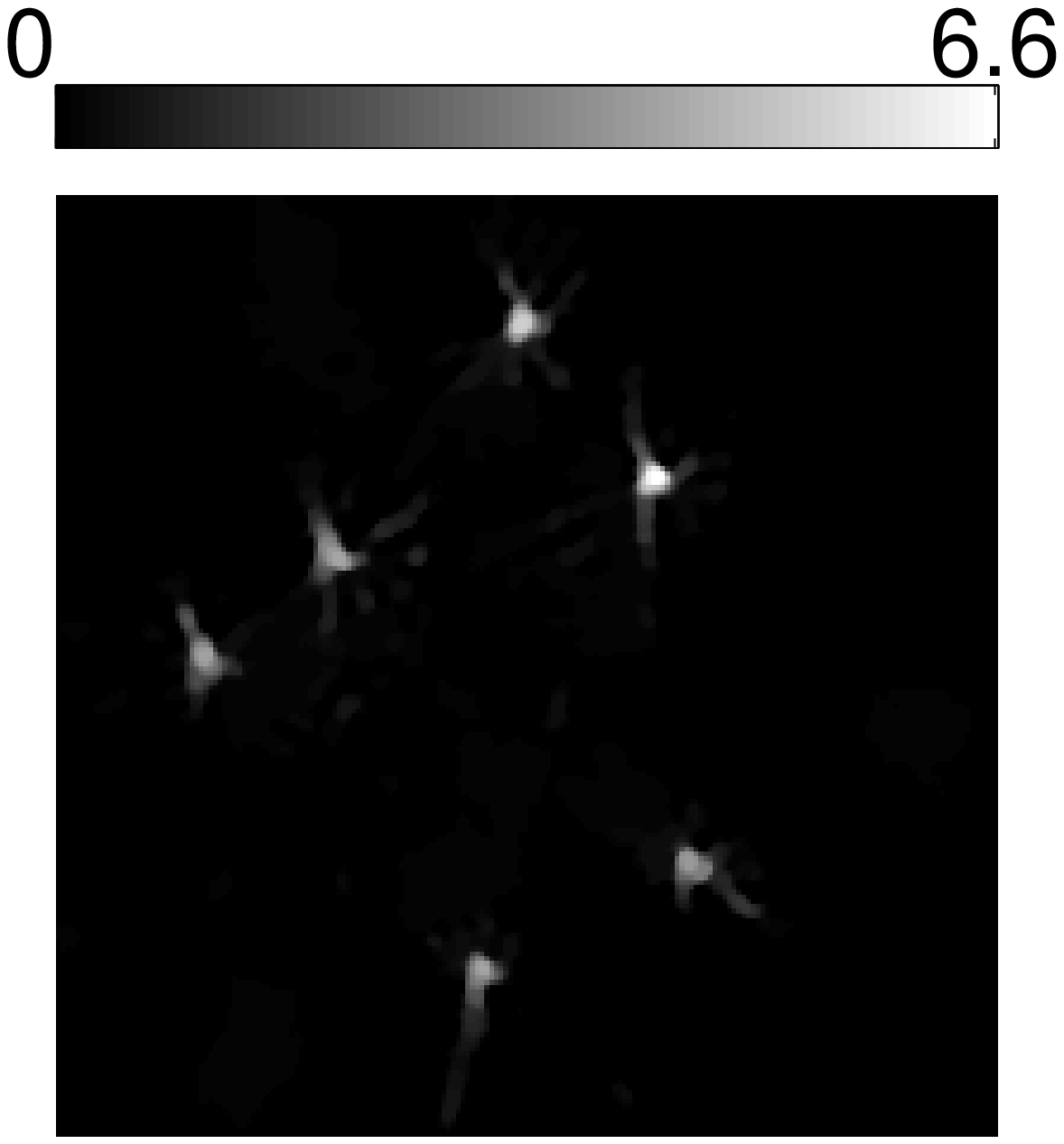}}
\end{center}
\caption{
Slices corresponding to the plane $z=-2.0$-mm through the 3D images of the six-tube phantom 
reconstructed from the $72$-view data by use of 
(a) the FBP algorithm with $f_c=3.7$-MHz; 
(b) the PLS-Q algorithm with $\alpha=1.0\times 10^{3}$; 
and
(c) the PLS-TV algorithm with $\beta=0.07$. 
All images are of size $19.8\times 19.8$-mm$^2$.
The ranges of the grayscale windows were determined by the minimum and the maximum values in each image.
\label{fig:SixTubeView72}
}
\end{figure}

\begin{figure}[ht]
\begin{center}
\subfigure[]
{\includegraphics[width=5.1cm]{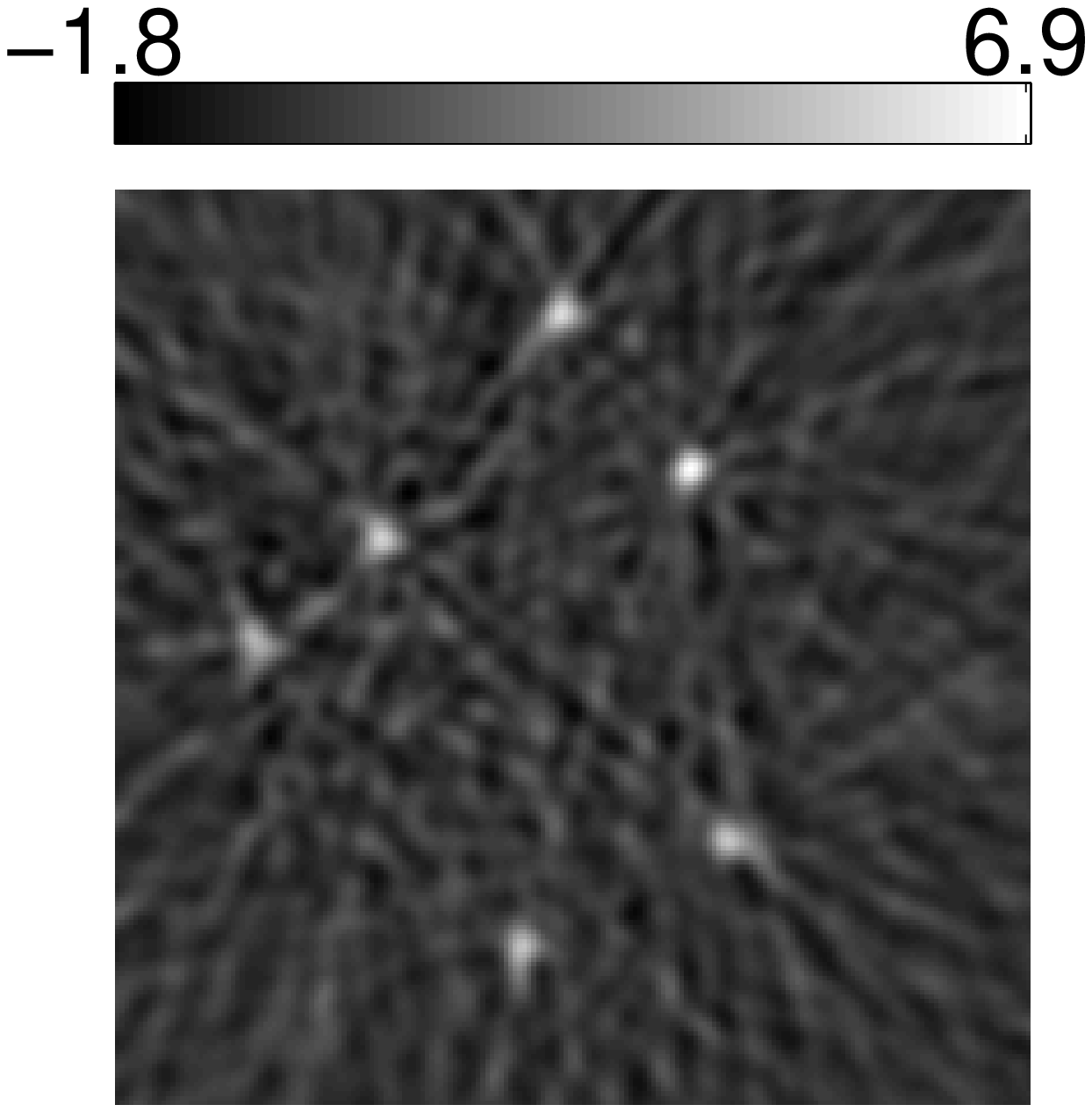}}
\subfigure[]
{\includegraphics[width=5.1cm]{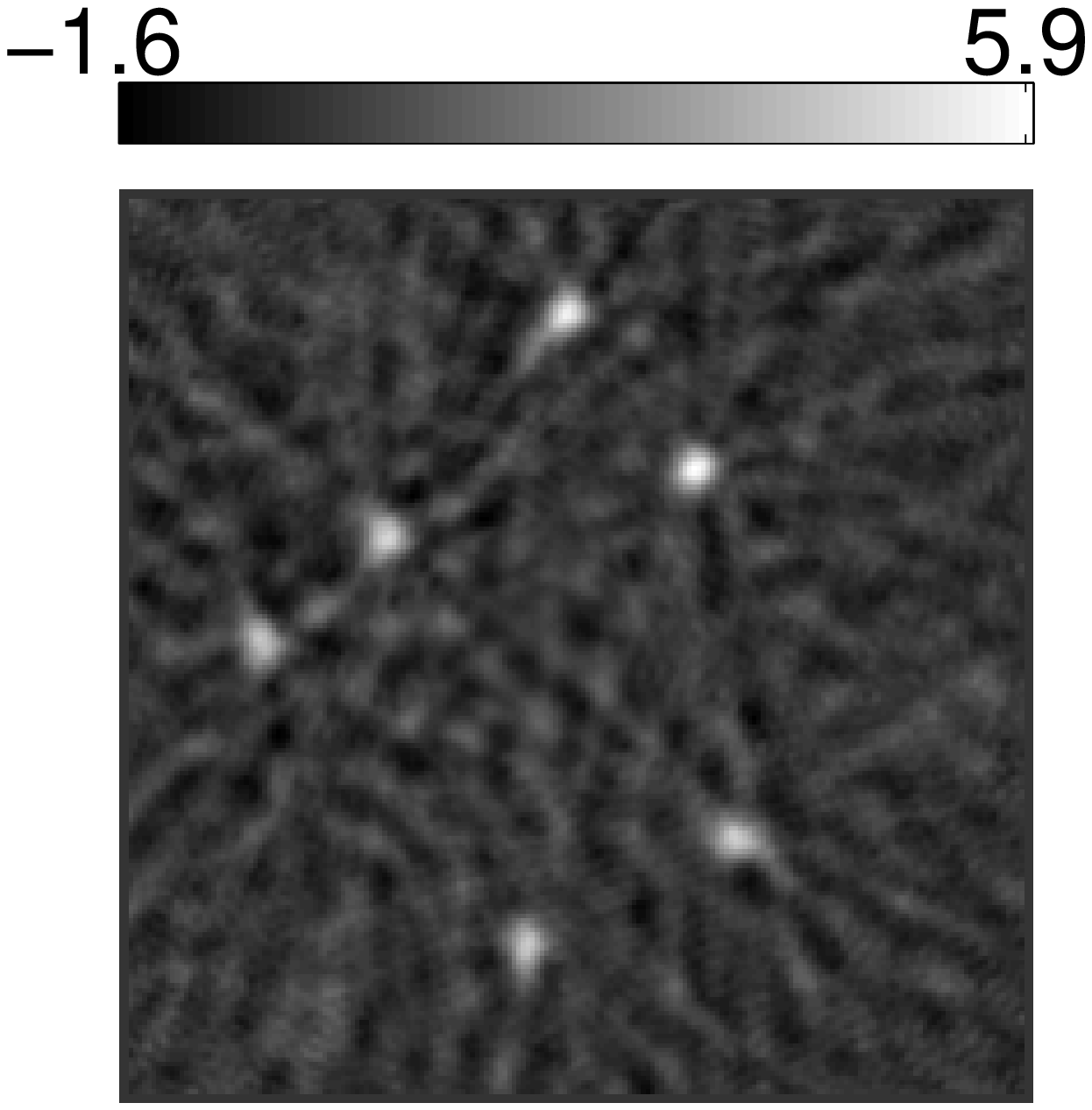}}
\subfigure[]
{\includegraphics[width=5.1cm]{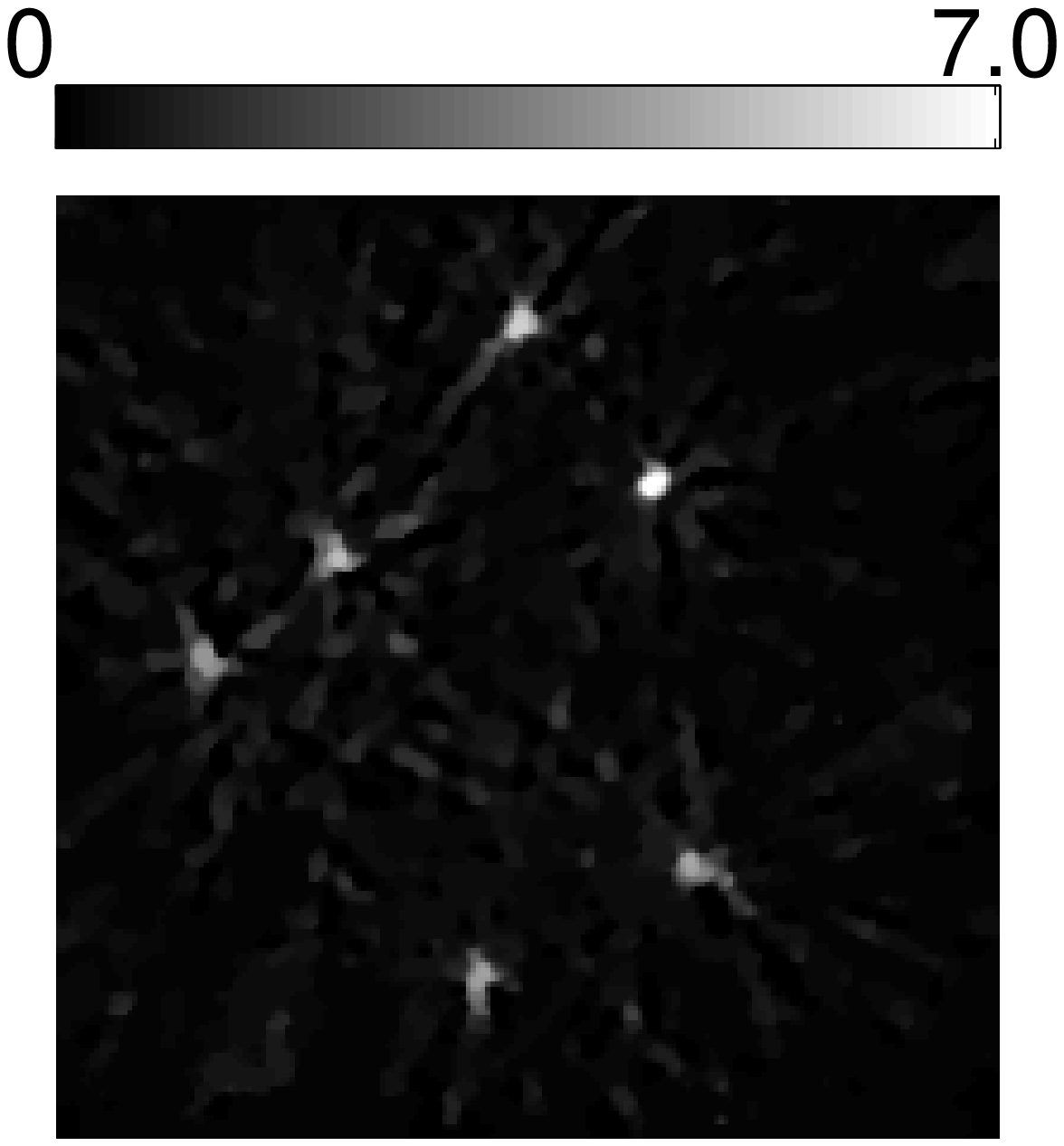}}
\end{center}
\caption{
Slices corresponding to the plane  $z=-2.0$-mm through the 3D images of the six-tube phantom
reconstructed from the $36$-view data by use of 
(a) the FBP algorithm with $f_c=3.3$-MHz; 
(b) the PLS-Q algorithm with $\alpha=7.0$; 
and
(c) the PLS-TV algorithm with $\beta=0.02$;
\label{fig:SixTubeView36}
All images are of size $19.8\times 19.8$-mm$^2$.
The ranges of the grayscale windows were determined by the minimum and the maximum values in each image.
}
\end{figure}

{\bf Reconstructed images from sparsely-sampled data sets:}
Images reconstructed from the 72-view data set and the 36-view data set 
are displayed in \Fref{fig:SixTubeView72} and \Fref{fig:SixTubeView36} 
respectively. 
The regularization parameters were selected such that the quantitative 
values of the images are within the similar range. 
As expected, from both data sets, the images reconstructed by use of PLS-TV 
algorithm appear to have higher spatial resolution as well as cleaner 
backgrounds, suggesting the PLS-TV algorithm can effectively mitigate data 
incompleteness in 3D OAT.

\subsection{Mouse whole-body imaging}
\begin{figure}[h]
\begin{center}
\subfigure[]{\includegraphics[width=5cm, trim=6cm 1.0cm 6cm 4cm, clip]{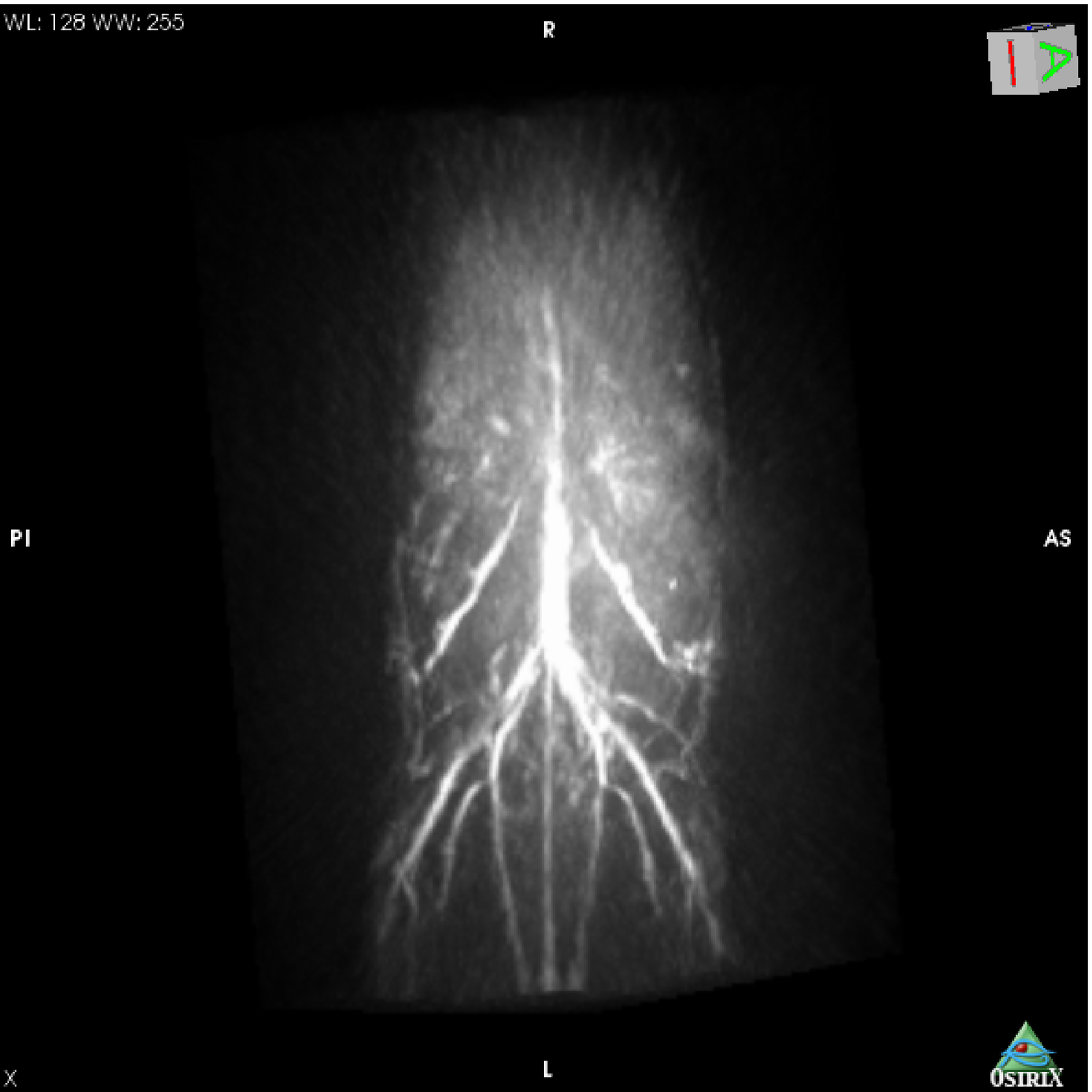}}
\subfigure[]
{\includegraphics[width=5.0cm,trim=6cm 1.0cm 6cm 4cm, clip]{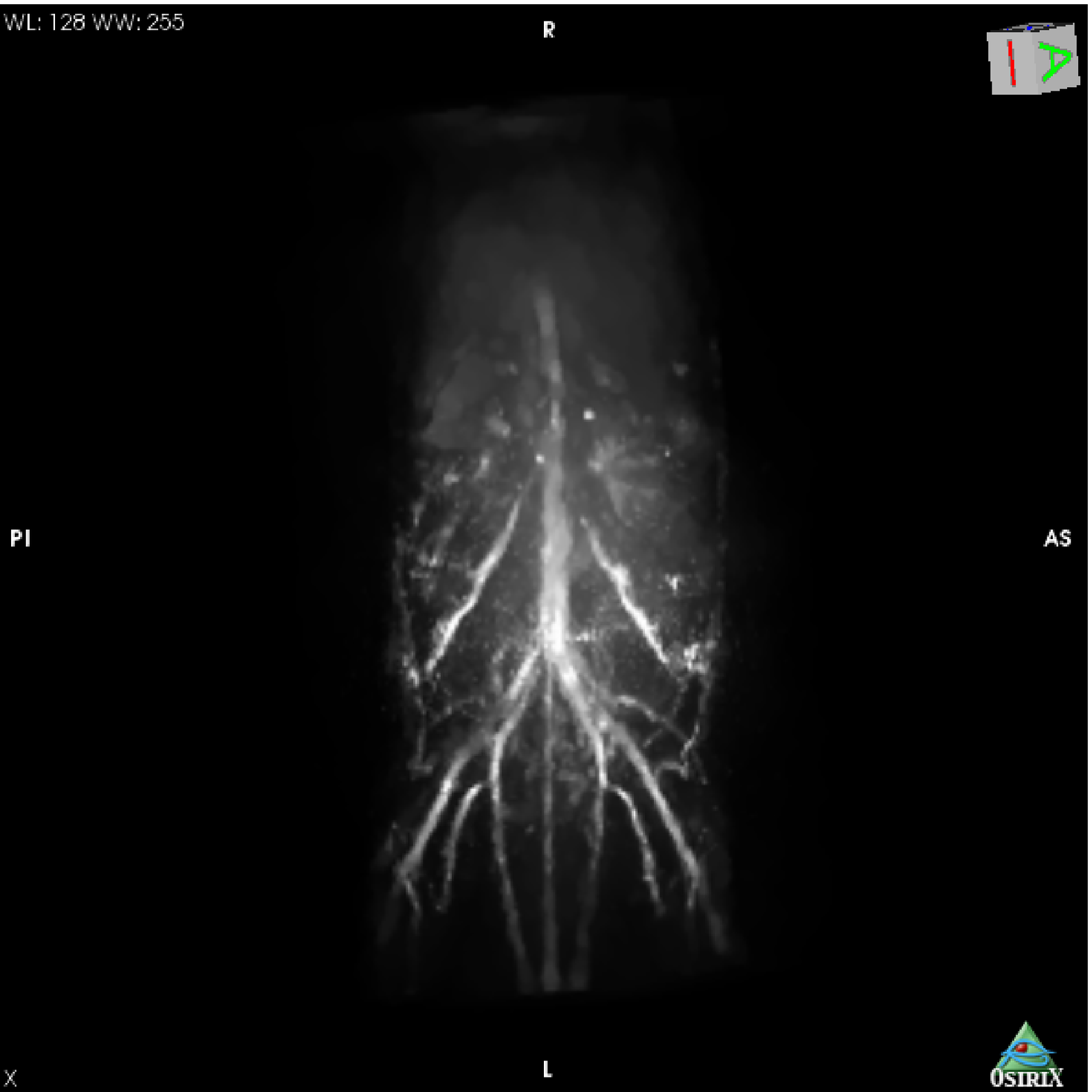}}
\end{center}
\caption{
MIP renderings of the 3D images of the mouse body 
reconstructed from the $180$-view data 
by use of 
(a) the FBP algorithm with $f_c=5$-MHz; 
and
(b) the PLS-TV algorithm with $\beta=0.05$;
The grayscale window is [0,12.0]. 
(QuickTime)
\label{fig:mouseV1803D}
}
\end{figure}

{\bf Visual inspection of reconstructed images from densely-sampled data sets:}
From the 180-view data set, the MIP images corresponding to the FBP and the PLS-TV algorithms 
appear to be very similar as shown in \Fref{fig:mouseV1803D}.
In contrast to the images of the six-tube phantom that have a uniform background, 
the mouse whole-body images have a diffuse background. 
The diffuse background is due to the measurement random noise as well as the 
capillaries that are beyond the resolution limit of the imaging 
system \cite{petermice:jbo}, thus carrying little information regarding the object. 
In general, the images reconstructed by the PLS-TV algorithm appear to have a cleaner 
background while revealing a sharper appearing body vascular tree. 
Besides, the left kidney in the images reconstructed by use of the PLS-TV algorithm 
appears to have a higher contrast than the image reconstructed by use of the FBP algorithm. 
In addition, comparing with the images reconstructed by use of direct backprojection 
from the raw data, (see figure 6 in \cite{petermice:spie}), 
both our algorithms appear to improve the continuity of blood vessels. 
We believe this is because our algorithms are based on an imaging model that incorporates 
the transducer SIR and EIR. 

\begin{figure}[ht]
\begin{center}
\subfigure[]{\includegraphics[width=5.5cm]{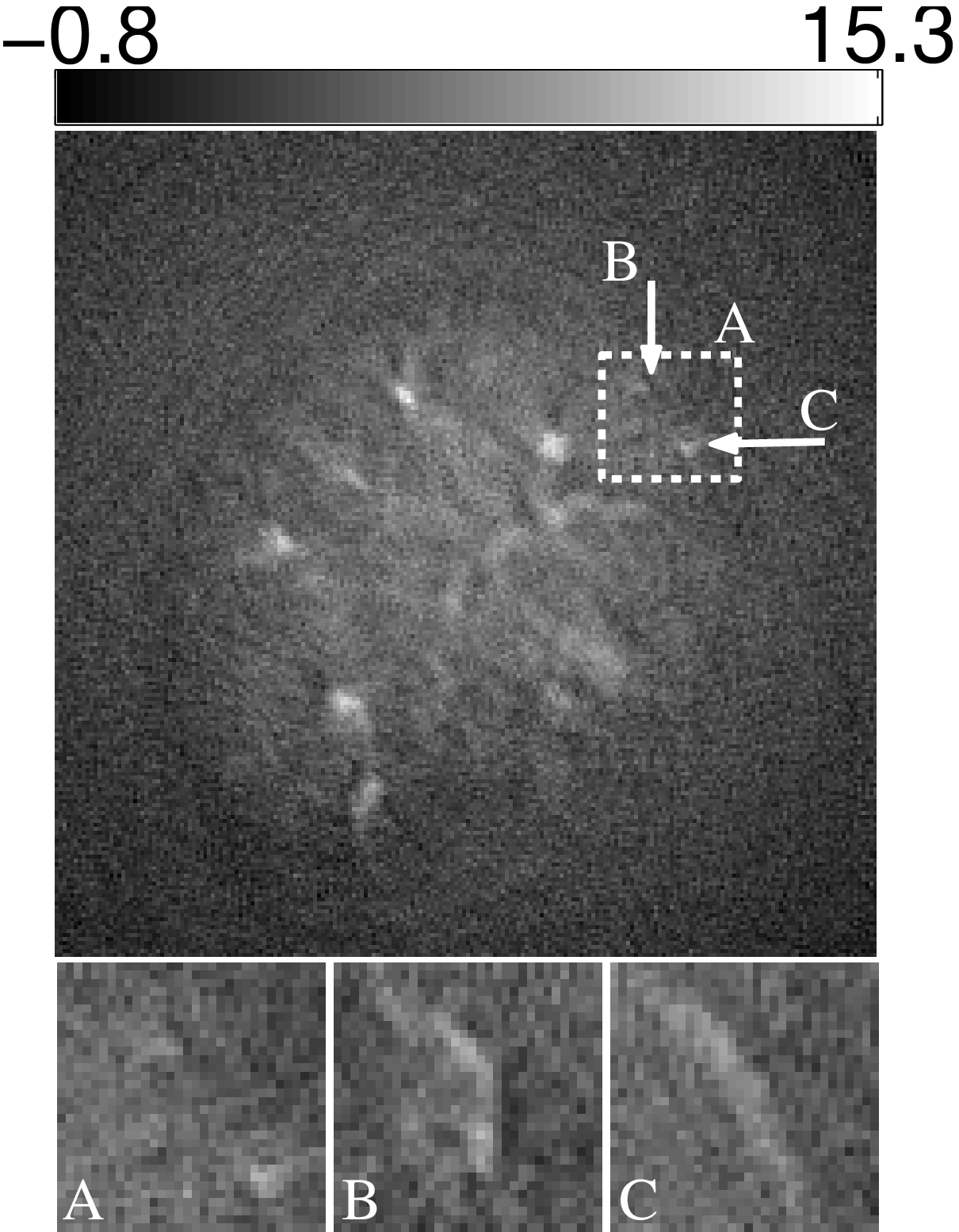}}
\hskip 1.0cm
\subfigure[]{\includegraphics[width=5.5cm]{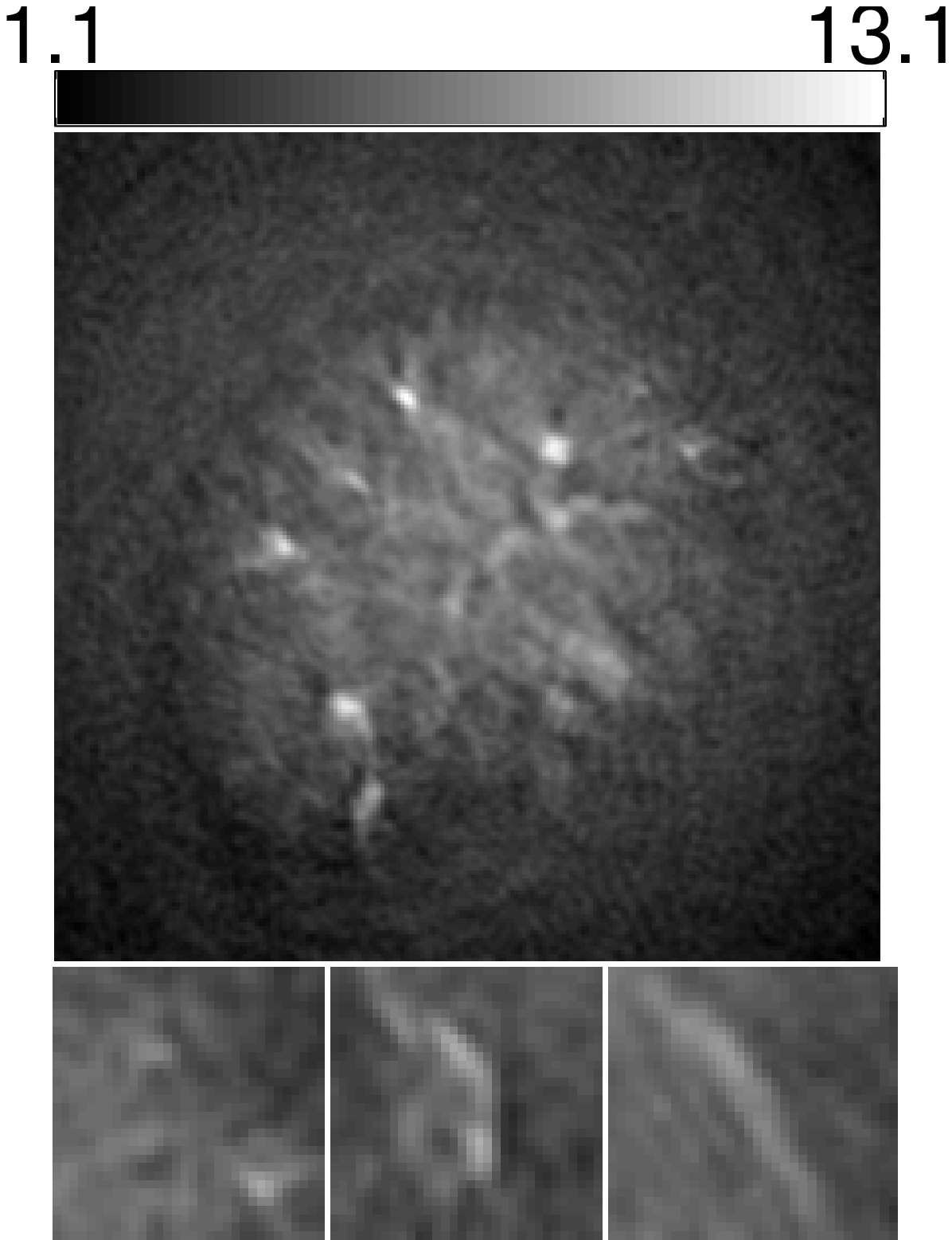}}\\
\subfigure[]{\includegraphics[width=5.5cm]{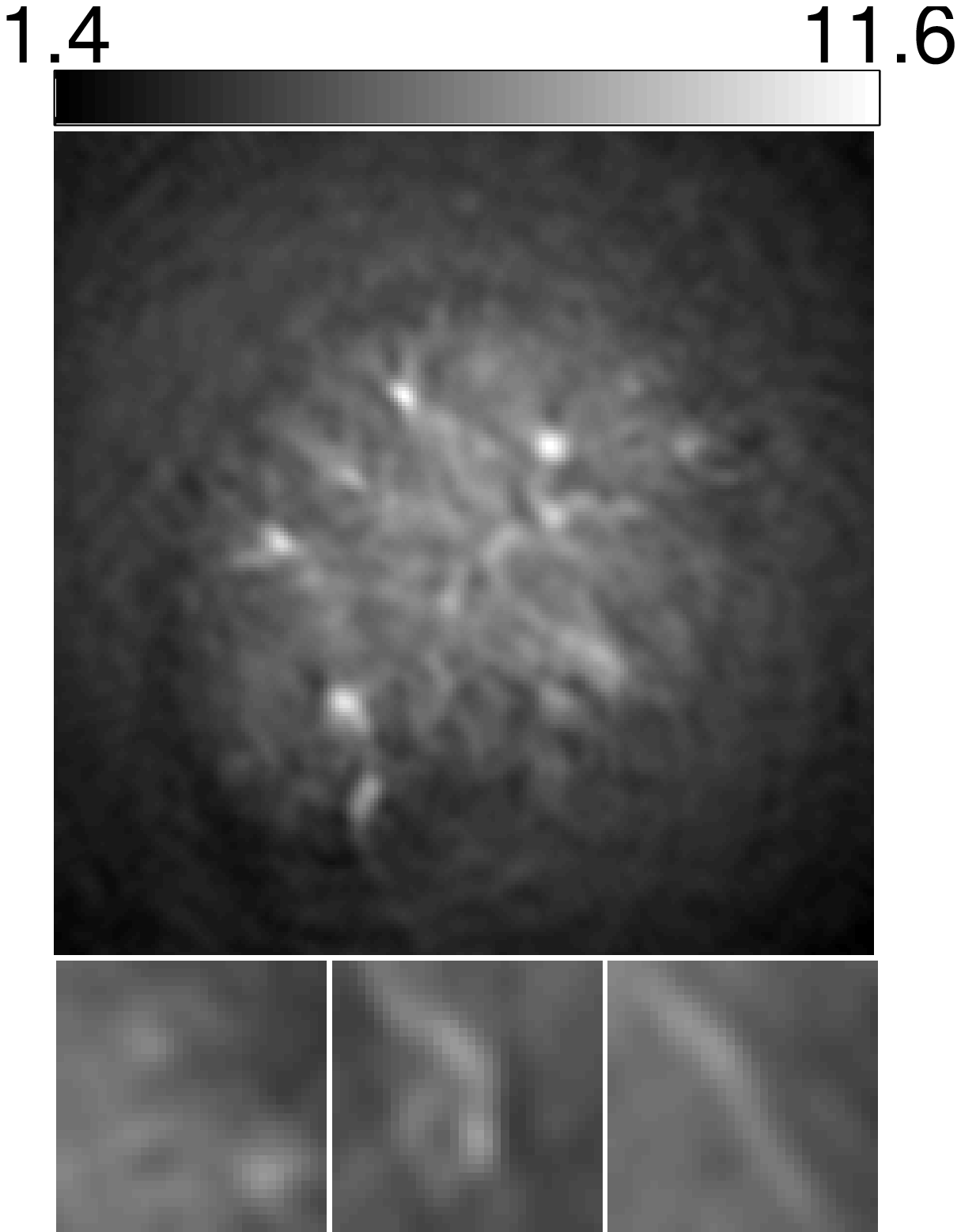}}
\hskip 1.0cm
\subfigure[]{\includegraphics[width=5.5cm]{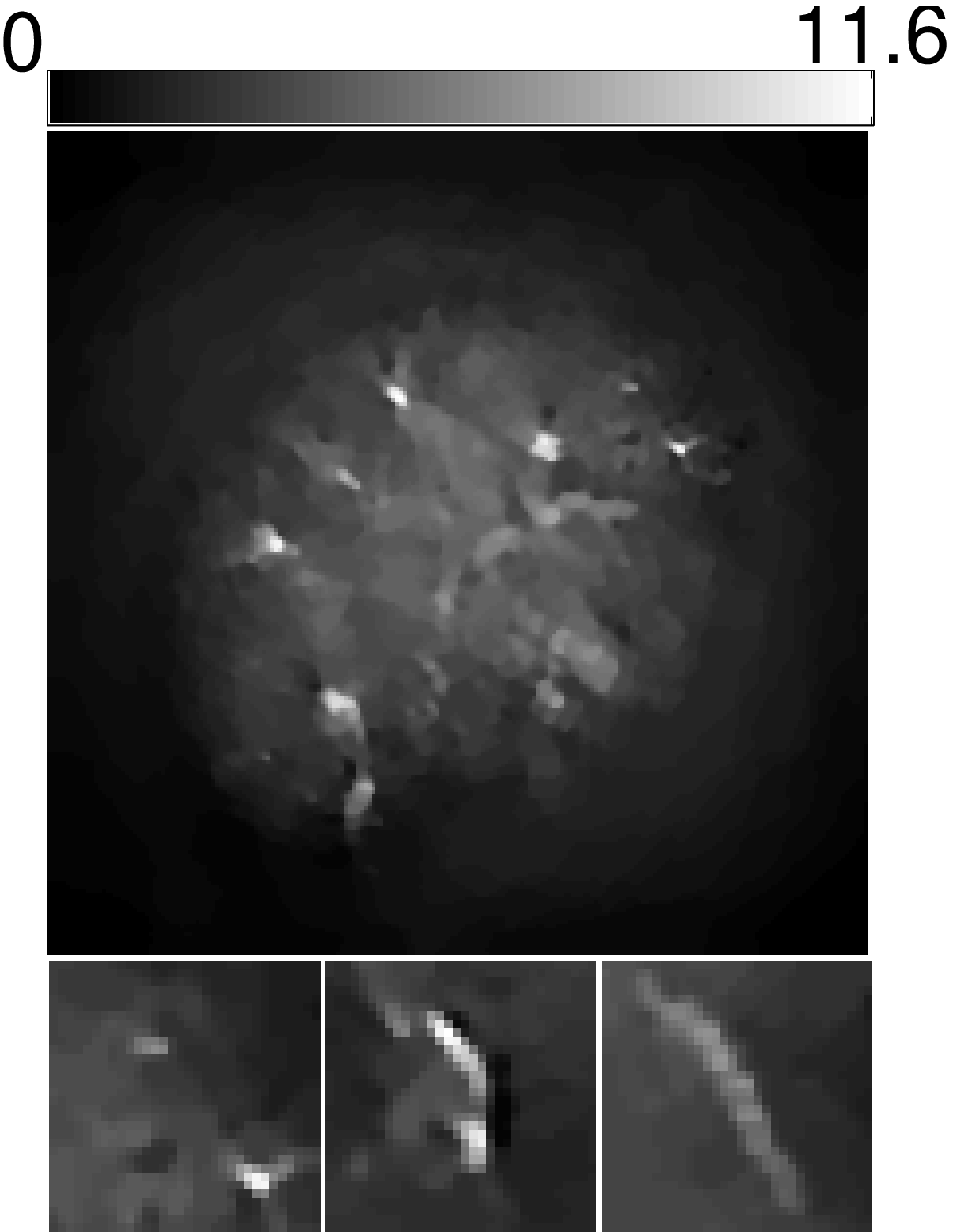}}
\end{center}
\caption{
Slices corresponding to the plane $z=-8.47$-mm through the 3D images of the mouse body
reconstructed from the $180$-view data by use of 
(a) the FBP algorithm with $f_c=8$-MHz;  
(b) the FBP algorithm with $f_c=5$-MHz;
(c) the FBP algorithm with $f_c=3$-MHz;
and
(d) the PLS-TV algorithm with $\beta=0.05$. 
The images are of size $29.4\times 29.4$-mm$^2$. 
The three zoomed-in images correspond to the ROIs of the dashed rectangle A, 
and the images on the orthogonal planes $x=8.47$-mm (Intersection line is along the arrow B), 
and $y=-3.29$-mm (Intersection line is along the arrow C), respectively. 
All zoomed-in images are of size $4.34\times 4.34$-mm$^2$. 
The ranges of the grayscale windows were determined by the minimum and the maximum values in each image. 
\label{fig:mouseV1802D}
}
\end{figure}

\begin{figure}[ht]
\begin{center}
\subfigure[]{\includegraphics[width=5.5cm]{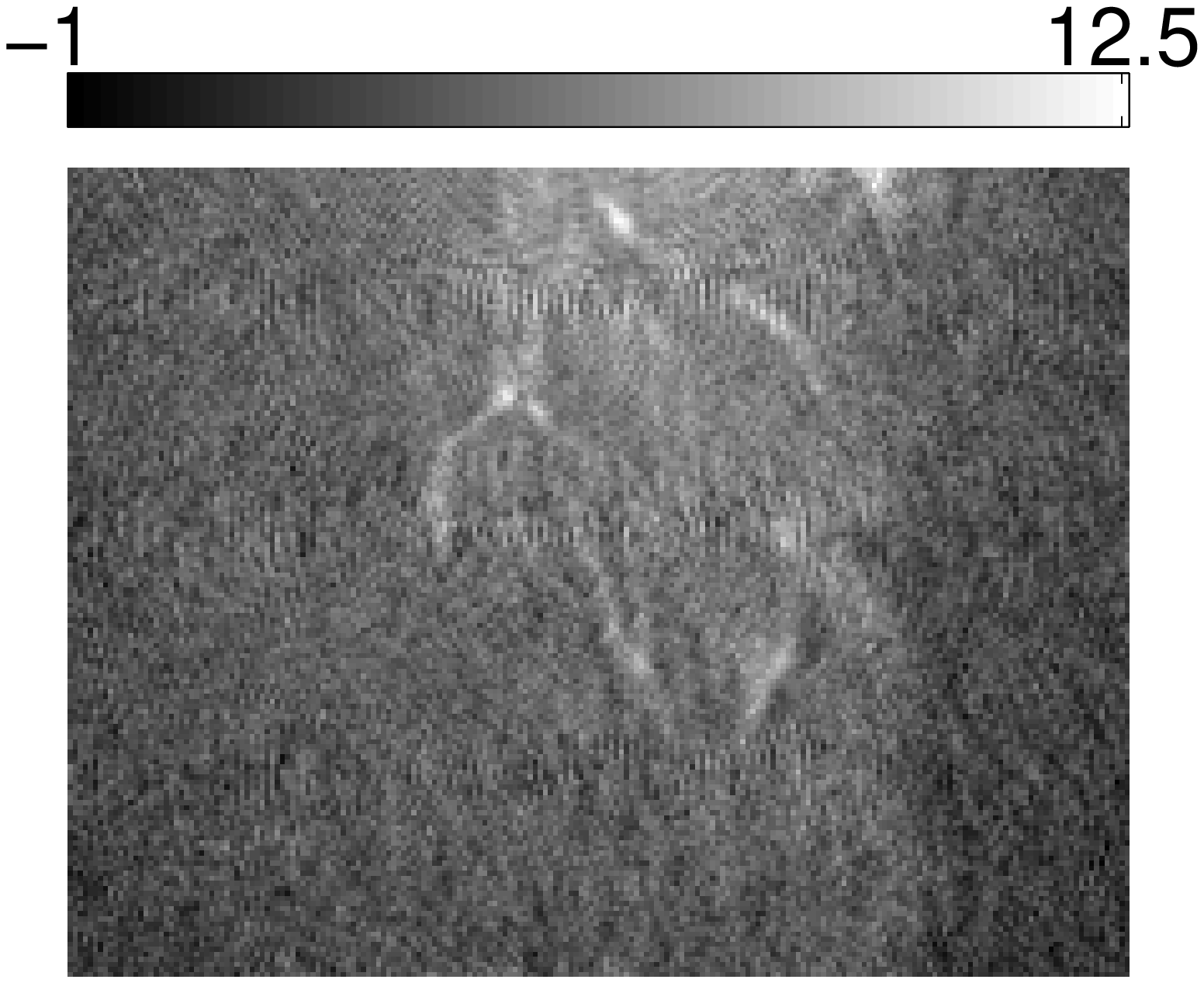}}
\hskip .5cm
\subfigure[]{\includegraphics[width=5.5cm]{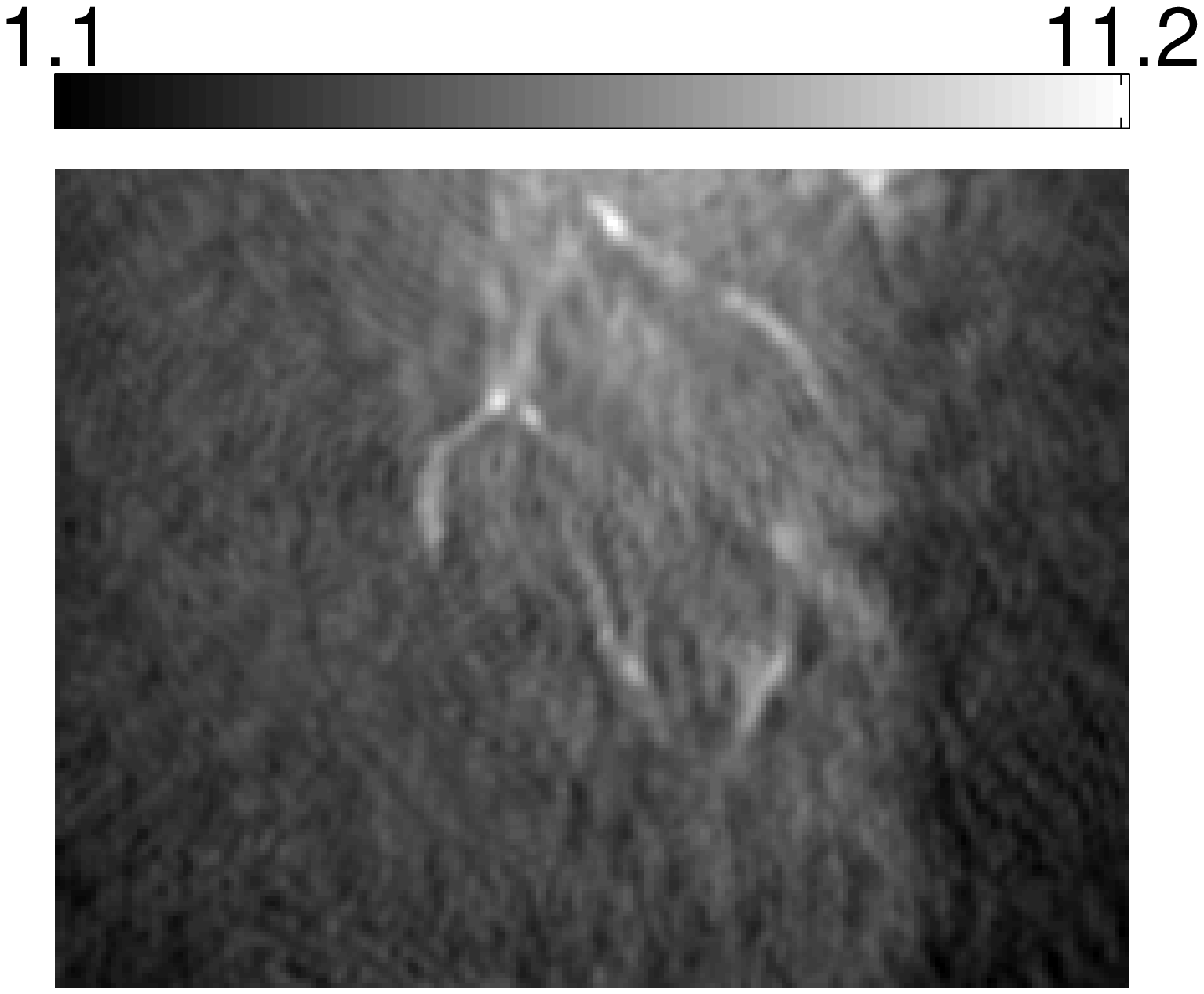}}\\
\subfigure[]{\includegraphics[width=5.5cm]{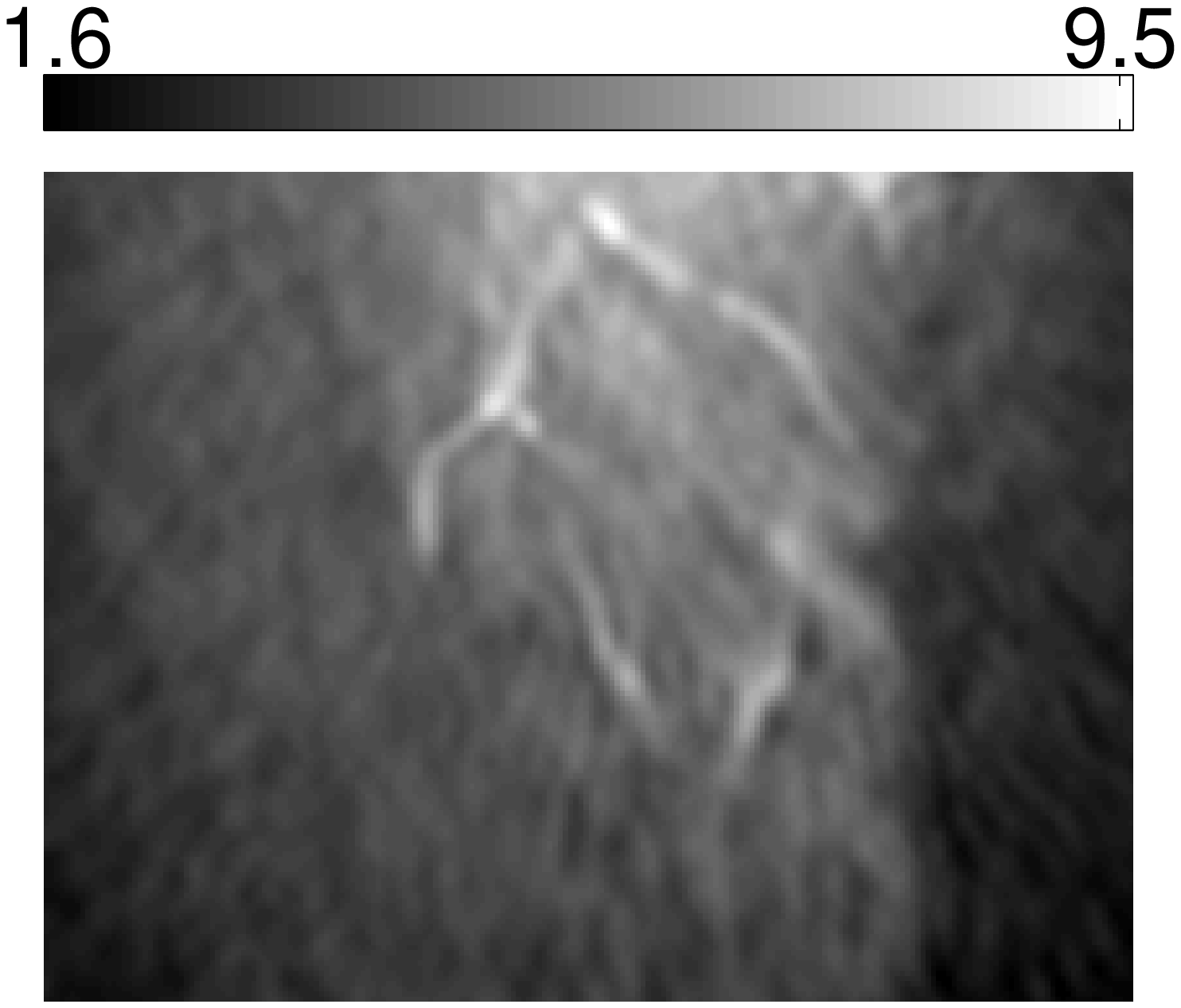}}
\hskip .5cm
\subfigure[]{\includegraphics[width=5.5cm]{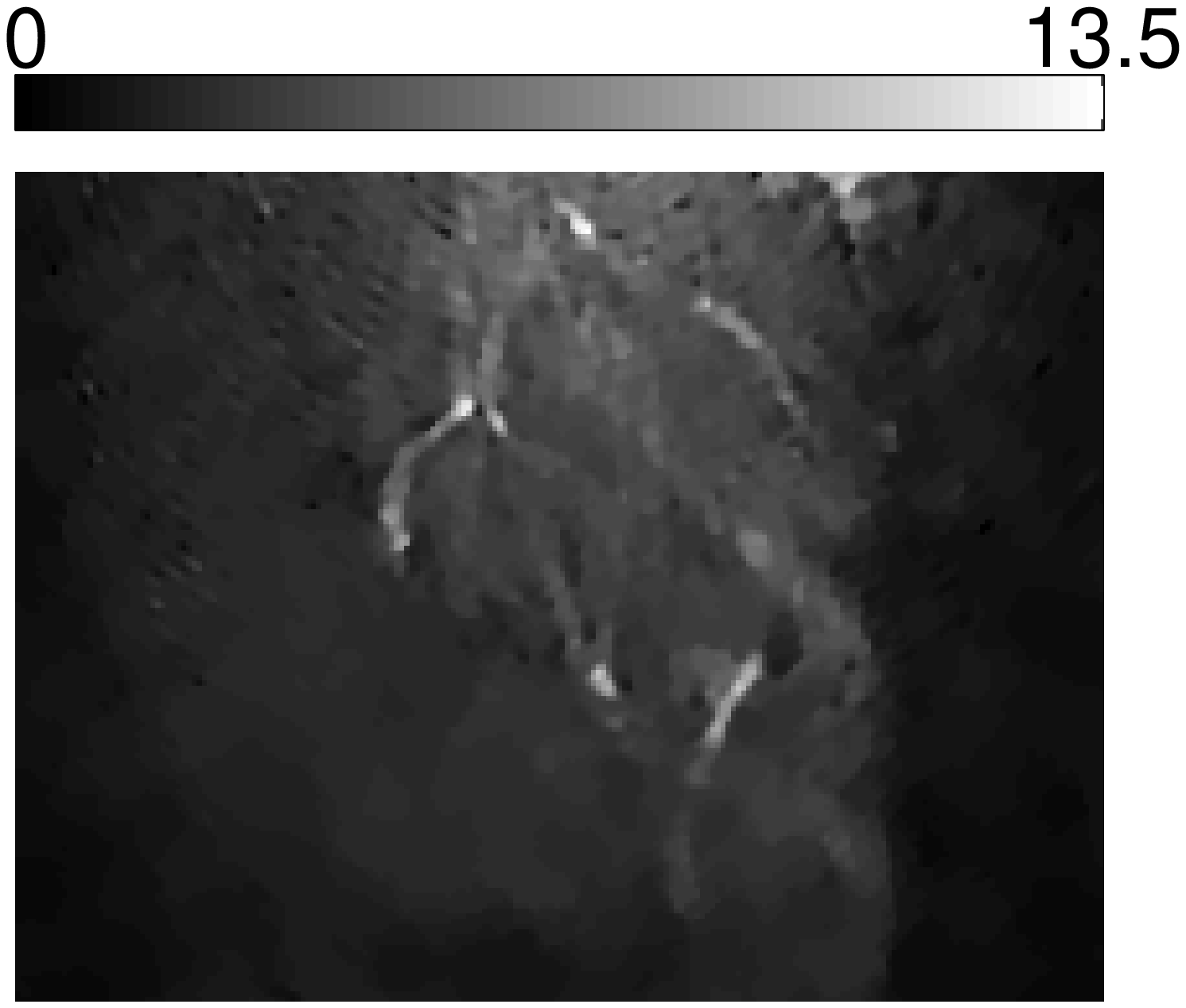}}
\end{center}
\caption{
Slices corresponding to the plane $y=-3.57$-mm through the 3D images of the mouse body 
reconstructed from the $180$-view data by use of 
(a) the FBP algorithm with $f_c=8$-MHz;  
(b) the FBP algorithm with $f_c=5$-MHz;
(c) the FBP algorithm with $f_c=3$-MHz;
and
(d) the PLS-TV algorithm with $\beta=0.05$. 
The images are of size $22.4\times 29.4$-mm$^2$. 
The ranges of the grayscale windows were determined by the minimum and the maximum values in each image. 
\label{fig:mouseV1802DY}
}
\end{figure}

Additional details are revealed in the 2D sectional images as shown in 
\Fref{fig:mouseV1802D} and \Fref{fig:mouseV1802DY}.  
Obviously, the contrast of the blood vessels in the images reconstructed by use of 
the PLS-TV algorithm are higher than those reconstructed by use of the FBP algorithm. 
In particular, the PLS-TV algorithm significantly enhanced the appearance of peripheral blood vessels. 
For example, within the ROI A in \Fref{fig:mouseV1802D}, two 
blood vessels B and C can be detected easily as two bright spots in 
zoomed-in image A. However, the two bright spots have much lower visual contrast in the 
images reconstructed by use the FBP algorithm. 
In addition, as shown in \Fref{fig:mouseV1802DY}, 
the PLS-TV algorithm effectively mitigates the foggy background as well as
noise with minimal sacrifice in image resolution. 
However, none of the images reconstructed by use the FBP algorithm 
has a background as clean as the images reconstructed by the PLS-TV algorithm. 

\begin{figure}[ht]
\begin{center}
\subfigure[]{\includegraphics[width=5.5cm]{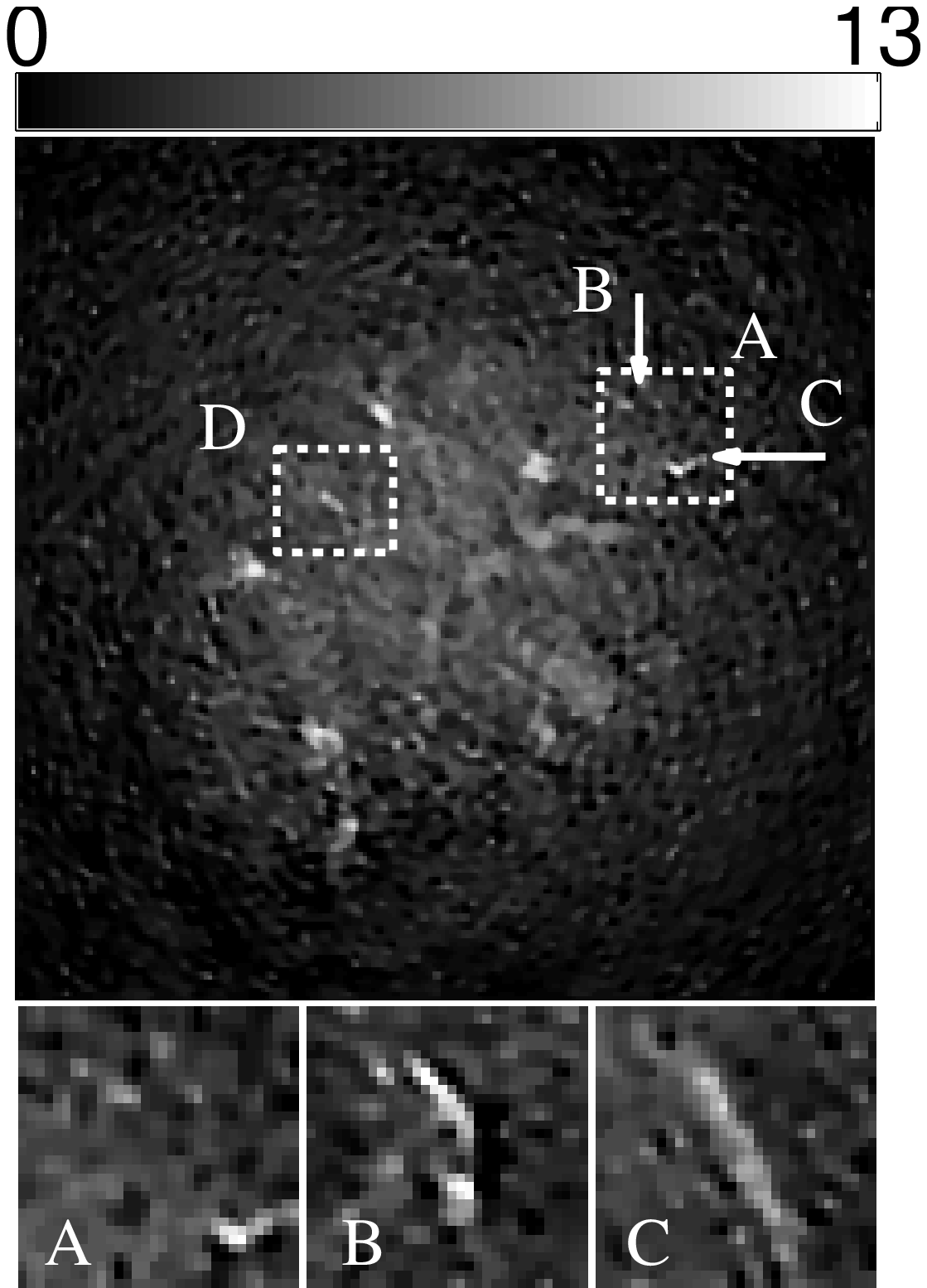}}
\hskip 1cm
\subfigure[]{\includegraphics[width=5.5cm]{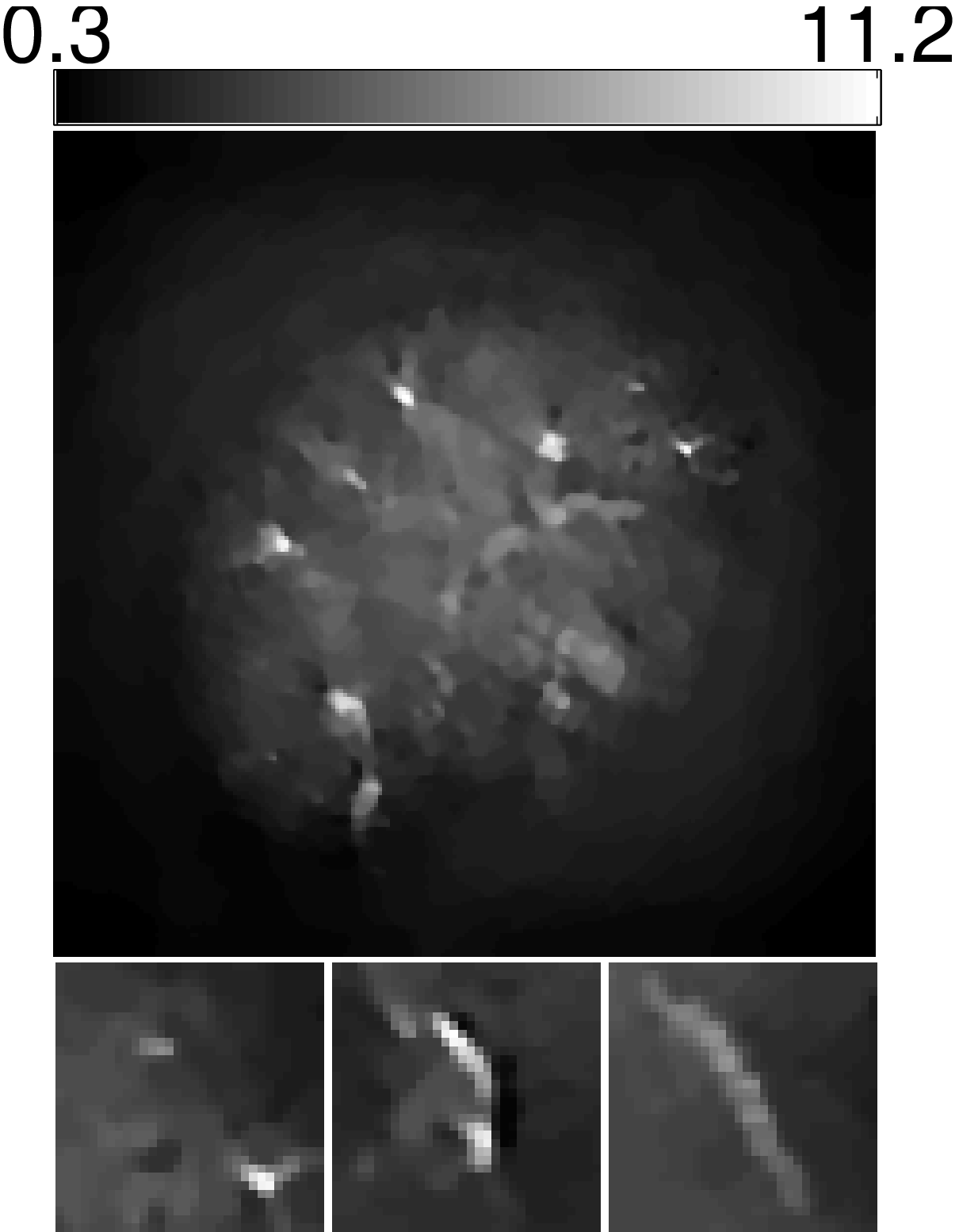}}
\end{center}
\caption{
Slices corresponding to the plane $z=-8.47$-mm through the 3D images of the mouse body 
reconstructed from the $180$-view data by use of the PLS-TV algorithm 
with 
(a) $\beta=0.01$;
and  
(b) $\beta=0.1$. 
The images are of size $29.4\times 29.4$-mm$^2$. 
The three zoomed-in images correspond to the ROIs of the dashed rectangle A,
and the images on the orthogonal planes $x=8.47$-mm (Intersection line is along the arrow B),
and $y=-3.29$-mm (Intersection line is along the arrow C), respectively.
All zoomed-in images are of size $4.34\times 4.34$-mm$^2$.
The ranges of the grayscale windows were determined by the minimum and the maximum values in each image.
\label{fig:mouseTVReg180}
}
\end{figure}

{\bf Qualitative comparison of regularization methods:}
\Fref{fig:mouseV1802D} and \Fref{fig:mouseV1802DY} demonstrate how the 
low-pass filter regularizes the FBP algorithm. 
Similar to the observations from the six-tube phantom, a large value of $f_c$ 
results in high spatial resolution, large signal intensities, 
and high noise level.
For the PLS-TV algorithm, besides the image corresponding to $\beta=0.05$ shown 
in \Fref{fig:mouseV1802D}-(d), 
images corresponding to $\beta=0.01$ and $\beta=0.1$ are displayed in 
\Fref{fig:mouseTVReg180}. 
Though the TV regularization also suppresses the background variance as well 
as the signal intensities when the regularization is enhanced,  
the TV regularization results in minimal sacrifice in spatial resolution. 

\begin{figure}[ht]
\centering
\subfigure[]{\includegraphics[width=7.0cm]{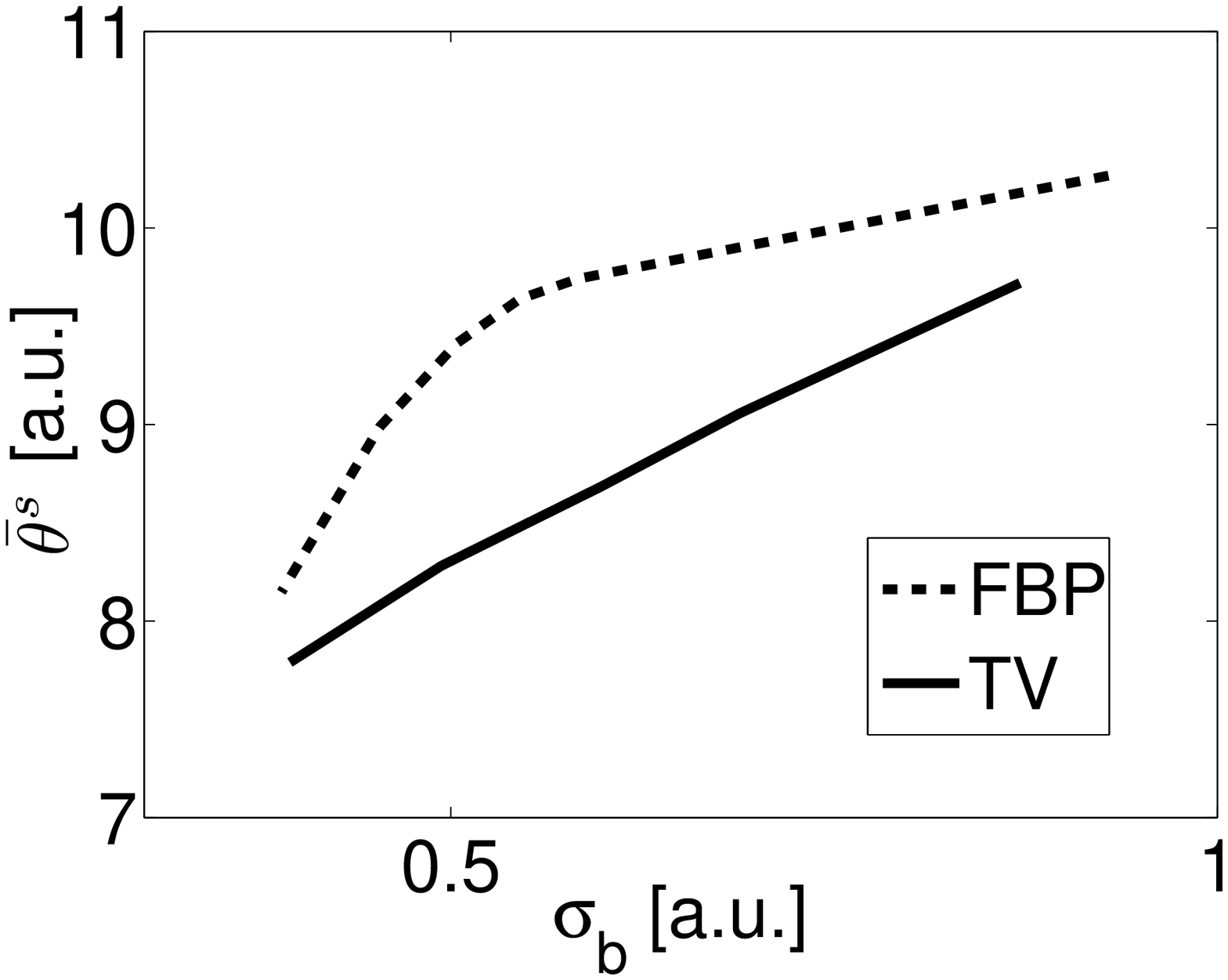}}
\hskip 0.5cm
\subfigure[]{\includegraphics[height=6cm, width=7.0cm]{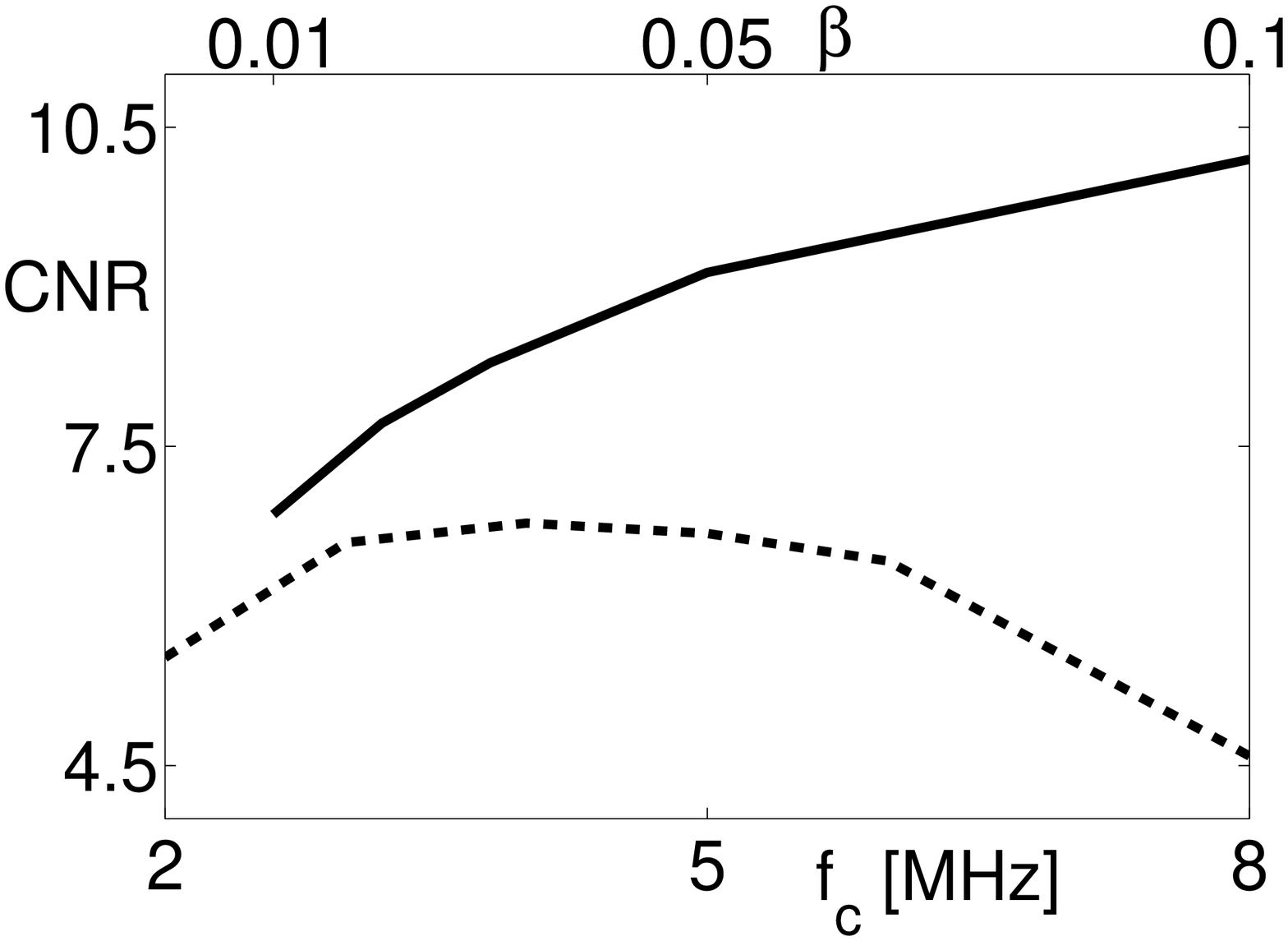}}
\caption{
(a) Signal intensity vs. standard deviation curves 
for the images reconstructed by use of the FBP (dashed line) and the PLS-TV (solid line) algorithms 
from the $180$-view data;
(b) CNR vs. the cutoff frequency curve for the FBP algorithm (dashed line) 
and CNR vs. the regularization parameter $\beta$ curve for the PLS-TV algorithm (solid line)
from the $180$-view data. 
\label{fig:MouseIntStd}
}
\end{figure}

{\bf Trade-off between signal intensity and noise level of 
reconstructed images:}
The s-ROI is defined to be voxels within a blood vessel that  
extends from $z=-9.87$-mm to $z=-7.07$-mm, including $20$ voxels. 
At the plane of $z=-8.74$-mm, the blood vessel is centered at the white 
dashed box D shown in \Fref{fig:mouseTVReg180}-(a). 
The signal intensities are plotted against the image standard deviations in
\Fref{fig:MouseIntStd}-(a).
This plot indicates that the signal intensity in the images reconstructed by 
use of the PLS-TV algorithm is lower than that of the FBP algorithm. 
This reveals that the PLS-TV algorithm can introduce image biases to achieve the 
same level of noise suppression. 
This observation is different than the previous observations from the six-tube
phantom, perhaps because the PLS-TV algorithm somehow promotes discontinuities 
in the diffuse background. 
Nevertheless, the CNR's of the images reconstructed by use of the PLS-TV algorithm 
are higher than those of the FBP algorithms for the regularization parameters
spanning a wide range as shown in \Fref{fig:MouseIntStd}-(b).

\begin{figure}[ht]
\centering
\includegraphics[width=7.0cm]{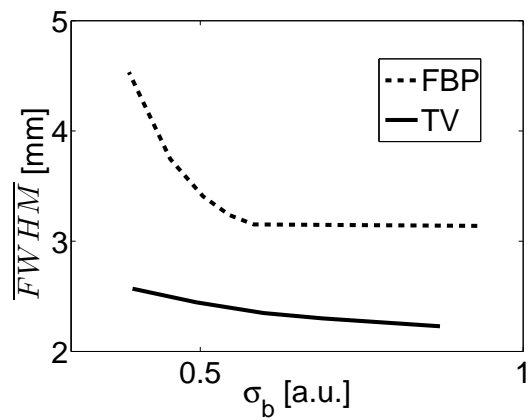}
\caption{
Image resolution vs. standard deviation curves for the images reconstructed by 
use of the FBP and PLS-TV algorithms from the $180$-view data. 
\label{fig:MouseResStd}
}
\end{figure}

{\bf Trade-off between image resolution and noise level of
reconstructed images:}
The plots of resolution against background standard deviation are shown in 
\Fref{fig:MouseResStd}. 
Similar to our observations from the six-tube phantom imaging, 
the spatial resolution of the images reconstructed by use of the PLS-TV 
algorithm is higher than that of the images reconstructed by use of the FBP 
algorithm when the images have the same background standard deviation. 
Also, the curve corresponding to the PLS-TV algorithm is flatter than that of the 
FBP algorithm, confirming that the TV regularization mitigates image noise 
with minimal sacrifice in spatial resolution. 

\begin{figure}[ht]
\begin{center}
\subfigure[]{\includegraphics[width=5.5cm]{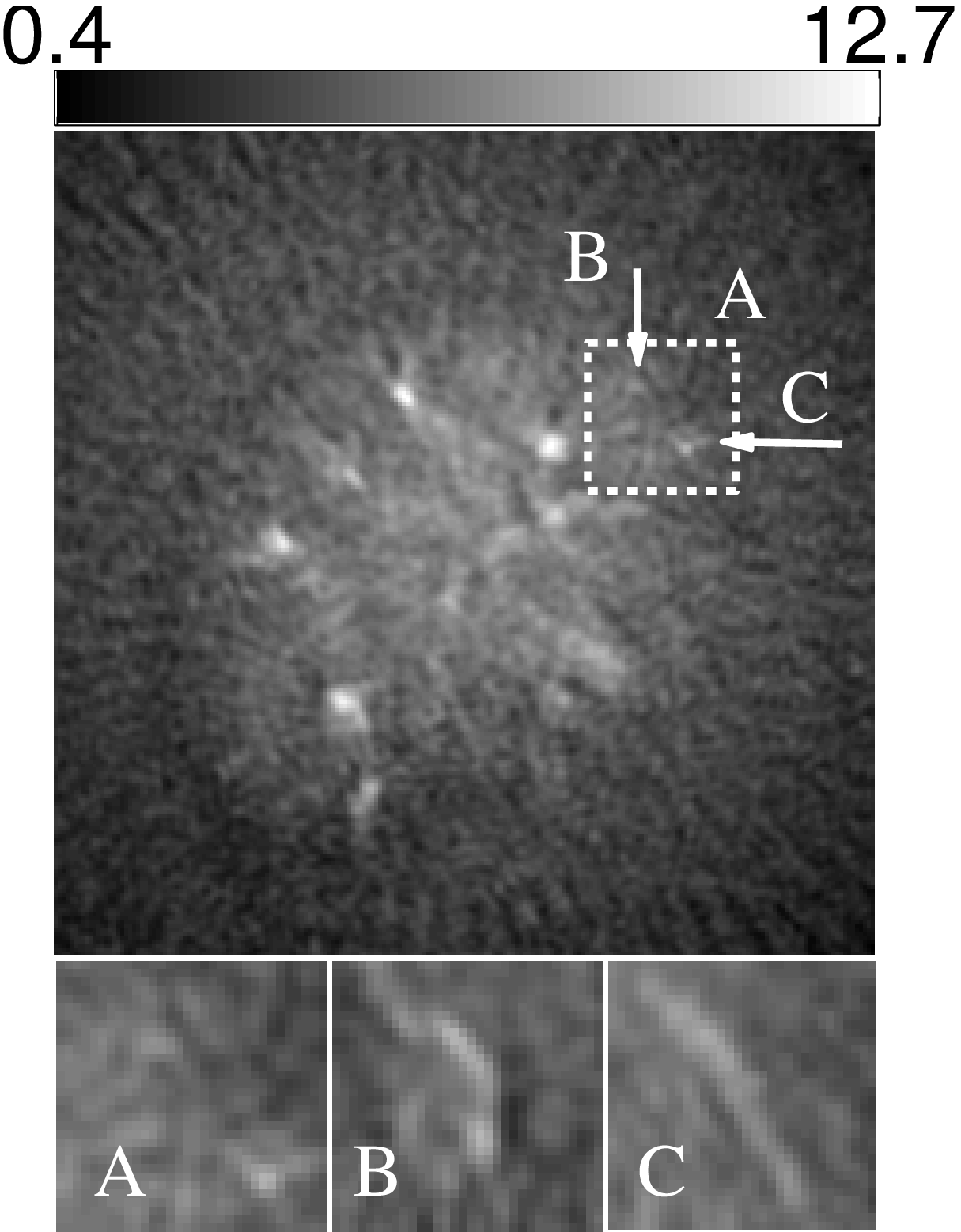}}
\hskip 1cm
\subfigure[]{\includegraphics[width=5.5cm]{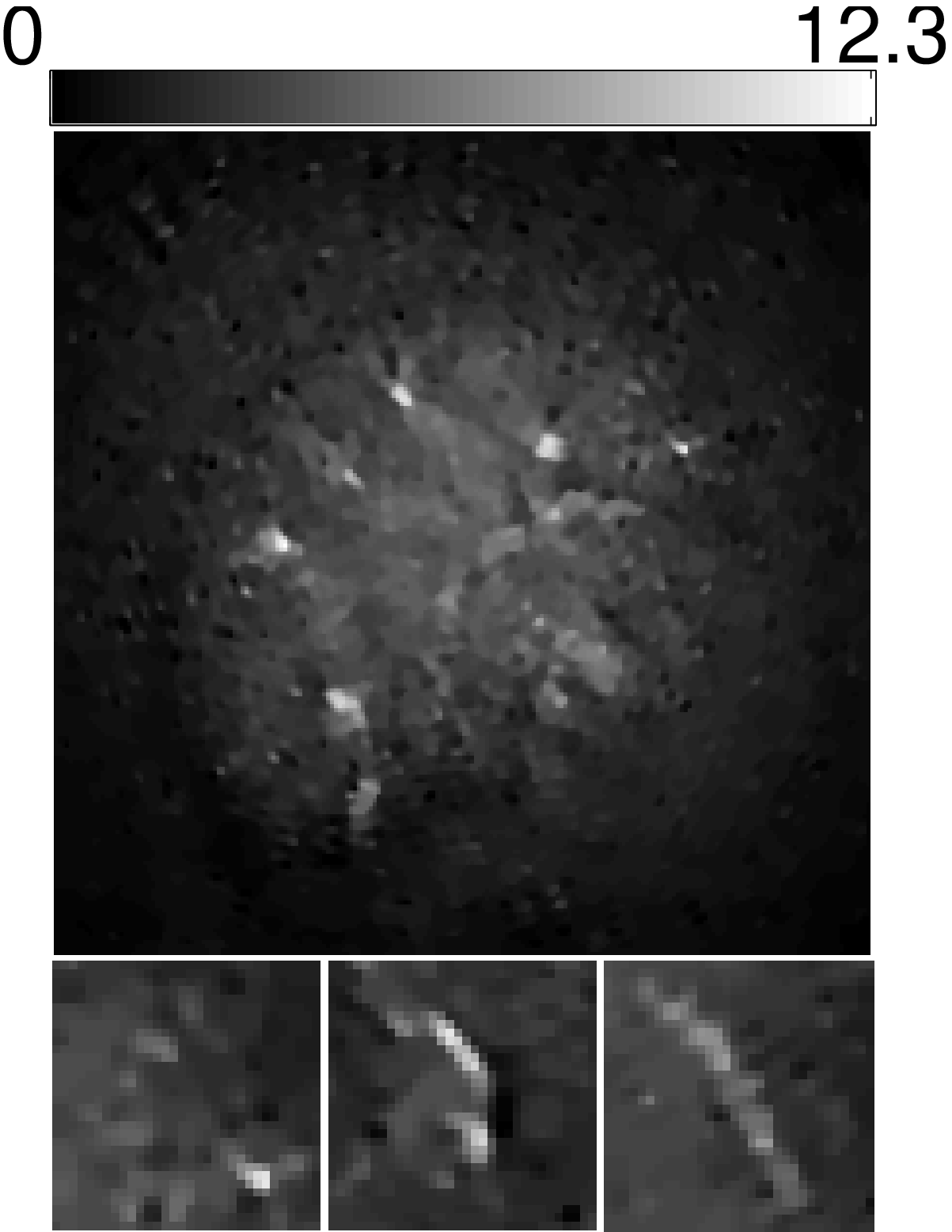}}\\
\subfigure[]{\includegraphics[width=5.5cm]{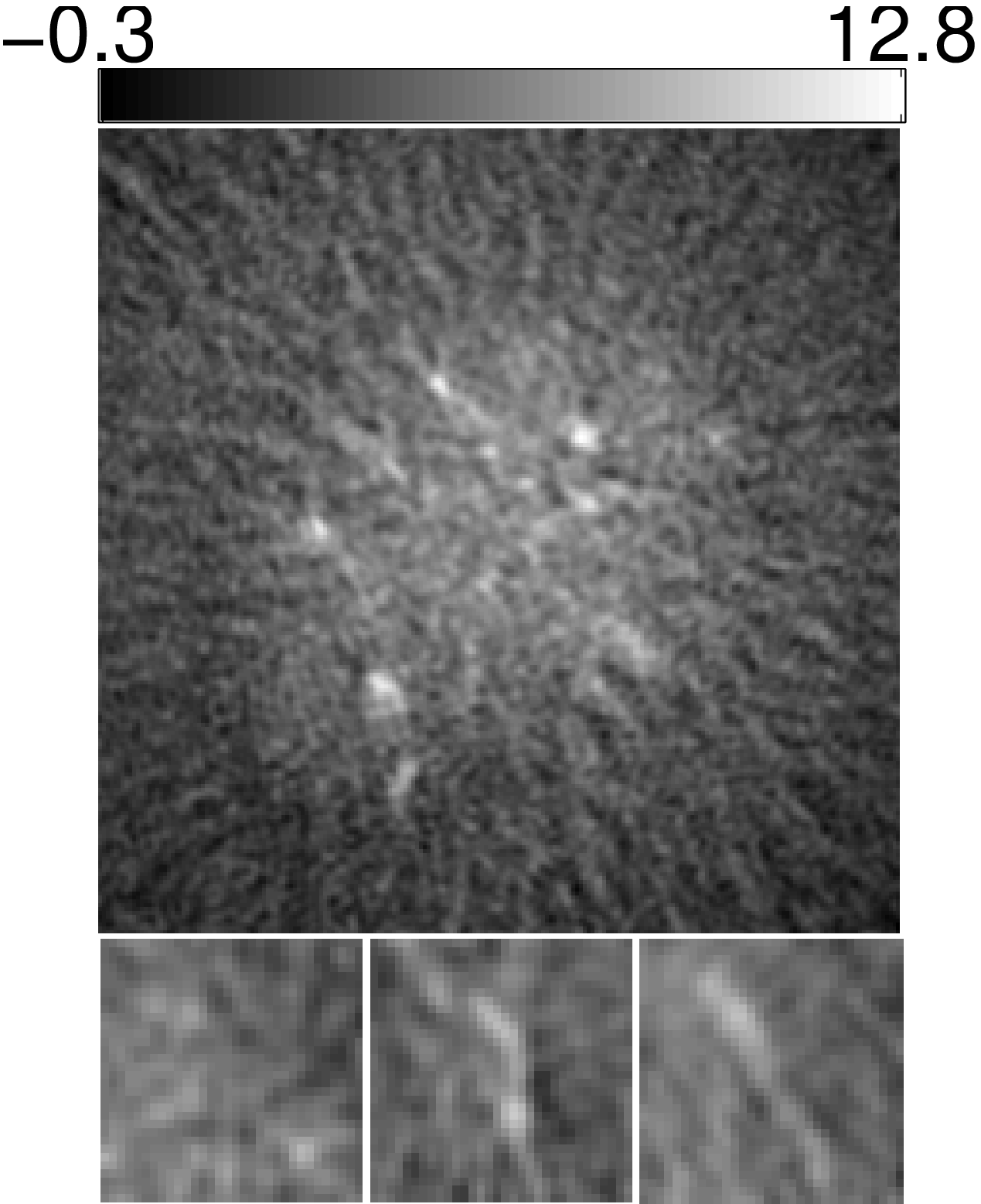}}
\hskip 1cm
\subfigure[]{\includegraphics[width=5.5cm]{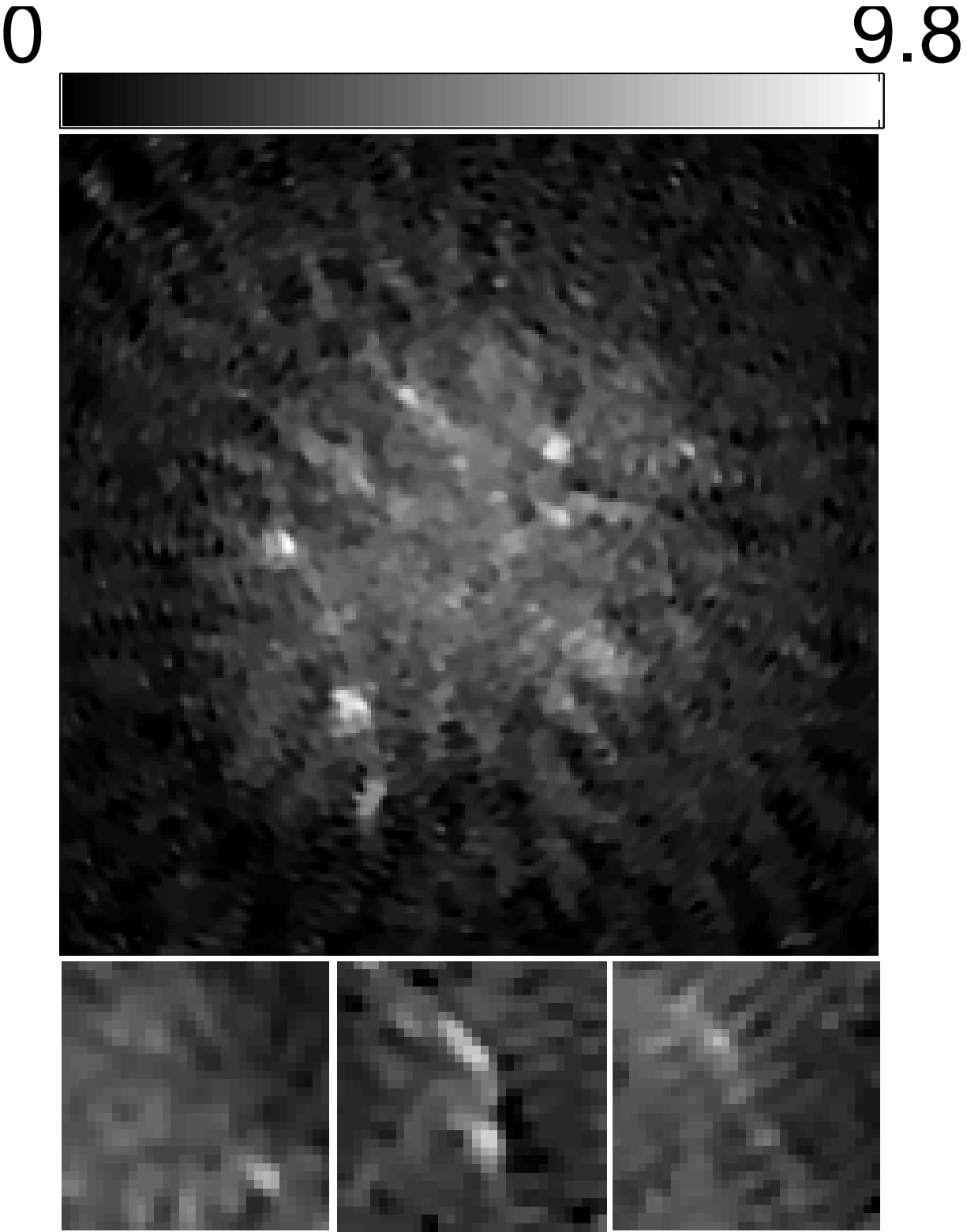}}
\end{center}
\caption{
Slices corresponding to the plane $z=-8.47$-mm through the 3D images of the mouse body
reconstructed from the $90$-view data (top row: a, b) and the $45$-view data 
(bottom row: c, d)  by use of 
(a) the FBP algorithm with $f_c=5$-MHz; 
(b) the PLS-TV algorithm with $\beta=0.03$;
(c) the FBP algorithm with $f_c=5$-MHz;
and
(d) the PLS-TV algorithm with $\beta=0.01$.
The images are of size $29.4\times 29.4$-mm$^2$. 
The three zoomed-in images correspond to the ROIs of the dashed rectangle A,
and the images on the orthogonal planes $x=8.47$-mm (Intersection line is along the arrow B)
and $y=-3.29$-mm (Intersection line is along the arrow C), respectively.
All zoomed-in images are of size $4.34\times 4.34$-mm$^2$.
The ranges of the grayscale windows were determined by the minimum and the maximum valuse of 
each image.  
\label{fig:mouseZ2D}
}
\end{figure}

\begin{figure}[ht]
\begin{center}
\subfigure[]{\includegraphics[width=5.5cm]{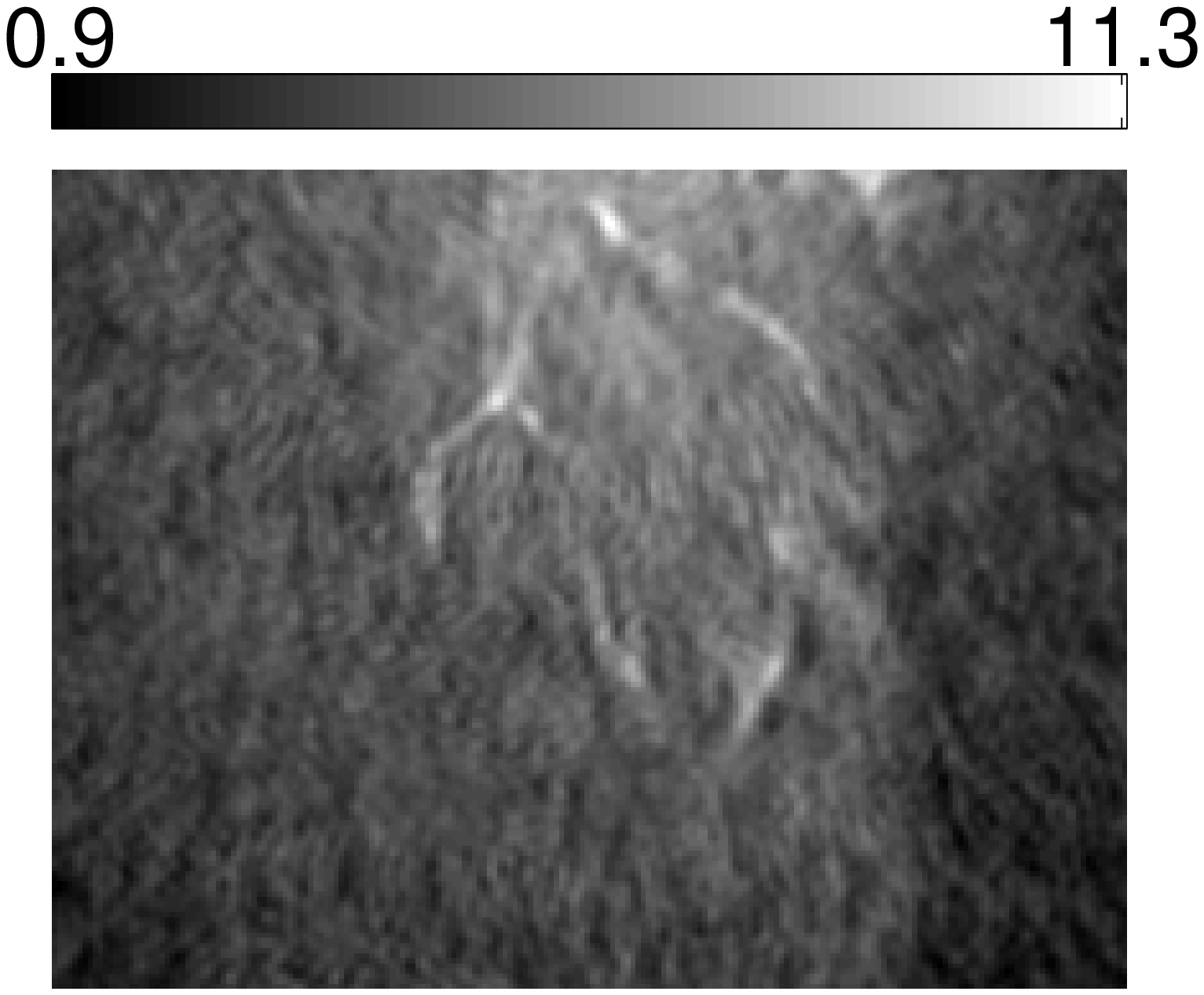}}
\hskip .5cm
\subfigure[]{\includegraphics[width=5.5cm]{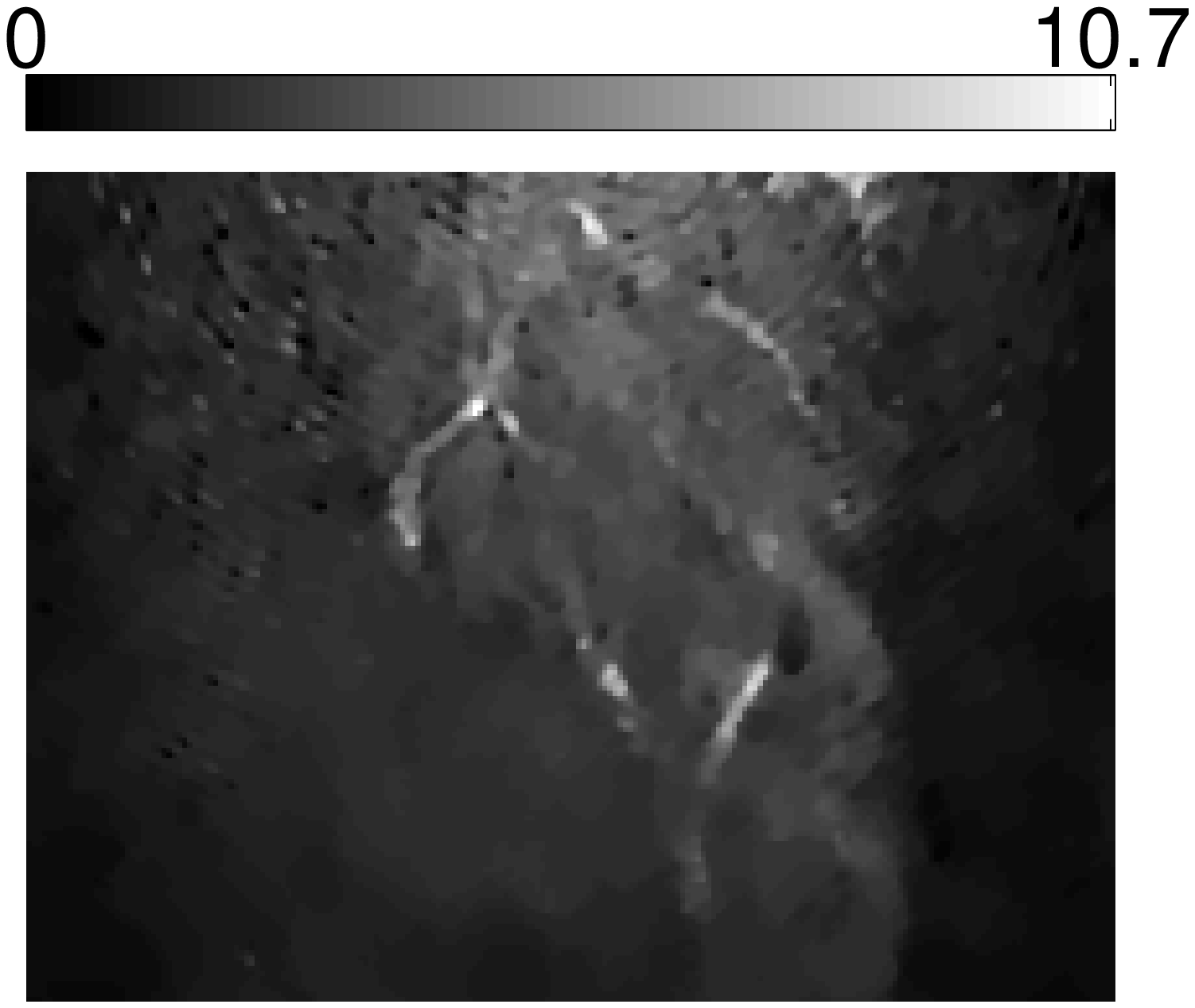}}\\
\subfigure[]{\includegraphics[width=5.5cm]{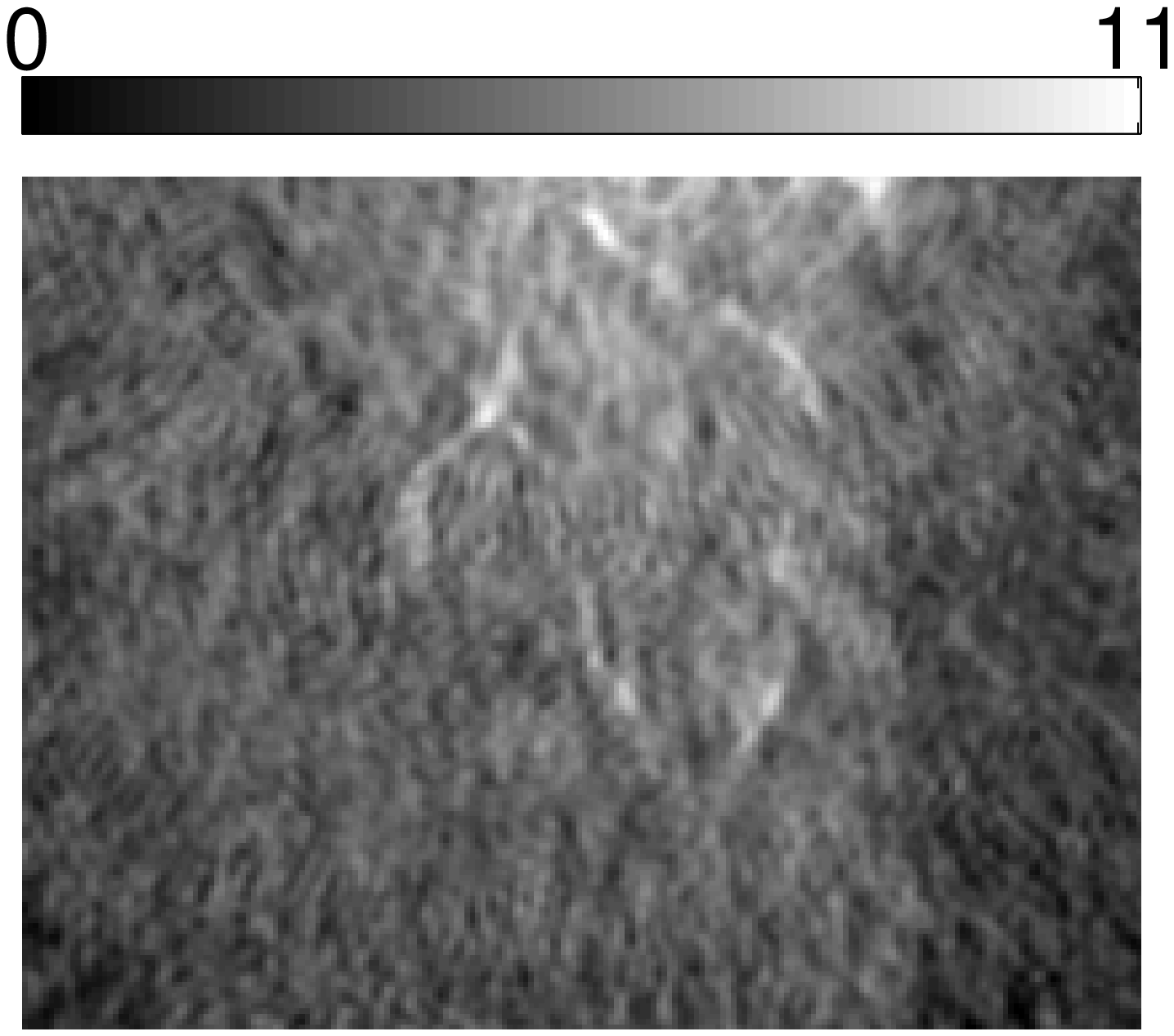}}
\hskip .5cm
\subfigure[]{\includegraphics[width=5.5cm]{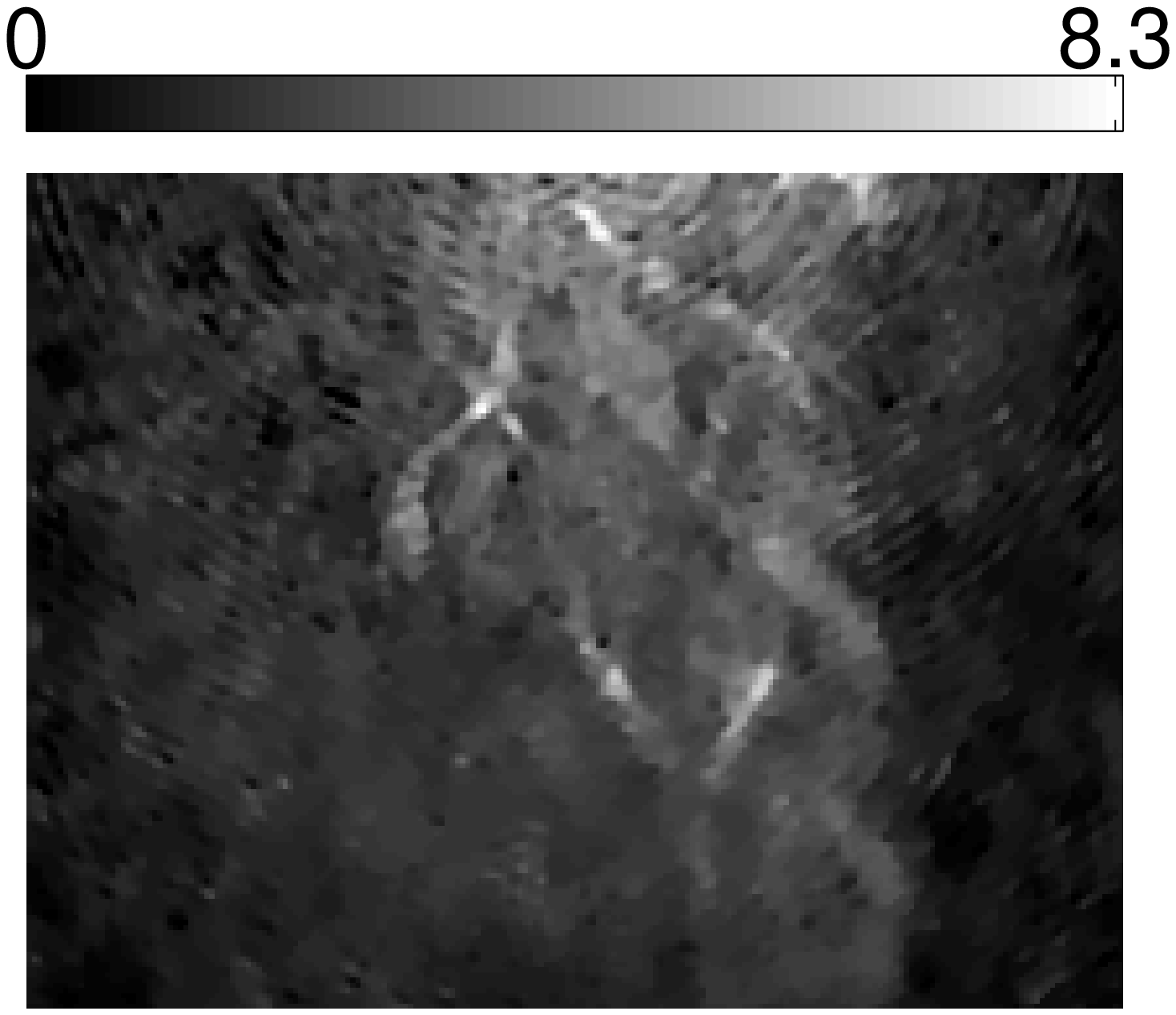}}
\end{center}
\caption{
Slices corresponding to the plane $y=-3.57$-mm through the 3D images of the mouse body 
reconstructed from the $90$-view data (top row: a, b) and the $45$-view data 
(bottom row: c, d) by use of 
(a) the FBP algorithm with $f_c=5$-MHz;
(b) the PLS-TV algorithm with $\beta=0.03$;
(c) the FBP algorithm with $f_c=5$-MHz;
and
(d) the PLS-TV algorithm with $\beta=0.01$.
The images are of size $22.4\times 29.4$-mm$^2$. 
The ranges of the grayscale windows were determined by the minimum and the maximum values
in each image. 
\label{fig:mouseY2D}
}
\end{figure}

{\bf Reconstructed images from sparsely-sampled data sets:}
\Fref{fig:mouseZ2D} and \Fref{fig:mouseY2D} show sectional 
images at different locations. 
Each figure contains subfigures reconstructed by use of the FBP and the PLS-TV 
algorithms from the $90$-view data set and the $45$-view data set. 
The observations are in general consistent with those corresponding to 
densely-sampled data sets; namely the images reconstructed by use of the PLS-TV algorithm
appear to have higher spatial resolution, higher contrast, and cleaner backgrounds. 

\section{Discussion and summary}
\label{sec:summary}
In this study, we investigated two iterative imaging reconstruction algorithms
for 3D OAT: the penalized least-squares (PLS) with a 
quadratic smoothness penalty (PLS-Q) and the PLS with a TV-norm penalty (PLS-TV).
To our knowledge, this was the first systematic investigation of 3D iterative image 
reconstruction for OAT animal imaging. 
Our results demonstrated the feasibility and advantages of 3D iterative image 
reconstruction algorithms for OAT. 
Specifically, the PLS-TV algorithm overall outperforms the FBP algorithm proposed by Finch et al.\ and the conventional iterative image reconstruction algorithm (e.g., PLS-Q) for reconstruction from incomplete data.
Although not reported here, we observed this result to also hold true when other mathematically equivalent FBP algorithms were employed \cite{Xu:2005bp}.

In OAT, the majority of studies of advanced image reconstruction algorithms have been
based on 2D imaging models
\cite{GuoCS:2010,OATTV09:Provost,Ntziachristos:2011,HJiangTV:2011}. 
For a 2D imaging model to be valid in practice, it is necessary to assume the focused 
transducers only receive in-plane acoustic signals. 
The accuracy of this assumption is still under investigation \cite{Ntziachristos:QPAT2012}. 
However, it is interesting to note that none of these studies compared the performances of 
2D analytic reconstruction algorithms with those of the iterative algorithms,
 although the 2D analytic reconstruction
algorithms have been proposed and proved to be mathematically exact 
\cite{Finch:07,Schulze:2DInvPro}. 
In this work, our studies are based on a 3D imaging model that incorporates 
ultrasonic transducer characteristics \cite{TMI:transmodel}, avoiding heuristic 
assumptions regarding the transducer response. 
Although the FBP algorithm neglects the SIR effect, when the region-of-interest is close 
to the center of the measurement sphere, the images reconstructed by use of the FBP algorithm 
from densely-sampled data provide an accurate reference image that permits quantitative evaluation of 
images reconstructed by use of the PLS-Q and PLS-TV algorithms when data are incomplete. 
On the other hand, from densely-sampled data, the images reconstructed by use of different 
algorithms are quantitatively consistent, further validating our 3D imaging model. 

The TV-norm regularization penalty has been intensively investigated within the context of 
 mature imaging modalities including X-ray computed tomography (CT) \cite{XPanIP2009,SidkyQuality:2011}. 
In a study of X-ray CT, the TV-norm regularized iterative reconstruction algorithm 
has been demonstrated to achieve the same image quality as those reconstructed by use of analytic reconstruction 
algorithms, while reducing the amount of data required to one sixth of that the latter requires
\cite{SidkyQuality:2011}. 
However, our images reconstructed from sparsely-sampled data sets by use of the PLS-TV algorithm contain 
clear differences from those reconstructed from densely-sampled data 
by use of the FBP algorithm. 
Moreover, from densely-sampled data, the images reconstructed by use of the PLS-TV algorithm 
also appear to be different from those reconstructed by use of the FBP algorithm. 
Note the streaklike artifacts in the six-tube phantom images reconstructed by use of the 
FBP algorithm in \Fref{fig:SixTubeDenseSample}, 
which remain visible even when the number of tomographic views is increased to $720$. 
These artifacts are likely due to the inconsistencies between the measured data 
and the numerical imaging model.
Such inconsistencies can be caused by unmodeled heterogeneities in the medium \cite{Chao:attenuation,Chao:abberration,Schoonover:11}, 
errors in the assumed transducer response, 
and uncharacterized noise sources \cite{xu:water,xu:edema}. 
These inconsistencies can prevent OAT reconstruction algorithms from working as effectively as 
their counterparts in mature imaging modalities such as X-ray CT that are well-characterized.

There remain several important topics for future studies that may further improve image
quality in 3D OAT.
Such topics include the development and investigation of more accurate imaging models that model the effects of  
acoustic heterogeneities and attenuation. 
Also, in this study, we employed an unweighted least-squares data fidelity term, 
which is equivalent to the maximum likelihood estimator assuming that the randomness in the measured data 
is due to additive Gaussian white noise \cite{Wernickbookchap21}. 
However, additive Gaussian white noise may not be a good approximation in practice 
\cite{Mandelis:PATSNR}. 
Identification of the noise sources and characterization of their second order statistical properties 
will facilitate  iterative reconstruction algorithms that may optimally reduce noise 
levels in the reconstructed images. 
Even though our reconstruction algorithms were implemented using GPUs, the reconstruction time were 
still on the order of hours, which is undesirable for future clinical imaging applications 
of 3D OAT. Therefore there remains a need for the development of accelerated reconstruction
algorithms and their evaluation for specific diagnostic tasks.

\ack     
M.A. Anastasio and K. Wang were supported in part by NIH award EB010049. 
The authors would like to thank Dr. E. Sidky and Dr. X. Pan for stimulating conversations
regarding the use of TV regularization that inspired this work. 
K. Wang would also like to thank Mr. W. Qi for recommending the FISTA algorithm that improved the convergence of the PLS-TV algorithm.
\section*{References}

\bibliographystyle{harvard}   
\bibliography{reflect}   

\clearpage
\newpage

\clearpage

\clearpage

\clearpage

\clearpage

\clearpage

\clearpage

\clearpage

\clearpage

\clearpage

\clearpage

\clearpage

\clearpage

\clearpage

\clearpage

\clearpage

\clearpage

\end{document}